\documentclass{amsart}
\usepackage{graphicx}
\graphicspath{ {./images/} }

\usepackage{amsmath,amssymb,amsthm}
\allowdisplaybreaks[4]

\usepackage{tikz}
\usetikzlibrary{arrows.meta,calc}

\usepackage{float}

\usepackage{xcolor}

\usepackage[backend=bibtex,url=false]{biblatex}
\usepackage{hyperref}
\usepackage{cleveref}

\usepackage{algpseudocode,algorithm}

\newcommand*\Let[2]{\State #1 $\gets$ #2}
\algrenewcommand\algorithmicrequire{\textbf{Input:}}

\newtheorem{theorem}{Theorem}[section]

\newtheorem{lemma}[theorem]{Lemma}
\newtheorem{corollary}[theorem]{Corollary}
\newtheorem{proposition}[theorem]{Proposition}

\newtheorem{remark}[theorem]{Remark}
\newtheorem{definition}[theorem]{Definition}
\newtheorem{example}[theorem]{Example}

\def\N{\mathbb{N}}
\def\Z{\mathbb{Z}}
\def\Q{\mathbb{Q}}
\def\R{\mathbb{R}}

\def\cA{\mathcal{A}}
\def\cW{\mathcal{W}}
\newcommand{\Cp}{C_p}

\DeclareMathOperator{\Diff}{Diff}
\DeclareMathOperator{\pal}{pal}

\title{Boundary of the Lagrange and Markov spectra near Hall's ray}

\author[H. Erazo]{Harold Erazo}
\address{Harold Erazo: SUSTech International Center for Mathematics, Shenzhen, Guangdong, People's Republic of China}
\email{harold.eraz@gmail.com}

\author[C. Matheus]{Carlos Matheus}
\address{Carlos Matheus: CMLS, \'Ecole Polytechnique, CNRS (UMR 7640),
91128, Palaiseau, France}
\email{matheus.cmss@gmail.com}

\author[C. G. Moreira]{Carlos Gustavo Moreira}
\address{Carlos Gustavo Moreira: SUSTech International Center for Mathematics, Shenzhen, Guangdong, People's Republic of China and IMPA, Rio de Janeiro, Brazil}
\email{gugu@impa.br}
\thanks{The third author was partially supported by CNPq and FAPERJ}

\date{\today}



\begin{document}

\begin{abstract}
The Lagrange and Markov spectra, denoted by $L$ and $M$, respectively, are closed subsets of the real line with intricate fractal geometry related to Diophantine approximations. In 1975, Freiman proved that both spectra contain the ray $[c_F,\infty)$, where $c_F=4.527829566\dots$, and that this ray is optimal, as both spectra avoid a tiny interval $(\nu_F,c_F)$. Using renormalization techniques and the thickness of Cantor sets, we prove that the boundaries $\partial L$ and $\partial M$ have positive Hausdorff dimensions arbitrarily close to $\nu_F$. 

Furthermore, in connection with Berstein's conjecture, which predicts that $[4.1,4.52]\subset L$, we exhibit a gap in $L$ near $4.08$, the first known gap in the region $L\cap[3.945,4.52]$. We also prove that the boundaries of both spectra have positive Hausdorff dimension near $4.5213$ and $4.08$.
\end{abstract}

\maketitle

\section{Introduction}

Arithmetic sums and stable intersections of Cantor sets are ubiquitous objects in Mathematics: for example, they played a key role in the seminal works of Palis, Takens, Newhouse, Yoccoz and the third author 
(cf. \cite{PalisTakensAnnals1987}, \cite{NewhouseIHES1979}, \cite{GuguThesis}, \cite{MoreiraYoccoz2010}) about homoclinic bifurcations and, as it was first pointed out to the third author by J. Palis more than 30 years ago, also in the study of Diophantine questions in Number Theory - namely, Palis, who was the Ph.D. advisor of the third author, mentioned the paper \cite{Ha47} by M. Hall, who contains the first non-trivial condition for an arithmetic sum of Cantor sets to contain a nontrivial interval, in the context of continued fractions with bounded partial quotients. In the present paper, we will follow this tradition by describing how arithmetic sums of Cantor sets explain a subtle behavior of the boundaries of the Lagrange and Markov spectra. 
 
More concretely, recall that the Lagrange and Markov spectra ${L}$ and ${M}$ are subsets of the real line related to
Diophantine approximation problems. They are classically defined by
$${L}:=\left\{k(\alpha):=\limsup\limits_{\substack{p,q\to\infty\\p,q\in\mathbb{Z}}} \frac{1}{|q(q\alpha-p)|}<\infty: \alpha\notin\mathbb{Q}\right\}$$
and 
$${M}:=\left\{\frac{1}{\inf\limits_{\substack{(x,y)\in\mathbb{Z}^2\\(x,y)\neq(0,0)}}|q(x,y)|}<\infty: q\in \mathcal{Q} \right\}$$ 
where $\mathcal{Q}=\{q(x,y)=ax^2+bxy+cy^2; a, b, c\in {\mathbb R}, \ b^2-4ac=1\}$. In other words, the Lagrange and Markov spectra are subsets of the real line consisting of the best constants of Diophantine approximations of certain irrational numbers and indefinite binary quadratic forms.

It is well known that $L\subset M\subset\R^+$ and that both spectra exhibit intricate fractal geometry. One of the first striking results concerning these sets, proved by Markov \cite{M:formes_quadratiques1, M:formes_quadratiques2} around 1880, is that they coincide below $3$ and that $L\cap(0,3)=M\cap(0,3)$ is discrete, with $3$ as its only accumulation point. On the other hand, M. Hall \cite{Ha47} proved in 1947 that both spectra contain a half-line $[c,\infty)$ for some $c>5$. In 1975, after extensive computations, Freiman \cite{Fr75} determined the optimal ray $[c_F,\infty)\subset L\subset M$ contained in both spectra and identified a gap $(\nu_F,c_F)$ of $M$, where
$c_F=4.527829566\dots$ and $\nu_F=4.527829538\dots$ are explicit algebraic numbers, with both endpoints belonging to $L$.

Between these discrete and continuous regimes, the structure of the spectra is considerably more complicated. The third author proved in \cite{geometricproperties} that the Hausdorff dimensions of their truncations coincide and define a continuous function
\begin{equation*}
d(t)=\dim_H(L\cap(-\infty,t))
=\dim_H(M\cap(-\infty,t))
\end{equation*}
for all $t\in\R$. Moreover, it was pointed out in \cite{geometricproperties} that
\begin{equation}
t_1=\min\{t\in\R:d(t)=1\}=3.334384\dots
\end{equation}
may mark another transition between two distinct regimes (the first six decimal digits of $t_1$ were determined in \cite{MMPV}). Indeed, the sets $L\cap[3,t_1]$ and $M\cap[3,t_1]$ are compact, totally disconnected, and have zero Lebesgue measure and in \cite{geometricproperties} it was asked whether the interiors of $L\cap[t_1,\nu_F]$ and $M\cap[t_1,\nu_F]$ accumulate at $t_1$. On the other hand, let $L^\prime$ denote the set of accumulation points of $L$. Since $L^\prime$ was shown in \cite{geometricproperties} to be a perfect set, it may exhibit local features reminiscent of a \emph{Cantorval} in the region $L^\prime\cap[t_1,\nu_F]$. Here, a Cantorval is understood as a perfect set with dense interior in which every gap is accumulated on at least one side by infinitely many gaps and nontrivial intervals. We use this term in a less restrictive sense than the definition given in \cite{MendesOliveira}. We remark that this expectation might be false in general, if the answer for the first question of \cite{geometricproperties} is negative, namely, if there is a drop in the local Hausdorff dimension of $L^\prime$.

This latter possibility of a drop (nonmonotonicity) of the local dimension, however, does hold for the set of accumulation points of the Markov spectrum throughout the regime $M^{\prime}\cap[t_1,\nu_F]$. Indeed, it was shown in \cite{Fractalgeometryofcomplement} that there exist $j_0<j_1$, with $j_0=3.709\dots$, such that
\begin{equation}
M\cap(j_0,j_1)=(M\setminus L)\cap(j_0,j_1)
\end{equation}
and
\begin{equation}
0.5<\dim_H(M\cap(j_0,j_1))<0.6.
\end{equation}

In 1973, Berstein \cite{BerConj}, then a student of Freiman, computationally investigated the region of the Markov spectrum between $4.1$ and $5.1$. More specifically, on page 77 of \cite{BerConj}, he explicitly identified two gaps of $M$ near $4.521$. On the same page, he conjectured that $[4.1,4.52]\subset L$ and that $[4.522,\infty)\subset L$. The latter conjecture was subsequently disproved by Freiman's work, whereas the former remains open. Although the gaps he identified justify the choice of $4.52$ as the right endpoint of the interval in his first conjecture, he provided no explanation for his choice of $4.1$ as the left endpoint.

\subsection{Boundary points between $4.521$ and $\nu_F$}

Recently, the last two authors, together with L. Jeffreys, proved in \cite{Newgaps} that $M$ has an infinite sequence of maximal gaps $(\nu_n,\mu_n)$ with $\nu_n\to\nu_F$, i.e., with $\nu_n,\mu_n\in M$. Their approach relies on renormalization techniques and local thickness of regular Cantor sets, ideas that appeared in the paper corresponding to the third author's Ph.D. thesis \cite[Theorem IV.1]{GuguThesis}. This result is a fundamental geometrical tool to extend the classical results by Palis and Takens (\cite{PalisTakensAnnals1987}) in order to get persistently positive density of hyperbolicity after the unfolding of certain open sets of homoclinic bifurcations associated to horseshoes with Hausdorff dimensions larger than one (the main result of the forementioned paper by Palis and Takens gives full density of hyperbolicity after the unfolding of homoclinic bifurcations in surfaces associated to horseshoes with Hausdorff dimensions smaller than one). In the present work, we develop this approach further to strengthen the results of \cite{Newgaps}. Given any subset $X\subset\R$, we denote by $\partial X$ its topological boundary with respect to the usual topology.

\begin{theorem}\label{thm:uncountable_nu_F}
    For any $\varepsilon>0$, the set $\partial L\cap\partial M\cap[\nu_F-\varepsilon,\nu_F]$ has positive Hausdorff dimension. In fact, it has Hausdorff dimension at least
    \begin{equation*}
        \frac{\log 2}{2\log\left(2\sqrt{243542} + 987\right)} = 0.04567\dots
    \end{equation*}
\end{theorem}

Let us recall briefly the Perron dynamical characterization of the spectra. Given a bi-infinite sequence $\underline{a}=(a_n)_{n\in\Z}\in\{1,2,\dots\}^\Z$, we define the following sum of continued fractions
\begin{equation*}
    \lambda_n(\underline{a}):=[a_n;a_{n+1},a_{n+2},\dots]+[0;a_{n-1},a_{n-2},\dots].
\end{equation*}
The Markov and Lagrange value of $\underline{a}$ are respectively
\begin{equation*}
    \ell(\underline{a})=\limsup_{n\to\infty}\lambda_n(\underline{a}), \qquad m(\underline{a})=\sup_{n\in\Z}\lambda_n(\underline{a}).
\end{equation*}

The Markov and Lagrange spectra can be characterized as follows
\begin{equation*}
    L = \{\ell(\underline{a})<\infty:\underline{a}\in(\N^*)^\Z\}, \qquad M = \{m(\underline{a})<\infty:\underline{a}\in(\N^*)^\Z\}.
\end{equation*}
Since our sequences $(a_n)_{n\in\Z}\in\{1,2,3,4\}^{\Z}$ will have bounded digits , we will write the bi-infinite sequence as $\dots a_{-2}a_{-1}a_0^*a_1a_2\dots$. We will use $\overline{abc}=abcabc\dots$ to denote periodicity.

From \cite[Theorem 3]{BerConj}, we know that Berstein proved that there is a gap of $M$ in the interval $(\nu_*,\mu_*)$ with $$\nu_* = \lambda_{0}(\overline{111313121313}331213^{*}13123\overline{313121313111}) = 4.521715440484506773...$$ and $$\mu_* = \lambda_{0}(\overline{111313121313}331213131231213^{*}1312\overline{313121313111}) = 4.521731667575794970... .$$
and there is a gap of $M$ in the interval $(\nu_{**},\mu_{**})$ with $$\nu_{**} = \lambda_{0}(\overline{1313^{*}1212}) = 4.521039005773764...$$ and $$\mu_{**} = \lambda_{0}(\overline{111313121313}11344431213^{*}1312\overline{313121313111}) = 4.5212652356170763....$$

Computer investigations led us to conjecture that the second gap identified by Berstein, namely $(4.52103\dots,4.52126\dots)$, is in fact the largest gap of $M\cap[4.521,\nu_F]$. Moreover, our numerical experiments suggest that the region $M\cap[4.52125,4.52135]$ has a particularly intricate structure, whereas the region $M\cap[4.52135,\nu_F]$ may contain large intervals. 

Motivated by these observations, we show that our approach can also be used to establish the uncountability of the boundaries of both spectra in other regions whose structure is not yet well understood. Since our method requires an irrationality condition at specific periodic points of the Gauss map, the above gaps found by Berstein do not work well because the periods that appear in the sequences associated with the borders of the gaps are exactly the same. Fortunately, through computer experiments we found abundant gaps nearby that satisfy our conditions. We will instead show that $M$ has gap from
\begin{equation}\label{eq:nu_4521}
    \nu=\lambda_0\left(\overline{131312131311}34431213^*1312312\overline{21313121}\right)=4.5212920725004\dots
\end{equation}
to
\begin{equation*}
    \mu=\lambda_0(\overline{131312131311}344431213^*13123123\overline{313121313111}) =4.521292947013\dots
\end{equation*}

Since the periods of $\nu$ are different, and in fact, the corresponding multipliers of the Gauss map are multiplicatively independent, together with a thickness condition satisfied by some suitable dynamically defined Cantor sets, we will show that 

\begin{theorem}\label{thm:uncountable_4521}
    For every $\varepsilon>0$, the set $\partial L\cap\partial M\cap[\nu-\varepsilon,\nu]$ has positive Hausdorff dimension, where $\nu$ is defined on \eqref{eq:nu_4521}. 
\end{theorem}

\subsection{Boundary and gaps near $4.1$}

Recall that Berstein \cite[Page 77]{BerConj} conjectured that $[4.1,4.52]\subset L$. However, at present, there is no known explicit gap in the region $M\cap[3.945,4.1]$. Through our numerical explorations we have been able to identify several tiny previously unknown gaps in this region, for instance there are gaps of $M$ containing 3.9453, 3.95155, 3.963358, 3.9703953, 4.0601, 4.080495, respectively. As the pictures suggest, finding these gaps is not easy because they are probably surrounded by very large intervals in $M$. We will provide explicit computations for the gap containing 4.080495 and show that even in these regions where is difficult to find gaps, the boundary of both spectra is also uncountable.

Let 
\begin{equation*}
    \nu' = m\left(\overline{3133}11213^*213212\overline{1333}\right) = 4.08049296860135092\dots
\end{equation*}
and 
\begin{equation}\label{eq:mu_408}
    \mu' = m\left(\overline{132311}3311213^*21321\overline{3313}\right) = 4.08049996235899758\dots
\end{equation}

\begin{theorem}\label{thm:uncountable_408}
    For every $\varepsilon>0$, the set $\partial L\cap\partial M\cap[\mu',\mu'+\varepsilon]$ has positive Hausdorff dimension, where $\mu'$ is defined on \eqref{eq:mu_408}. 
\end{theorem}

\subsection{A new portion of $M_3\setminus L_3$}

Another topic investigated by Freiman and Berstein was the study of the difference set $M\setminus L$. In 1968, Freiman \cite{Fr68} gave the first example of a point belonging to the Markov spectrum but not to the Lagrange spectrum. Since then, several works have been devoted to understanding the structure of $M\setminus L$. For instance, the last two authors \cite{Fractalgeometryofcomplement} proved that $\dim_H(M\setminus L)<1$, which implies that $L$ and $M$ have the same interior. A significant open problem is to determine $\sup(M\setminus L)$. The largest known examples, found in \cite{LGH}, lie near $3.942$, but it is expected that $M\setminus L$ contains points between $4.521$ and $4.527$.

Although we have been unable to solve this problem, we have found examples in a more restrictive setting. Consider the Lagrange and Markov spectra restricted to continued fractions with digits in $\{1,\dots,N\}$:
\begin{equation*}
    L_N := \{\ell(\underline{a}):\underline{a}\in\{1,\dots,N\}^\Z\}, \quad
    M_N := \{m(\underline{a}):\underline{a}\in\{1,\dots,N\}^\Z\}.
\end{equation*}

The main result of \cite{DynamicalCharacterization} establishes that, for every $N\neq 3$,
\begin{equation*}
    L_N = L\cap\left(-\infty,\sqrt{N^2+4N}\right], \quad
    M_N = M\cap\left(-\infty,\sqrt{N^2+4N}\right].
\end{equation*}
This equality fails for $N=3$, since $M_3$ has substantially larger gaps than $M$. Indeed, it was shown in \cite[Page 850]{DynamicalCharacterization} that
\begin{equation*}
    \left(\frac{4\sqrt{462}}{19},\frac{154+\sqrt{21}}{35}\right)
    =(4.52509\dots,4.5309307\dots)
\end{equation*}
is a maximal gap of $M_3$. Nevertheless, this interval contains many points of both $L$ and $M$, including points belonging to the beginning of Hall's ray. We prove the following result.
\begin{theorem}\label{thm:M_3_minus_L_3}
    $(M_3\setminus L_3)\cap(4.5249,4.525)\neq\varnothing$.
\end{theorem}

The existence of such examples reinforces our expectation that $M\setminus L$ contains points between $4.521$ and $4.527$. However, the examples provided by \Cref{thm:M_3_minus_L_3} do not appear to help resolve this problem, since the entire interval $[4.5249,4.525]$ might be contained in $L$ (and hence in $M$).

\subsection{Organization of the paper}
In \Cref{sec:general_criteria} we will provide general criteria for exhibiting positive Hausdorff dimension boundary for the arithmetic sums of regular Cantor set. In \Cref{sec:uncountable_nu_F} we will prove \Cref{thm:uncountable_nu_F} based on the previous results of \cite{Newgaps}. In \Cref{sec:uncountable_near_45213} we will prove \Cref{thm:uncountable_4521}. In \Cref{sec:gap_near_408} we will prove \Cref{thm:uncountable_408}. In \Cref{sec:M_3_minus_L_3} we will state the lemmas necessary to proof \Cref{thm:M_3_minus_L_3} without providing details. Finally, in \Cref{sec:algorithm}, we will explain the algorithm used to obtain the new plots for the Markov spectrum which vastly influenced this paper.

\section{A general criteria for uncountability of the boundary of arithmetic sums of regular Cantor sets} \label{sec:general_criteria}

\subsection{Pointwise thickness}
A dynamical (or regular) Cantor set $K$ is a compact set defined by a Markov partition $\mathcal{P}=\{I_1,I_2,...,I_s\}$, a topologically mixing $C^{r}$-expanding map $\psi:I_1\cup I_2\cup ...\cup I_s\rightarrow I$, where $I$ is the convex hull of the disjoint union of intervals $I_1\cup I_2\cup ...\cup I_s$ and

\begin{enumerate}
	\item[1)] for any $j=1,2,...,s$ the $\psi(I_j)$ is the convex hull of an union $\cup_{t=1}^{R(j)}I_{m_t}$
	
	\item[2)] the set $K$ is given by $$K=\bigcap_{n\ge 0}\psi^{-n}(I_1\cup I_2\cup... \cup I_s).$$

    \item[3)] for each $j=1,2,\dots,s$ and $N$ sufficiently big, $\psi^N(K\cap I_j)=K$ (this means that $\psi|_{K}$ is topologically mixing).
\end{enumerate}

We say that a finite word $\underline{a}=(a_1,\dots,a_n)\in\{1,\dots,s\}^n$ is admissible if $I_{a_{i+1}}\subset\psi(I_{a_i})$ for $i=1,\dots,n-1$. Given an admissible finite word $\underline{a}=(a_1,\dots,a_n)$, we will denote the corresponding interval of the construction of $K$ by
\begin{equation*}
    I(\underline{a})=I_{a_1}\cap\psi^{-1}(I_{a_2})\dots\cap\psi^{-(n-1)}(I_{a_n}).
\end{equation*}

In this paper we will always consider regular Cantor sets with regularity at least $r=1+\varepsilon^\prime$. A remarkable fact about $C^{1+\varepsilon^\prime}$, $\varepsilon^\prime>0$, regular Cantor sets is the \emph{bounded distortion property}, which in rough terms says that relative distances on a microscopic scale are at most distorted by a constant factor with respect to relative distances corresponding on a larger scale, and increasingly less distorted if we decrease the size of the major scale. Given $\delta>0$, there exists $C(\delta)$ decreasing on $\delta$, such that if $x,y,z\in K$ are such that $|\psi^i(y)-\psi^i(x)|\leq \delta$ and $|\psi^i(z)-\psi^i(x)|\leq \delta$ for $i=0,1,\dots,n$, then
\begin{equation}\label{eq:bounded_distortion}
	e^{-C(\delta)}\frac{|z-x|}{|y-x|}\leq\frac{|\psi^n(x)-\psi^n(y)|}{|\psi^n(y)-\psi^n(x)|}\leq e^{C(\delta)}\frac{|z-x|}{|y-x|}.
\end{equation}

\begin{remark}\label{rem:smaller_copies}
    In general, if $K$ is a $C^{r}$ dynamically defined Cantor set defined by an expanding map $\psi:I_1\cup\dots\cup I_s\to I$ and $f$ is a $C^r$ map defined on some neighborhood of $I$ with $|f^\prime|$ uniformly bounded by below in $I$, then we can exhibit $f(K)$ as a dynamically defined Cantor set by conjugating the dynamics, namely, first we consider a refined Markov partition given by the intervals $I(\underline{a})$ with all $\underline{a}\in\{1,\dots,s\}^N$ admissible and by defining $\tilde{I}(\underline{a})=f(I(\underline{a}))$ and $\tilde{\psi}=f\circ\psi^N\circ f^{-1}$. Since $|\tilde{\psi}^\prime(f(x))|=|(\psi^N)^\prime(x)|\frac{|f^\prime(\psi^N(x))|}{|f^\prime(x)|}$, for some large $N$ the map $\tilde{\psi}$ is expanding. 
\end{remark}

\begin{example}\label{example}
    Let us describe the kind of Cantor sets that appear in the study of the Markov and Lagrange spectra. Let $\mathcal{A}=\{1,\dots,A\}$ be an alphabet, $\vec{r}=(r_1,\dots,r_k)$ be a vector of lengths and let $F$ be a set of forbidden words
\begin{equation*}
    F=F(\vec{r})=\{w_1^{(1)}\dots w_{r_1}^{(1)}\in\mathcal{A}^{r_1},\dots,w_1^{(k)}\dots w_{r_k}^{(k)}\in\mathcal{A}^{r_k}\}.
\end{equation*}
The \emph{Gauss-Cantor} set defined by $F$ is the set $X_F\subset[0,1]$ of continued fractions $[0;a_1,a_2,\dots]$ such that $(a_1,a_2,\dots)\in\mathcal{A}^\N$ does not contain any of the subwords of $F$. 

We can exhibit such a Cantor set as a dynamically defined Cantor set. First, let us consider $R=\max_{1\leq i\leq k}r_i$. Given a finite word $w=(w_1,\dots,w_n)\in\{1,\dots,A\}^n$, we define the following interval
\begin{equation*}
    I(w)=\textrm{convex hull of } \{[0;w_1,\dots,w_n,x]\in X_F\}.
\end{equation*}
The domain of the Markov partition of $X_F$ is the set of intervals $I(u)$ where $u\in\mathcal{A}^{R-1}$ does not contain any of the words of $F$ and $\psi=g$ restricted to the union of those intervals. 
\end{example}

\begin{definition}\label{def:pointwise_thickness}
A gap $U$ of a Cantor set $K$, is a connected component of $\R\setminus K$. Given a bounded gap $U$ of $K$, its left and right bridges, denoted by $L_U$ and $R_U$, are the maximal closed intervals sharing an endpoint with the boundary of $U$ and that contain no gap of length greater or equal than $|U|$.

\begin{center}
\begin{tikzpicture}[scale=1.5]
    \draw (-3.5,0) -- (-3,0);
    \node[red] at (-3,0) {(};
    \node[red] at (-1,0) {)};

    \draw (-1,0) -- (0,0) node[midway, above] {$L_U$};

    \node[red] at (0,0) {(};
    \node[red] at (1,0) {)};
    \node[above,yshift=0.5em] at (0.5,0) {$U$};

    \draw (1,0) -- (2,0) node[midway, above] {$R_U$};
    \node[red] at (2,0) {(};
    \node[red] at (4,0) {)};

    \draw (4,0) -- (4.5,0);
\end{tikzpicture}
\end{center}

Given $p\in\partial U^\prime$ a boundary point of a gap of $K$, define the \emph{pointwise thickness} of $K$ at $p$ by
\begin{align*}
    \tau_p(K,p)&=
\liminf_{\substack{\min U\to p,\\ p=\min L_U}}
\frac{|L_U|}{|U|}, && \text{if } p =\max\partial U^\prime, \\
\tau_p(K,p)&=
\liminf_{\substack{\max U\to p,\\ p=\max R_U}}
\frac{|R_U|}{|U|},  && \text{if } p =\min\partial U^\prime.
\end{align*}
where the infimum is taken over bounded gaps $U$ of $K$.
\end{definition}

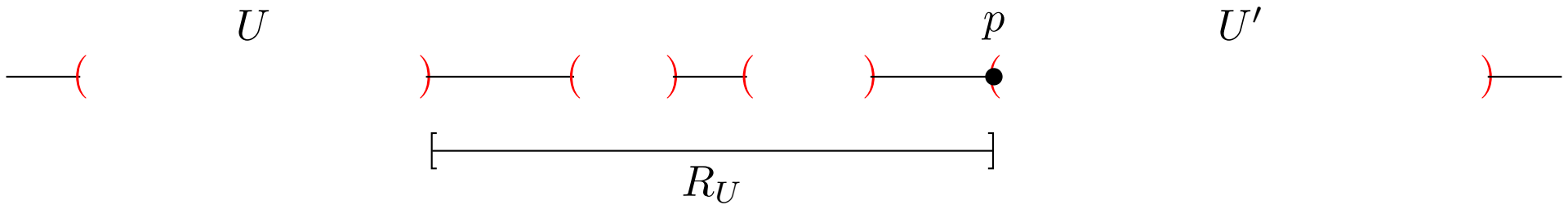
\begin{figure}
\centering
\begin{tikzpicture}[scale=2]
\draw (-3.5,0) -- (-3.2,0);

\node[red] at (-3.2,0) { ( };
\node[red] at (-1.8,0) { ) };
\node[above,yshift=0.5em] at (-2.5,0) {$U$};

\draw (-1.8,0) -- (-1.2,0);
\node[red] at (-1.2,0) { ( };
\node[red] at (-0.8,0) { ) };

\draw (-0.8,0) -- (-0.5,0);
\node[red] at (-0.5,0) { ( };
\node[red] at (0,0) { ) };

\draw (0,0) -- (0.5,0);
\node[red] at (0.5,0) { ( };
\node[red] at (2.5,0) { ) };
\node[above,yshift=0.5em] at (1.5,0) {$U^\prime$};

\draw (2.5,0) -- (2.8,0);

\coordinate (A) at (0.5,0);
\fill (A) circle (1pt);
\node[above,yshift=0.5em] at (A) {$p$};

\draw[{Bracket[width=3mm]}-{Bracket[width=3mm]}] (-1.78,-0.3) -- (0.5,-0.3) node[midway, below] {$R_U$};

\end{tikzpicture}
\caption{The pointwise thickness of a point $p=\min\partial U^\prime$ is the limit of ratios $|R_U|/|U|$.}
\end{figure}

Observe that for a dynamically defined Cantor set $K$ defined by some expanding map $\psi$ with Markov partition $\{I_1,\dots,I_s\}$, any boundary point of a gap $U$ of $K$ is preperiodic for $\psi$. To see this, notice that the Markov condition implies that if we call $E=\bigcup_{j=1}^s \partial I_j$ then $\psi(E)\subset E$, so the points of $E$ are preperiodic. On the other hand, for gaps such that $U\subset I_1\cup\dots\cup I_s$, there is some iterate $\ell_p$ such that $\psi^{\ell_p}(U)$ is a primary gap, i.e., a bounded connected component of $\R\setminus(I_1\cup\dots\cup I_s)$, so in particular for those boundary points $p\in\partial U$ we obtain $\psi^{\ell_p}(p)\in E$.

From the previous paragraph we see that there is a finite set of periodic points $\{q_1,\dots,q_m\}\subset K$ with $\psi^{k_i}(q_i)=q_i$, $k_i\geq 1$ such that for any boundary point $p\in\partial U$ of a gap of the Cantor set there is some $\ell_p$ and $1\leq i\leq m$ such that $\psi^{\ell_p}(p)=q_i$. The next lemma shows that $\tau_p(K,p)=\tau_p(K,q)$.

\begin{lemma}
Let $p\in\partial U^\prime$ be a boundary point of a gap $U^\prime$ of $K$. Let $f:[a,b]\to\R$ be a $C^1$ diffeomorphism defined on an interval such that $a<p<b$. Then we have that $\tau_p(f(K\cap[a,b]),f(p))=\tau_p(K,p)$. 
\end{lemma}

\begin{proof}
    By continuity of $f^{\prime}$ near $p$, given $\tilde{\varepsilon}>0$ we can chose $\tilde{\delta}>0$ such that $e^{-\tilde{\varepsilon}}\leq|f^{\prime}(x)|/|f^{\prime}(y)|\leq e^{\tilde{\varepsilon}}$ for all $|x-p|<\tilde{\delta}$ and $|y-p|<\tilde{\delta}$. 
    
    Given $\varepsilon>0$, let $U$ be a sufficiently small gap of $K$ with $p\in\partial B_U$ for $B_U\in\{L_U,R_U\}$ with $|U|+|B_U|<\tilde{\delta}$ and $\frac{|B_U|}{|U|}\leq\tau_p(K,p)+\varepsilon$. By the mean value theorem, we have that the gap $f(U)$ of $f(K\cap[a,b])$ satisfies \[\frac{|f(B_U)|}{|f(U)|}\leq e^{\tilde{\varepsilon}}\frac{|B_U|}{|U|}\leq e^{\tilde{\varepsilon}}(\tau_p(K,p)+\varepsilon).\] 

    We do not necessarily have that $f(B_U)$ is equal to $B_{f(U)}$, however, if this is not the case, then $f(B_U)$ contains a gap $\widehat{U}$ of $f(K\cap[a,b])$ with $|\widehat{U}|\geq|f(U)|$, $q_i\in B_{\widehat{U}}$ and clearly $|B_{\widehat{U}}|<|f(B_U)|$. In particular replacing $f(U)$ by $\widehat{U}$ and $f(B_U)$ by $B_{\widehat{U}}$ decreases even more the quotient $|f(B_U)|/|f(U)|$. The previous argument shows that $\tau_p(f(K\cap[a,b]),f(p))\leq \tau_p(K,p)$. Since we can apply the same argument to $f^{-1}$ we obtain the equality.
\end{proof}

In certain situations, for example when we know the limit geometries of our Cantor set, it is possible to explicitly compute the pointwise thickness. The reason is that the limit geometry map linearizes the dynamics at periodic points, as we will see below. We will follow the exposition of \cite{Gazette} and remark that the difference between the limit geometries defined on \cite{MoreiraYoccoz2010} is only a matter of scale, since in \cite{MoreiraYoccoz2010} the limit geometries expand to $[0,1]$ while in the exposition below they expand to a respective interval of the initial Markov partition.

By definition, if $(a,b)$ is admissible, then $\psi$ is a diffeomorphism between a sub-inteval $I(a,b)\subset I(a)$ and $I(b)$, so we denote $f_{a,b}=\left(\psi|_{I(a,b)}\right)^{-1}$ as this inverse branch of $\psi$. In general, if $\underline{a}=(a_0,\dots,a_n)$ is an admissible word, then we denote $f_{\underline{a}}=f_{a_0,a_1}\circ\dots\circ f_{a_{n-1},a_n}$ to denote the inverse branch of $\psi^n$ that sends $I(a_n)$ to a subinterval $I(\underline{a})$ of $I(a_0)$. 

Given an admissible $\underline{\theta}=(\theta_{-n})_{n\leq0}\in\{1,\dots,s\}^{\Z^{-}}$, let $\underline{\theta}^{(n)}=(\theta_{-n},\dots,\theta_0)$ and let $B(\underline{\theta}^{(n)}):I(\underline{\theta}^{(n)})\to I(\theta_0)$ be the unique onto affine transformation such that $k^{\underline{\theta}^{(n)}}:=B(\underline{\theta}^{(n)})\circ f_{\underline{\theta}^{(n)}}$ is a orientation preserving diffeomorphim. Note that in this way $k^{\underline{\theta}^{(n)}}$ is defined on $I(\theta_0)$ taking values in $I(\theta_0)$, that is, $k^{\underline{\theta}^{(n)}}: I(\theta_0)\rightarrow I(\theta_0).$ Moreover, The image of the set \(K \cap I(\theta_0)\) is an enlarged copy of the set \(K \cap I(\underline{\theta}^{(n)})\). 

\begin{proposition}{\cite{Sullivan}}
Let $K$ be a dynamically defined Cantor set of class $C^r$, $r\in(1,\infty)$. For any admissible $\underline{\theta}=(\theta_{-n})_{n\leq0}\in\{1,\dots,s\}^{\Z^{-}}$, there is a map $k^{\underline{\theta}}\in\Diff^r_{+}(I(\theta_0))$ such that the sequence $k^{\underline{\theta}^{(n)}}$ converges to $k^{\underline{\theta}}$ in $\Diff^{r^\prime}_{+}(I(\theta_0))$ for any $r^\prime<r$, uniformly in $\underline{\theta}$.
\end{proposition}

We define the \emph{limit geometry} of $K$ associated to $\underline{\theta}$ the set
\begin{equation*}
K^{\underline{\theta}}=k^{\underline{\theta}}\left(K\cap I(\theta_0)\right)
\end{equation*}
where $k^{\underline{\theta}}=\lim_{n\to\infty}k^{\underline{\theta}^{(n)}}$.

Given some admissible finite word $\underline{a}=(\theta_0,a_1,\dots,a_n)$ , let \[\underline{\theta}\ast\underline{a}=(\dots,\theta_{-\ell},\dots,\theta_0,a_1,\dots,a_n).\] Let $F^{\underline{\theta}}_{\underline{a}}$ be the affine map that sends $I(a_n)$ onto $k^{\underline{\theta}}(I(\underline{a}))$  with the same orientation as $f_{\underline{a}}$. In particular 
\begin{equation}\label{eq:limit_geometry_relation}
    k^{\underline{\theta}}\circ f_{\underline{a}}=F^{\underline{\theta}}_{\underline{a}}\circ k^{\underline{\theta}\ast\underline{a}}.
\end{equation}

Let $p\in\partial U^\prime$ be a boundary point of a gap $U^\prime$ of $K$. Let $\ell_p\geq 0$ be such that $\psi^{\ell_p}(p)=q$. Then, we have that $\tau_p(K,p)=\tau_p(K,q)$. Now suppose that $q$ is associated to an admissible periodic sequence $\underline{\theta}$ with period $(\theta_0,\dots,\theta_{k-1})$.  In particular, letting $\underline{a}=(\theta_0,\dots,\theta_{k-1},\theta_0)$ we have that $f_{\underline{a}}(q)=q$ and $\underline{\theta}\ast\underline{a}=\underline{\theta}$.  Now let us recall, as in \Cref{rem:smaller_copies}, that we can realize the limit geometry $K^{\underline{\theta}}$ by conjugating the dynamics. More precisely by \eqref{eq:limit_geometry_relation} we have that $k^{\underline{\theta}}\circ f_{\underline{a}}\circ (k^{\underline{\theta}})^{-1}=F^{\underline{\theta}}_{\underline{a}}$. Since $k^{\underline{\theta}}(q)$ is the fixed point of the linear map $F^{\underline{\theta}}_{\underline{a}}$,  we can compute $\tau_p(K,q)$ by computing finitely many quotients of the form $|B_{\widehat{U}}|/|\widehat{U}|$ where $\widehat{U}=k^{\underline{\theta}}(U)$ and $\widehat{B}_{\widehat{U}}$ is the bridge of $K^{\underline{\theta}}$ with $k^{\underline{\theta}}(q)\in\partial \widehat{B}_{\widehat{U}}$. More precisely, we only have to consider those gaps $\widehat{U}$ of $K^{\underline{\theta}}$ contained in $I(\theta_0)$ but not in $k^{\underline{\theta}}(I(\underline{a}))$ and such that $k^{\underline{\theta}}(q)\in\partial \widehat{B}_{\widehat{U}}$.

\begin{example}\label{example_1}
From \cite{Gazette}, we know that if $C_N=\{[0;a_1,a_2,\dots]:a_i\in\{1,\dots,N\},\forall i\geq 1\}$, then the limit geometries are the M\"{o}bius transformations
\begin{equation*}
    k^{\underline{\theta}}(x)=\mu_{\underline{\theta}}+\frac{1+\beta_{\underline{\theta}}\nu_{\underline{\theta}}}{1+\beta_{\underline{\theta}}x}(x-\mu_{\underline{\theta}}), \quad x\in[\mu_{\underline{\theta}},\nu_{\underline{\theta}}]=I(\theta_0)
\end{equation*}
where $\underline{\theta}\in\{1,\dots,N\}^{\Z^{-}}$, $\mu_{\underline{\theta}}=[0;\theta_0,\overline{1,N}]$, $\nu_{\underline{\theta}}=[0;\theta_0,\overline{N,1}]$ and $\beta_{\underline{\theta}}=[0;\theta_{-1},\theta_{-2},\dots]$. Indeed, we have that $f_{\underline{\theta}^{(n)}}(x)=[0;\theta_{-n},\dots,\theta_{-1}+x]$, so if we take $n$ odd then the previous map is increasing we can take 
\begin{equation*}
    B(\underline{\theta}^{(n)})(x)=\frac{|I(\theta_0)|}{|I(\underline{\theta}^{(n)})|}(x-\mu_{\underline{\theta}})+\mu_{\underline{\theta}}.
\end{equation*}
In particular, if we denote $p_n/q_n = [0;\theta_{-n},\dots,\theta_{-1}]$, $p_{n-1}/q_{n-1} = [0;\theta_{-n},\dots,\theta_{-2}]$ and $\beta_{n+1}=q_{n-1}/q_n$, then we have by the theory of continued fractions that
\begin{equation*}
    f_{\underline{\theta}^{(n)}}(x)-f_{\underline{\theta}^{(n)}}(\mu_{\underline{\theta}}) = \frac{\frac{1}{\mu_{\underline{\theta}}}-\frac{1}{x}}{q_n^2(\frac{1}{x}+\beta_{n+1})(\frac{1}{\mu_{\underline{\theta}}}+\beta_{n+1})}=\frac{x-\mu_{\underline{\theta}}}{q_n^2(1+\beta_{n+1}x)(1+\beta_{n+1}\mu_{\underline{\theta}})}.
\end{equation*}
Therefore
\begin{align*}
    B(\underline{\theta}^{(n)})\circ f_{\underline{\theta}^{(n)}}(x)&=\frac{|I(\theta_0)|}{|I(\underline{\theta}^{(n)})|}(f_{\underline{\theta}^{(n)}}(x)-\mu_{\underline{\theta}})+\mu_{\underline{\theta}} \\
    &=\frac{\nu_{\underline{\theta}}-\mu_{\underline{\theta}}}{f_{\underline{\theta}^{(n)}}(\nu_{\underline{\theta}})-f_{\underline{\theta}^{(n)}}(\mu_{\underline{\theta}})}(f_{\underline{\theta}^{(n)}}(x)-\mu_{\underline{\theta}})+\mu_{\underline{\theta}} \\
    &=\frac{1+\beta_{n+1}\nu_{\underline{\theta}}}{1+\beta_{n+1}x}(x-\mu_{\underline{\theta}})+\mu_{\underline{\theta}}
\end{align*}
Since $\beta_{n+1}=[0;\theta_{-1},\dots,\theta_{-n}]$, taking $n\to\infty$ yields the claimed formula.

For $\max C_N = [0;\overline{1,N}]$ we have that $f_{\underline{a}}(x)=[0;1,N+x]$ is the relevant contraction map and $F^{\underline{\theta}}_{\underline{a}}$ is the affine map that sends $I(1)=[[0;1,\overline{1,N}],[0;\overline{1,N}]]$ onto $k^{\underline{\theta}}(I(\underline{a}))$ where $I(1,N,1)=I(\underline{a})=f_{\underline{a}}(I(1))=[[0;1,N,1,\overline{1,N}],[0;\overline{1,N}]]$. 

We will compute the pointwise thickness for $N=3$ and $\max C_N$. We need to consider those gaps $\widehat{U}$ inside of $k^{\underline{\theta}}(I(1))$ but not in $k^{\underline{\theta}}(I(1,3,1))$ such that the right border of their right bridge is equal to $\max C_3 = [0;\overline{1,3}]$. There are only three such gaps which are given explicitly by
\begin{align*}
    \widehat{U}_1 &= \left(k^{\underline{\theta}}([0;1,1,\overline{1,3}]),k^{\underline{\theta}}([0;1,2,\overline{3,1}])\right), \\
    \widehat{U}_2 &= \left(k^{\underline{\theta}}([0;1,2,\overline{1,3}]),k^{\underline{\theta}}([0;1,3,\overline{3,1}])\right), \\
    \widehat{U}_3 &= \left(k^{\underline{\theta}}([0;1,3,2,\overline{3,1}]),k^{\underline{\theta}}([0;1,3,1,\overline{1,3}])\right).
\end{align*}

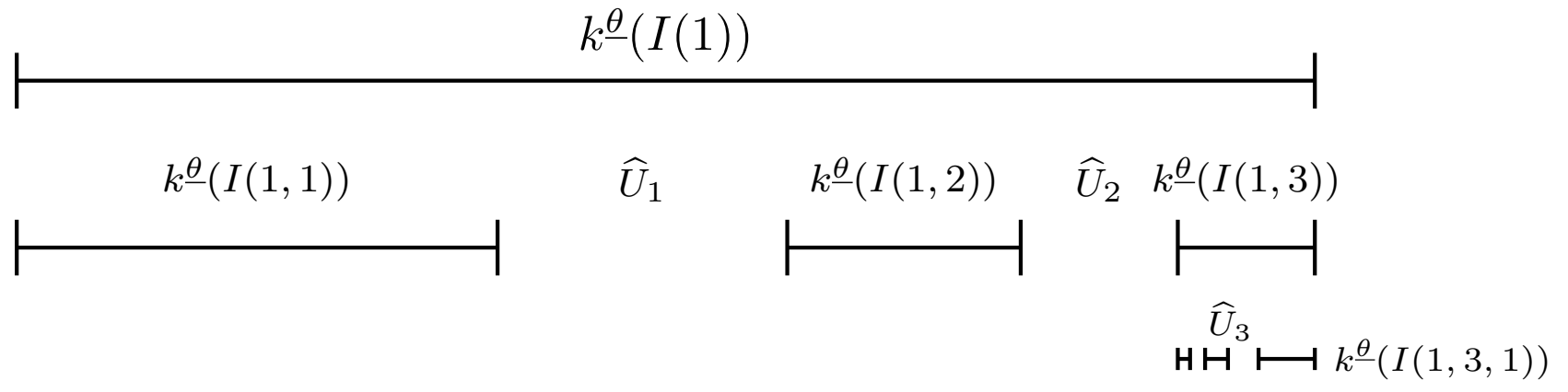
\begin{figure}
    \centering
	\begin{tikzpicture}[line width=0.8pt,scale=40]

    \draw (0.5582,0) -- (0.7912,0);
    \draw (0.5582,0.005) -- (0.5582,-0.005);
    \draw (0.7912,0.005) -- (0.7912,-0.005);

    \draw (0.5582,-0.03) -- (0.6445,-0.03);
    \draw (0.5582,-0.025) -- (0.5582,-0.035);
    \draw (0.6445,-0.025) -- (0.6445,-0.035);

    \draw (0.6965,-0.03) -- (0.7384,-0.03);
    \draw (0.6965,-0.025) -- (0.6965,-0.035);
    \draw (0.7384,-0.025) -- (0.7384,-0.035);

    \draw (0.7666,-0.03) -- (0.7912,-0.03);
    \draw (0.7666,-0.025) -- (0.7666,-0.035);
    \draw (0.7912,-0.025) -- (0.7912,-0.035);

    \draw (0.7666,-0.05) -- (0.7688,-0.05);
    \draw (0.7666,-0.048) -- (0.7666,-0.052);
    \draw (0.7688,-0.048) -- (0.7688,-0.052);

    \draw (0.7715,-0.05) -- (0.7756,-0.05);
    \draw (0.7715,-0.048) -- (0.7715,-0.052);
    \draw (0.7756,-0.048) -- (0.7756,-0.052);

    \draw (0.7811,-0.05) -- (0.7912,-0.05);
    \draw (0.7811,-0.048) -- (0.7811,-0.052);
    \draw (0.7912,-0.048) -- (0.7912,-0.052);

    \node[above] at (0.6747,0) {$k^{\underline{\theta}}(I(1))$};
    \node[above, font = \footnotesize] at (0.60135,-0.025) {$k^{\underline{\theta}}(I(1,1))$};
    \node[above, font = \footnotesize] at (0.71745,-0.025) {$k^{\underline{\theta}}(I(1,2))$};
    \node[above, font = \footnotesize] at (0.7789,-0.025) {$k^{\underline{\theta}}(I(1,3))$};

    \node[above, font = \footnotesize] at (0.6705,-0.025) {$\widehat{U}_1$};
    \node[above, font = \footnotesize] at (0.7525,-0.025) {$\widehat{U}_2$};
    \node[above, font = \scriptsize] at (0.77605,-0.05) {$\widehat{U}_3$};

    \node[right, font = \scriptsize] at (0.7912,-0.05) {$k^{\underline{\theta}}(I(1,3,1))$};

    \end{tikzpicture}
    \caption{Limit geometry of $C_3$ at $\max C_3$ and gaps whose right bridge contains the maximum.}
\end{figure}

To justify this, observe that the lengths of these gaps are
\begin{align*}
    |\widehat{U}_1|&=k^{\underline{\theta}}([0;1,2,\overline{3,1}])-k^{\underline{\theta}}([0;1,1,\overline{1,3}])=\frac{-157\sqrt{21} + 723}{68}=0.05199\dots \\
    |\widehat{U}_2|&=k^{\underline{\theta}}([0;1,3,\overline{3,1}])-k^{\underline{\theta}}([0;1,2,\overline{1,3}])=\frac{-2027\sqrt{21} + 9303}{500}=0.0282\dots \\
    |\widehat{U}_3|&=k^{\underline{\theta}}([0;1,3,1,\overline{1,3}])-k^{\underline{\theta}}([0;1,3,2,\overline{3,1}])=\frac{-139\sqrt{21} + 637}{4}=0.00549\dots
\end{align*}

Since
\begin{align*}
    |k^{\underline{\theta}}(I(1,1))|&=\frac{107\sqrt{21} - 483}{85} = 0.0863\dots \\
    |k^{\underline{\theta}}(I(1,1,1))|&=\frac{6883\sqrt{21} - 31437}{3145}=0.033\dots \\
    |k^{\underline{\theta}}(I(1,1,2))|&=\frac{-79}{15}\sqrt{21} + \frac{483}{20} = 0.0151\dots 
\end{align*}
the interval $k^{\underline{\theta}}(I(1,1))$ contains no gaps larger than $\widehat{U}_1$. Since
\begin{align*}
    |k^{\underline{\theta}}(I(1,2))|&=\frac{98\sqrt{21} - 447}{50} = 0.0418\dots \\
    |k^{\underline{\theta}}(I(1,2,1))|&=\frac{1133\sqrt{21} - 5187}{300}=0.0168\dots 
\end{align*}
the interval $k^{\underline{\theta}}(I(1,2))$ contains no gaps larger than $\widehat{U}_1$ and $\widehat{U}_2$. Since
\begin{align*}
    |k^{\underline{\theta}}(I(1,3))|&=\frac{443\sqrt{21} - 2027}{125} = 0.0246\dots \\
    |k^{\underline{\theta}}(I(1,3,1))|&=\frac{43\sqrt{21} - 197}{5}=0.0101\dots \\
    |k^{\underline{\theta}}(I(1,3,2))|&=\frac{661\sqrt{21} - 3029}{20}=0.00412\dots \\
    |k^{\underline{\theta}}(I(1,3,3))|&=\frac{-6857\sqrt{21} + 31423}{125}=0.0022\dots
\end{align*}
the interval $k^{\underline{\theta}}(I(1,3))$ contains no gaps larger than $\widehat{U}_1$ and $\widehat{U}_2$ and moreover the larger gap contained on it is $\widehat{U}_3$.

Therefore to compute the pointwise thickness we only need to compute the following three quotients
\begin{align*}
    \frac{|\widehat{B}_{\widehat{U}_1}|}{|\widehat{U}_1|}&=\frac{k^{\underline{\theta}}([0;\overline{1,3}])-k^{\underline{\theta}}([0;1,2,\overline{3,1}])}{k^{\underline{\theta}}([0;1,2,\overline{3,1}])-k^{\underline{\theta}}([0;1,1,\overline{1,3}])}=\frac{4\sqrt{21}+3}{5}=1.822\dots \\
    \frac{|\widehat{B}_{\widehat{U}_2}|}{|\widehat{U}_2|}&=\frac{k^{\underline{\theta}}([0;\overline{1,3}])-k^{\underline{\theta}}([0;1,3,\overline{3,1}])}{k^{\underline{\theta}}([0;1,3,\overline{3,1}])-k^{\underline{\theta}}([0;1,2,\overline{1,3}])}=\frac{4\sqrt{21}}{21}=0.872\dots \\
    \frac{|\widehat{B}_{\widehat{U}_3}|}{|\widehat{U}_3|}&=\frac{k^{\underline{\theta}}([0;\overline{1,3}])-k^{\underline{\theta}}([0;1,3,1,\overline{1,3}])}{k^{\underline{\theta}}([0;1,3,1,\overline{1,3}])-k^{\underline{\theta}}([0;1,3,2,\overline{3,1}])}=\frac{8\sqrt{21}+28}{5}=1.847\dots
\end{align*}
Therefore $\tau_p(C_3,\max C_3)=4\sqrt{21}/21$. Notice that this number is greater than (usual) thickness of $C_3$ which is $(-6+4\sqrt{21})/15=0.822\dots$
\end{example}

\begin{remark}
    From the fact that eventually periodic continued fractions are quadratic irrationals, it follows that the pointwise thickness of Gauss-Cantor sets are always algebraic numbers over $\Q$. 
\end{remark}

\subsection{Main lemmas}
Now let us explain what is the mechanism behind the proof of the result of \cite{Newgaps}, which in turn is based on \cite[Theorem IV.1]{GuguThesis}. The idea is that provided some irrationality condition holds at the maximum of the respective Cantor sets and there are gaps whose bridges are adjacent to these maximum and for which the product of the thickness ratios is small, by density we will be able to choose arbitrarily small and controlled scales for which we can repeat a non intersecting pattern infinitely many times (see \Cref{fig:scale_recurrence_figure}).

In order to obtain estimates for the Hausdorff dimension, we will use the fact that irrational rotations equidistribute on the circle. Considering the usual distance $d(x,y)=|x-y| \pmod{1}$ on $\R\backslash\Z$, for an irrational $\rho\in(0,1)$ define the covering radius of the first $N+1$ points to be
\begin{equation*}
    R_N(\rho)=\sup_{x\in \R\backslash\Z}\min_{0\leq n\leq N} d(x,n\rho).
\end{equation*}
Define the $\delta$-covering time to be
\begin{equation}\label{eq:delta_covering_time}
    N_{\delta}(\rho)=\min\{N\geq 0: R_N(\rho)\leq\delta\}.
\end{equation}

\begin{lemma}\label{lem:main_lemma}
	Let $K$ and $\widetilde{K}$ be dynamically defined Cantor sets with Markov partitions $\{I_1,\dots,I_r\}$ and $\{\tilde{I}_1,\dots,\tilde{I}_{\tilde{r}}\}$ and expanding maps $\psi$ and $\tilde{\psi}$, respectively. Let $U^\prime$ be a bounded gap of $K$ and let $p=\min\partial U^\prime\in K$ be such that 
    \begin{equation}\label{eq:only_one_point_of_intersection_condition}
        (x,\tilde{x})\in K\times\widetilde{K},~ x+\tilde{x}=p+\max\widetilde{K} \quad\implies\quad (x,\tilde{x})=(p,\max\widetilde{K}).
    \end{equation}
    
    Let $q\in K$, $\tilde{q}\in\widetilde{K}$ and $\ell,\tilde{\ell},k,\widetilde{k}$ be non-negative integers such that 
    \begin{equation*}
        \psi^\ell(p)=q, \quad \psi^k(q)=q, \quad\text{and}\quad \tilde{\psi}^{\tilde{\ell}}(\max\widetilde{K})=\tilde{q}, \quad  \tilde{\psi}^{\widetilde{k}}(\tilde{q})=\tilde{q}.
    \end{equation*}
    Suppose that
    \begin{equation}\label{eq:irrationality_condition}
        \frac{\log|(\psi^k)^\prime(q)|}{\log|(\tilde{\psi}^{\widetilde{k}})^\prime(\tilde{q})|}\in\R\setminus\Q.
    \end{equation}
    Let us denote by $I(\underline{a})$, $I(\underline{b})$ and $\tilde{I}(\underline{\tilde{a}})$, $I(\underline{\tilde{b}})$ be intervals of the construction of $K$ and $\widetilde{K}$ respectively such that 
    \begin{equation*}
        h_1(q)=p, \quad  f_1(q)=q, \quad\text{and}\quad h_2(\tilde{q})=\max\widetilde{K}, \quad f_2(\tilde{q})=\tilde{q}.
    \end{equation*}
    where the inverse branches are
    \begin{equation*}
        h_1=(\psi^{\ell}|_{I(\underline{a})})^{-1},  f_1=(\psi^{k}|_{I(\underline{b})})^{-1}, \quad\text{and}\quad h_2=(\tilde{\psi}^{\tilde{\ell}}|_{I(\underline{\tilde{a}})})^{-1}, f_2=(\tilde{\psi}^{\widetilde{k}}|_{I(\underline{\tilde{b}})})^{-1}.
    \end{equation*}
    Suppose that $U\subset I(\underline{b})$ and $\widetilde{U}\subset I(\underline{\tilde{b}})$ are gaps of $K$ and $\widetilde{K}$, respectively, such that for the adjacent interval $B_U\subset I(\underline{b})$ of $U$ with $q\in\partial B_U$ and the adjacent interval $B_{\widetilde{U}}\subset I(\underline{\tilde{b}})$ of $\widetilde{U}$ with $\tilde{q}\in\partial B_{\widetilde{U}}$, it holds that
    \begin{equation}\label{eq:bridges_gap_quotient}
        \liminf_{n\to\infty}\frac{|f_1^n(B_U)|}{|f_1^n(U)|}\frac{|f_2^n(B_{\widetilde{U}})|}{|f_2^n(\widetilde{U})|}<1.
    \end{equation}
    Then $\partial(K+\widetilde{K})\cap[p+\max\widetilde{K}-\varepsilon,p+\max\widetilde{K}]$ is infinite for all $\varepsilon>0$. More precisely, there are infinitely many positive integers $(m,n)$ such that $K+\widetilde{K}$ has a gap from
    \begin{equation*}
        \mu_n=\min(\partial (h_1\circ f_1^n(U)))+\max\widetilde{K},
    \end{equation*}
    to
    \begin{equation*}
        \nu_n=\min (h_1\circ f_1^n(B_U)) + \min (h_2\circ f_2^m(B_{\widetilde{U}})).
    \end{equation*}
    Moreover,
    \begin{equation}\label{eq:only_one_point_of_intersection_condition_smaller_gap}
        (x,\tilde{x})\in K\times\widetilde{K},~ x+\tilde{x}=\mu_n \quad\implies\quad (x,\tilde{x})=(\min(\partial (h_1\circ f_1^n(U))),\max\widetilde{K}).
    \end{equation}

    Furthermore, if we denote $n_1<n_2<n_3<\dots$ a sequence of positive integers such that $(\mu_{n_i},\nu_{n_i})$ is a gap, then we have that 
    \begin{equation}\label{eq:bounded_steps}
        n_{i+1}-n_i \leq N_{\delta^\prime}(\rho)
    \end{equation}
    where $N_{\delta^\prime}(\rho)$ is the $\delta^\prime$-covering time defined on \eqref{eq:delta_covering_time} and for some large $N$
    \begin{equation}\label{eq:delta_prime_size}
        \delta^\prime = \frac{1}{2\log|(\tilde{\psi}^{\widetilde{k}})^\prime(\tilde{q})|}\log\left(\frac{|f_1^N(U)|}{|f_2^N(B_{\widetilde{U}})|}\cdot\frac{|f_1^N(\widetilde{U})|+|f_2^N(B_{\widetilde{U}})|}{|f_1^N(U)|+|f_1^N(B_{U})|}\right), \qquad \rho = \frac{\log|(\psi^k)^\prime(q)|}{\log|(\tilde{\psi}^{\widetilde{k}})^\prime(\tilde{q})|}
    \end{equation}
\end{lemma}

\begin{proof}
    Observe that the hypothesis \eqref{eq:only_one_point_of_intersection_condition} implies that there is an $\varepsilon>0$ such that $(K+\tilde{K})\cap I_{\varepsilon}\subset \left((K\cap I(\underline{a}))+\widetilde{K}\right)\cap I_{\varepsilon}$ where we denote $I_{\varepsilon}=[p+\max\widetilde{K}-\varepsilon,p+\max\widetilde{K}]$. Otherwise there would be a decreasing sequence of $\varepsilon_n>0$ and $(x_n,\tilde{x}_n)\in K\times\widetilde{K}$ such that $x_n+\tilde{x}_n\in I_{\varepsilon}$ but $x_n\not\in I(\underline{a})$, since both $K\cap(\R\setminus I(\underline{a})$ and $\tilde{K}$ are compact by taking a subsequence as $n\to\infty$ we would find $x_{\infty}\in K\cap(\R\setminus I(\underline{a})$ and $\tilde{x}\in\widetilde{K}$ satisfying \eqref{eq:only_one_point_of_intersection_condition}, a contradiction.

    Because of the condition \eqref{eq:bridges_gap_quotient}, we can replace $U$, $B_{U}$, $\widetilde{U}$, $B_{\widetilde{U}}$ by $f_1^{N}(U)$,  $f_1^{N}(B_{U})$, $f_2^{N}(\widetilde{U})$, $f_2^{N}(B_{\widetilde{U}})$ (for some large $N$) to assume that $U$, $B_{U}$, $\widetilde{U}$, $B_{\widetilde{U}}$ have sufficiently small diameter, say $\max\{|U|+|B_{U}|,|\widetilde{U}|+|B_{\widetilde{U}}|\}<\delta$ for some $\delta>0$ we will choose later.
    
    Let us denote
    \begin{equation*}
        U_n = h_1\circ f_1^n(U), J_n = h_1\circ f_1^n(B_U), \quad\text{and}\quad \widetilde{U}_m = h_2 \circ f_2^m(\widetilde{U}), \widetilde{J}_m = h_2 \circ f_2^m(B_{\widetilde{U}}).
    \end{equation*}
    By \eqref{eq:bridges_gap_quotient}, and by the irrationality condition \eqref{eq:irrationality_condition}, by making $\delta>0$ arbitrarily small, we can obtain a sufficiently small $C(\delta)$ such that we can choose positive integers $m,n$ with
    \begin{equation}\label{eq:infinitely_many_gaps_condition}
        e^{C(\delta)}\frac{|h_1^\prime(q)|}{|h_2^\prime(\tilde{q})|}\cdot\frac{|B_{U}|}{|\widetilde{U}|}<\frac{|f_2^\prime(\tilde{q})|^m}{|f_1^\prime(q)|^n}<e^{-C(\delta)}\frac{|h_1^\prime(q)|}{|h_2^\prime(\tilde{q})|}\cdot\frac{|U|}{|B_{\widetilde{U}}|}.
    \end{equation}
    Therefore
    \begin{align*}
        \frac{|\widetilde{J}_m|}{|U_n|}&<e^{C(\delta)}\frac{|h_2^\prime(\tilde{q})|}{|h_1^\prime(q)|}\frac{|f_2^\prime(\tilde{q})|^m}{|f_1^\prime(q)|^n}\frac{|B_{\widetilde{U}}|}{|U|}<1, \\
        \frac{|J_n|}{|\widetilde{U}_m|}&<e^{C(\delta)}\frac{|h_1^\prime(q)|}{|h_2^\prime(\tilde{q})|}\frac{|f_1^\prime(q)|^n}{|f_2^\prime(\tilde{q})|^m}\frac{|B_U|}{|\widetilde{U}|}<1.
    \end{align*}

    \begin{figure}
    \centering
	\begin{tikzpicture}[line width=0.8pt]

    \draw (-5,0) -- (-0.5,0);
    \draw (-5,0.18) -- (-5,-0.18);
    
    \node at (-0.5,0) {\big (};
    \node at (1.2,0) {\big )};
    
    \draw (1.2,0) -- (2,0);
    \draw (2,0.18) -- (2,-0.18);
    
    \node[above] at (0.35,0.25) {$U_n$};
    \node[above] at (1.6,0.25) {$J_n$};
    
    \draw (-0.1,-1.15) -- (0.9,-1.15);
    \draw (-0.1,-1.33) -- (-0.1,-0.97);
    
    \node at (0.9,-1.15) {\big (};
    \node at (3,-1.15) {\big )};
    
    \draw (3,-1.15) -- (4.8,-1.15);
    \draw (4.8,-1.33) -- (4.8,-0.97);
    
    \node[below] at (0.4,-1.25) {$t-\widetilde{J}_m$};
    \node[below] at (2,-1.25) {$t-\widetilde{U}_m$};
    
    \end{tikzpicture}
    \caption{The exact moment $t-\widetilde{J}_m$ touches the left border of $U_n$, the right border of $J_n$ is far from the boundary of $t-\widetilde{U}_m$, provided that $|J_n|+|U_n|<|\widetilde{J}_m|+|\widetilde{U}_m|$. Moreover for that parameter $t = \min\partial U_n+\max\widetilde{J}_m$ the intersection occurs at exactly one point.}
    \label{fig:scale_recurrence_figure}
    \end{figure}
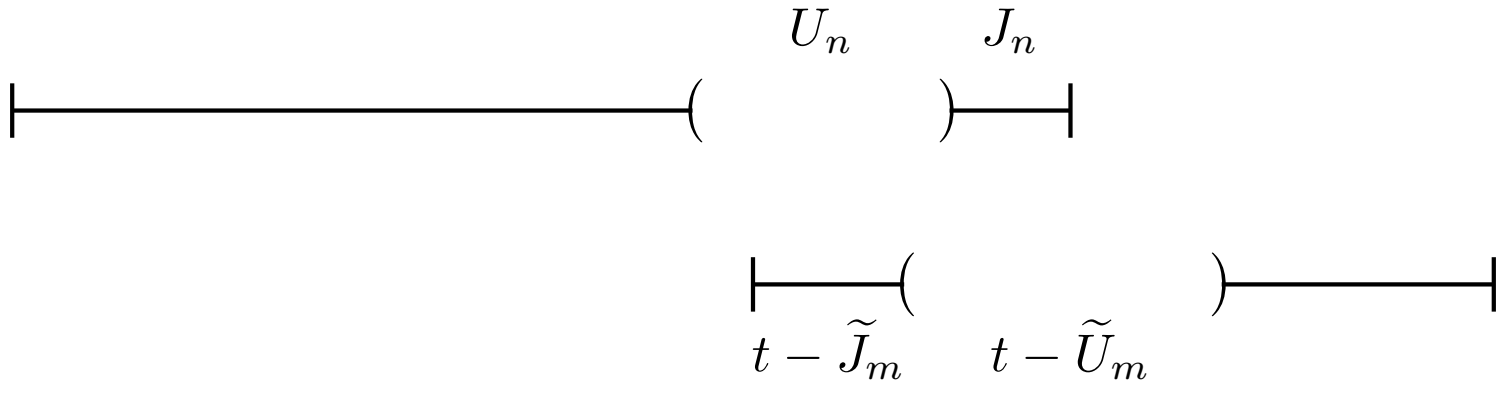

    These inequalities (see \Cref{fig:scale_recurrence_figure}) imply that $K+\widetilde{K}$ has a gap from 
    \begin{equation*}
        \max\{\max J_n+\min\partial \widetilde{U}_m,\max \widetilde{J}_m+\min\partial U_n\}
    \end{equation*}
    to $\min J_n+\min \widetilde{J}_m$ (observe that these points belong to $K+\widetilde{K}$). Since all these borders converge to $\max K$ and $\max\widetilde{K}$, we find infinitely many such gaps by taking arbitrarily large values of $m$ and $n$. 

    We can in fact control that the left border of the gaps is precisely $\max\widetilde{J}_m+\min\partial U_n$. The irrationality condition \eqref{eq:irrationality_condition} allow us to choose positive integers $m$ and $n$ satisfying \eqref{eq:infinitely_many_gaps_condition} and 
    \begin{equation}\label{eq:infinitely_many_gaps_condition2}
        e^{2C(\delta)}\frac{|h_1^\prime(q)|}{|h_2^\prime(\tilde{q})|}\cdot\frac{|B_U|+|U|}{|B_{\widetilde{U}}|+|\widetilde{U}|}<\frac{|f_2^\prime(\tilde{q})|^m}{|f_1^\prime(q)|^n}<e^{-C(\delta)}\frac{|h_1^\prime(q)|}{|h_2^\prime(\tilde{q})|}\cdot\frac{|U|}{|B_{\widetilde{U}}|}.
    \end{equation}
    Indeed, notice that \eqref{eq:bridges_gap_quotient} implies that 
    \begin{equation*}
        \frac{|B_U|+|U|}{|B_{\widetilde{U}}|+|\widetilde{U}|}\cdot\frac{|B_{\widetilde{U}}|}{|U|} = \frac{|B_U|/|U|+1}{|\widetilde{U}|/|B_{\widetilde{U}}|+1}<1,
    \end{equation*}
    and
    \begin{equation*}
        \frac{|B_U|+|U|}{|B_{\widetilde{U}}|+|\widetilde{U}|}\cdot\frac{|\widetilde{U}|}{|B_{U}|} = \frac{|U|/|B_U|+1}{|B_{\widetilde{U}}|/|\widetilde{U}|+1}>1
    \end{equation*}
    Therefore, we have that \eqref{eq:infinitely_many_gaps_condition2} implies condition \eqref{eq:infinitely_many_gaps_condition}. Moreover, the condition \eqref{eq:infinitely_many_gaps_condition2} is equivalent to require that 
    \begin{equation*}
        n\cdot\frac{\log|(\psi^k)^\prime(q)|}{\log|(\tilde{\psi}^{\widetilde{k}})^\prime(\tilde{q})|}-m 
    \end{equation*}
    belong to an interval of size
    \begin{equation*}
        \frac{1}{\log|(\tilde{\psi}^{\widetilde{k}})^\prime(\tilde{q})|}\left(-3C(\delta)+\log\left(\frac{|U|}{|B_{\widetilde{U}}|}\right) - \log\left(\frac{|B_U|+|U|}{|B_{\widetilde{U}}|+|\widetilde{U}|}\right)\right).
    \end{equation*}
    Since we can chose $\delta>0$ small enough so that $C(\delta)>0$ is small enough, we can assume that the previous interval is always of size at least \eqref{eq:delta_prime_size}. 
    
    Since by bounded distortion we have
    \begin{align*}
        |J_n|+|U_n|&=|f_1^n(U\cup J)|<e^{C(\delta)}|h_1^\prime(q)||f_1^\prime(q)|^n(|B_U|+|U|) \\
        |\widetilde{J}_m|+|\widetilde{U}_m|&=|f_2^m(\widetilde{U}\cup \widetilde{J})|>e^{-C(\delta)}|h_2^\prime(\tilde{q})||f_2^\prime(\tilde{q})|^m(|B_{\widetilde{U}}|+|\widetilde{U}|),
    \end{align*}
    by \eqref{eq:infinitely_many_gaps_condition2} we have
    \begin{equation*}
        \frac{|J_n|+|U_n|}{|\widetilde{J}_m|+|\widetilde{U}_m|}<e^{2C(\delta)}\frac{|h_1^\prime(q)|}{|h_2^\prime(\tilde{q})|}\frac{|f_1^\prime(q)|^n}{|f_2^\prime(\tilde{q})|^m}\frac{|B_U|+|U|}{|B_{\widetilde{U}}|+|\widetilde{U}|}<1.
    \end{equation*}
    
    The last condition implies that the right border of $\widetilde{J}_m$ touches the left border of $U_n$ much before the right border of $J_n$ touches the left border of $\widetilde{U}_m$ in \Cref{fig:scale_recurrence_figure} (we are imaging that the Cantor set $t-K$ is moving to the left, i.e., that we are decreasing $t$). Moreover, since the first point of intersection between the sets is unique, we obtain immediately \eqref{eq:only_one_point_of_intersection_condition_smaller_gap}. 
\end{proof}

\begin{remark}
    In particular, if $\tau_{p}(K,q)\cdot\tau_{p}(\widetilde{K},\tilde{q})<1$, the condition \eqref{eq:bridges_gap_quotient} is satisfied for some pair of gaps $U$ and $\widetilde{U}$ and conversely the condition \eqref{eq:bridges_gap_quotient} imply that $\tau_{p}(K,q)\cdot\tau_{p}(\widetilde{K},\tilde{q})<1$.
\end{remark}

\begin{figure}[ht]
\centering
\begin{tikzpicture}
    \tikzset{
        pics/mydiagram/.style args={#1/#2}{
            code={
                \node at (0,0) {$($};
                \node at (0.4,0) {$)$};
    
                \node at (1.2,0) {$($};
                \node at (1.6,0) {$)$};
    
                \node at (2.4,0) {$($};
                \node at (2.7,0) {$)$};
    
                \node[transform shape] at (3.2,0) {$\cdots$};
    
                \node at (4,0) {$($};
                \node at (6,0) {$)$};
        
                \node[transform shape=false]
                    at ([yshift=#2]5,0.5) {$#1$};
            }
        }
    }
    
    \pic[scale=0.8] at (0,0) {mydiagram={W_1}/0.5mm};
    \pic[scale=0.5] at (6,0) {mydiagram={W_2}/4mm};
    \pic[scale=0.3] at (10,0) {mydiagram={W_3}/10mm};
    
    \node at (12.5,0) {$\cdots$};
\end{tikzpicture}
\caption{Each gap has infinitely many gaps to the left approaching its left border and each new gap appearing this way has the same property.}
\label{fig:uncountable_boundary}
\end{figure}
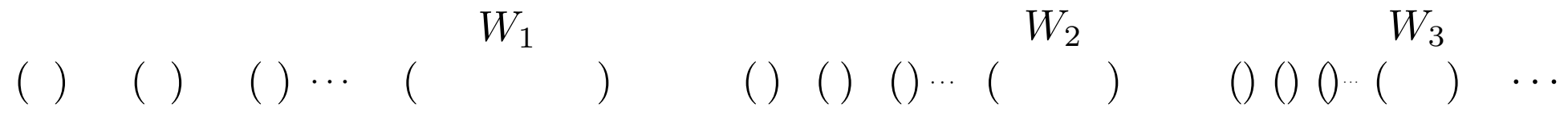

Now we will explain what is the mechanism for improving the situation from infinitely many gaps to positive Hausdorff dimension, and thus, uncountable boundary. The key idea is to construct an infinite family of gaps (of the arithmetic sum) with the property that for any gap $(\nu,\mu)$ in the family there are infinitely many gaps $(\nu_n,\mu_n)$ in the same family such that they are approaching the left border $\nu$, namely, $\nu_n<\mu_n<\nu_{n+1}<\mu_{n+1}<\dots<\nu$ for all $n$ and $\nu_n,\mu_n\to \nu$ as $n\to\infty$ (see \Cref{fig:uncountable_boundary}). Clearly this implies that such a set has uncountable boundary. The mechanism behind this construction is to apply the previous lemma recursively to each new gap provided we have a recursive pattern in our scales and periodic points. Since we can attain such a recursive patter in bounded times \eqref{eq:bounded_steps} (because of the equidistribution of the irrational rotation in finite time) we will be able to obtain positive Hausdorff dimension of the collection of left borders $\nu_n$ by using the mass distribution principle.

\begin{proposition}\label{prop:mass_distribution}
    Let $\mu$ be a mass distribution of a set $F$ and suppose that for some $d\geq 0$, there are $c>0$ and $\varepsilon>0$ such that $\mu(U)\leq c|U|^d$, for all intervals $U$ with $|U|\leq\varepsilon$. Then $\mathcal{H}^d(F)\geq\mu(F)/c$ and $d\leq\dim_H F\leq\underline{\dim}_BF\leq\overline{\dim}_BF$.
\end{proposition}

\begin{lemma}\label{lem:specialized_lemma}
    Let $\{q_1,\dots,q_m\}\subset K$ be a finite set of periodic points with $\psi^{k_i}(q_i)=q_i$, $k_i\geq 1$ such that for any boundary point $p\in\partial U$ of a gap of the Cantor set $K$ there is some $\ell_p$ and $1\leq i\leq m$ such that $\psi^{\ell_p}(p)=q_i$. 

    Let $\tilde{q}\in\widetilde{K}$ and $\tilde{\ell},\widetilde{j}$ be non-negative integers such that $\tilde{\psi}^{\tilde{\ell}}(\max\widetilde{K})=\tilde{q}$ and $\tilde{\psi}^{\widetilde{k}}(\tilde{q})=\tilde{q}$. Let us fix intervals $I(\underline{\tilde{a}})$ and $I(\underline{\tilde{b}})$ of the construction of $\widetilde{K}$ and a gap $\widetilde{U}\subset I(\underline{\tilde{b}})$ such that $(\tilde{\psi}^{\tilde{\ell}}|_{I(\underline{\tilde{a}})})^{-1}(\tilde{q})=\max\widetilde{K}$ and $(\tilde{\psi}^{\widetilde{k}}|_{I(\underline{\tilde{b}})})^{-1}(\tilde{q})=\tilde{q}$. 
    
    Suppose that $p_0=\max K$ and $p_1,\dots,p_l\in K$, $l\geq 1$, with associated periodic points $\psi^{\ell_{p_i}}(p_i)=q_{e_i}$ for some $1\leq e_i\leq m$ are such that:
    \begin{itemize}
        \item there are $I(\underline{a_i})$ and $I(\underline{b_i})$ intervals of the construction of $K$ such that  
        \begin{equation*}
            (\psi^{\ell_{p_i}}|_{I(\underline{a_i})})^{-1}(q_{e_i})=p_i, \quad\text{and}\quad (\psi^{k_{e_i}}|_{I(\underline{b_i})})^{-1}(q_{e_i})=q_{e_i}, \quad\text{for } i=0,\dots,l-1.
        \end{equation*}
        \item there are gaps $U_1\subset I(\underline{b_0}),\dots,U_l\subset I(\underline{b_{l-1}})$ of $K$ and adjacent intervals $B_{U_1}\subset I(\underline{b_0}),\dots,B_{U_l}\subset I(\underline{b_{l-1}})$ (i.e, $U_i\cap B_{U_i}=\varnothing$ and $\partial U_i\cap B_{U_i}\neq\varnothing$) such that 
        \begin{equation*}
            q_{e_i}\in \partial B_{U_{i+1}}, \, \, \, \partial(U_{i+1}\cup B_{U_{i+1}})=\{q_{e_i},p_{i+1}\} \quad\text{for } i=0,\dots,l-1.
        \end{equation*}
        \item 
        \begin{equation}\label{eq:irrationality_condition_2}
            \frac{\log|(\psi^{k_{e_i}})^\prime(q_{e_i})|}{\log|(\tilde{\psi}^{\widetilde{k}})^\prime(\tilde{q})|}\in\R\setminus\Q, \quad\text{for } i=0,\dots,l-1.
        \end{equation}
        \item 
        \begin{equation}\label{eq:bridges_gap_quotient_2}
            \liminf_{n\to\infty}\frac{|(\psi^{k_{e_i}}|_{I(\underline{b_i})})^{-n}(B_{U_{i+1}})|}{|(\psi^{k_{e_i}}|_{I(\underline{b_i})})^{-n}(U_{i+1})|}\frac{|(\tilde{\psi}^{\tilde{k}}|_{I(\underline{\tilde{b}})})^{-n}(B_{\widetilde{U}})|}{|(\tilde{\psi}^{\tilde{k}}|_{I(\underline{\tilde{b}})})^{-n}(\widetilde{U})|}<1, \quad\text{for } i=0,\dots,l-1.
        \end{equation}
    \end{itemize}
    Finally, suppose that $q_{e_l}=q_{e_0}$. Then $\partial(K+\widetilde{K})\cap[\max(K+\widetilde{K})-\varepsilon,\max(K+\widetilde{K})]$ has positive Hausdorff dimension for all $\varepsilon>0$. In particular, if $\tau_{p}(K,q_{e_i})\cdot\tau_{p}(\widetilde{K},\tilde{q})<1$, $i=0,\dots,l-1$, then such gaps exists. 
\end{lemma}

\begin{proof}
    By \Cref{lem:main_lemma} there is a sequence of gaps $W_1,W_2,W_3,\dots$ of $K+\widetilde{K}$ with $\min\partial W_j\to\max(K+\widetilde{K})$ and such that there are positive integers $n_j$ such that
    \begin{equation}\label{eq:gaps_W_j}
        \min \partial W_j = (\psi^{\ell_{p_0}}|_{I(\underline{a_0})})^{-1}\circ(\psi^{k_{e_0}}|_{I(\underline{b_0})})^{-n_j}(p_{1})+\max\widetilde{K}.
    \end{equation}
    Moreover, the point $p=(\psi^{\ell_{p_0}}|_{I(\underline{a_0})})^{-1}\circ(\psi^{k_{e_0}}|_{I(\underline{b_0})})^{-n_j}(p_{1})$ is the left border of a gap of $K$, namely the gap $(\psi^{\ell_{p_0}}|_{I(\underline{a_0})})^{-1}\circ(\psi^{k_{e_0}}|_{I(\underline{b_0})})^{-n_j}(U_1)$. Furthermore, by \eqref{eq:only_one_point_of_intersection_condition_smaller_gap} we have
    \begin{equation*}
        (x,\tilde{x})\in K\times\widetilde{K},~ x+\tilde{x}=p+\max\widetilde{K} \quad\implies\quad (x,\tilde{x})=(p,\max\widetilde{K}).
    \end{equation*}
    In particular by hypothesis there is $\ell_p=\ell_{p_0}+k_{e_0}n_j+\ell_{p_1}$ such that $\psi^{\ell_p}(p)=q_{e_1}$.

    Now we will use \Cref{lem:main_lemma} recursively. Suppose that $V$ is any gap of $K$ such that $p=\min\partial V$ satisfies $\psi^{\ell_p}(p)=q_{e_j}\in\{q_{e_0},q_{e_1},\dots,q_{e_{l-1}}\}$ for some $\ell_p\geq 0$ and
    \begin{equation*}
        (x,\tilde{x})\in K\times\widetilde{K},~ x+\tilde{x}=p+\max\widetilde{K} \quad\implies\quad (x,\tilde{x})=(p,\max\widetilde{K}).
    \end{equation*}
    We can apply \Cref{lem:main_lemma} directly to this gap to show that $K+\widetilde{K}$ contains infinitely many gaps of the form $(\mu,\nu)$ where
    \begin{equation*}
        \mu=\min(\partial (h_1\circ f_1^n(U_j))))+\max\widetilde{K},
    \end{equation*}
    and
    \begin{equation*}
        \nu=\min (h_1\circ f_1^n(B_{U_j})) + \min (h_2\circ f_2^m(B_{\widetilde{U}})).
    \end{equation*}
    where 
    \begin{equation*}
        h_1=(\psi^{\ell_p}|_{I(\underline{a})})^{-1},  f_1=(\psi^{k_{e_j}}|_{I(\underline{b_j})})^{-1}, \quad\text{and}\quad h_2=(\tilde{\psi}^{\tilde{\ell}}|_{I(\underline{\tilde{a}})})^{-1}, f_2=(\tilde{\psi}^{\widetilde{k}}|_{I(\underline{\tilde{b}})})^{-1}, 
    \end{equation*}
    and $I(\underline{a})$ is an interval of the construction of $K$ such that $h_1(q_{e_j})=p$. By construction the point $\hat{p}=(\psi^{\ell_{p}}|_{I(\underline{a})})^{-1}\circ(\psi^{k_{e_j}}|_{I(\underline{b_j})})^{-n}(p_{j})$ is the left border of a gap of $K$, namely the gap $(\psi^{\ell_{p}}|_{I(\underline{a})})^{-1}\circ(\psi^{k_{e_j}}|_{I(\underline{b_j})})^{-n}(U_j)$. Furthermore, by \eqref{eq:only_one_point_of_intersection_condition_smaller_gap} we have
    \begin{equation*}
        (x,\tilde{x})\in K\times\widetilde{K},~ x+\tilde{x}=\hat{p}+\max\widetilde{K} \quad\implies\quad (x,\tilde{x})=(\hat{p},\max\widetilde{K}).
    \end{equation*}
    In particular by hypothesis there is $\ell_{\hat{p}}=\ell_{p}+k_{e_j}n+\ell_{p_j}$ such that $\psi^{\ell_{\hat{p}}}(\hat{p})=q_{e_{j+1}}$. In conclusion, the gap $(\psi^{\ell_{p}}|_{I(\underline{a})})^{-1}\circ(\psi^{k_{e_j}}|_{I(\underline{b_j})})^{-n}(U_j)$ has the same properties as the gap $V$, so we finish the recursion.

    \begin{figure}
    \centering
	\begin{tikzpicture}[line width=0.8pt]

    \draw (-2,0) -- (-0.5,0);
    \draw (-2,0.18) -- (-2,-0.18);

    \node at (-0.5,0) {\big (};
    \node at (1.2,0) {\big )};
    
    \draw (1.2,0) -- (2,0);
    \draw (2,0.18) -- (2,-0.18);

    \draw (5.5,0.18) -- (5.5,-0.18);
    \draw (5.5,0) -- (5.8,0);
    \draw[thick, densely dotted] (5.8,0) -- (6.5,0);
    \draw (6.5,0) -- (6.8,0);
    \draw (6.8,0.18) -- (6.8,-0.18);
    
    \node[above] at (0.35,0.25) {$V$};
    \node[above] at (1.6,0.25) {$J$};
    
    \draw (-0.1,-1.15) -- (0.9,-1.15);
    \draw (-0.1,-1.33) -- (-0.1,-0.97);
    
    \node at (0.9,-1.15) {\big (};
    \node at (3,-1.15) {\big )};
    
    \draw (3,-1.15) -- (4.8,-1.15);
    \draw (4.8,-1.33) -- (4.8,-0.97);

    \draw (7.5,-1.33) -- (7.5,-0.97);
    \draw (7.5,-1.15) -- (7.8,-1.15);
    \draw[thick, densely dotted] (7.8,-1.15) -- (8.5,-1.15);
    \draw (8.5,-1.15) -- (8.8,-1.15);
    \draw (8.8,-1.33) -- (8.8,-0.97);
    
    \node[below] at (0.4,-1.25) {$t-B_{\widetilde{U}}$};
    \node[below] at (2,-1.25) {$t-\widetilde{U}$};

    \fill (-0.5,0) circle (2pt);
    \node[above,yshift=1em] at (-0.5,0) {$p$};
    \fill (6.8,0) circle (2pt);
    \node[above,yshift=1em] at (6.8,0) {$\max K$};

    \fill (-0.1,-1.15) circle (2pt);
    \node[left,xshift=-0.5em] at (-0.1,-1.15) {$t-\max\widetilde{K}$};

    \fill (8.8,-1.15) circle (2pt);
    \node[below,yshift=-0.5em] at (8.8,-1.15) {$t-\min\widetilde{K}$};
    
    \end{tikzpicture}
    \caption{We can chose intervals $U,\widetilde{U},J,\widetilde{J}$ arbitrarily small. }
    \end{figure}
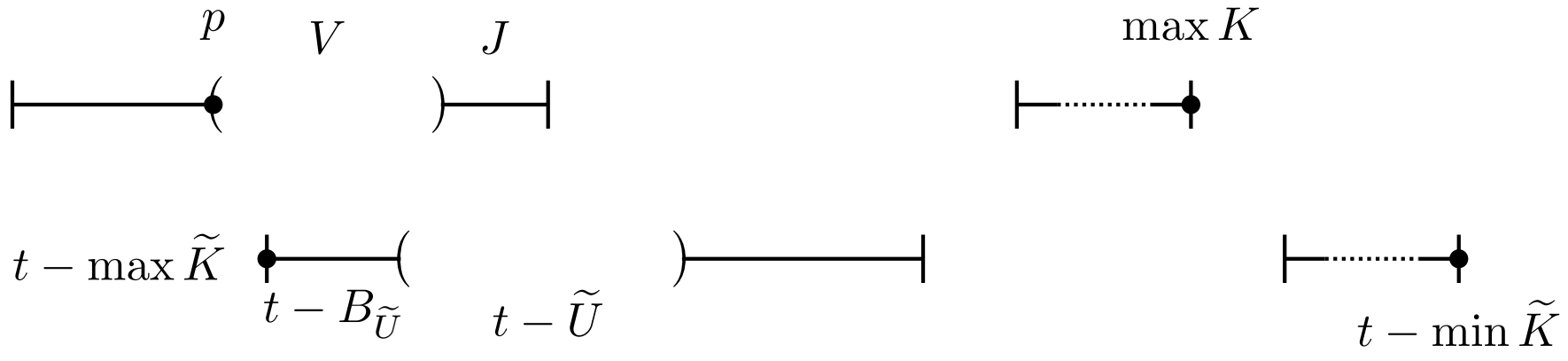

    Therefore we obtain uncountably many points of the form \\
    \begin{multline*}
        (\psi^{\ell_{p_0}}|_{I(\underline{a_0})})^{-1}\circ(\psi^{k_{e_0}}|_{I(\underline{b_0})})^{-n_1}\circ(\psi^{\ell_{p_1}}|_{I(\underline{a_1})})^{-1}\circ(\psi^{k_{e_1}}|_{I(\underline{b_1})})^{-n_2}\circ\dots \\
        \circ(\psi^{\ell_{p_j}}|_{I(\underline{a_j})})^{-1}\circ(\psi^{k_{e_j}}|_{I(\underline{b_j})})^{-n_j}(p_{j+1})
    \end{multline*}
    where we have extended cyclically $p_j = p_{j-l}$ and $\underline{a_j}=\underline{a_{j-l}}$ for $j\geq l+1$ and $\underline{b_j}=\underline{b_{j-l}}$  and $k_{e_j}=k_{e_{j-l}}$ for $j\geq l$. Moreover, because of the condition \eqref{eq:bounded_distortion}, the numbers $n_1,n_2,\dots$ lie in a syndetic set of choices, with step uniformly bounded by some absolute $N$, that depends on the multipliers \eqref{eq:irrationality_condition_2} and the initial gaps and adjacent intervals \eqref{eq:bridges_gap_quotient_2}. We can apply the mass transference principle to bound by below the Hausdorff dimension of such a set. Indeed, let us consider three consecutive gaps $W_1$, $W_2$ and $W_3$ at the first step in our construction, equivalently, consider three consecutive choices for $n_1$. We give equal mass of $1/2$ to the initial intervals $I_0 = [\max\partial W_1,\min\partial W_2]$ and $I_1=[\max\partial W_2,\min\partial W_3]$. Inside of each interval $I_0$ and $I_1$, we again chose 3 consecutive gaps (equivalently, 3 consecutive choices for $n_2$) and give equal mass to the four intervals so created. Given a finite word $(w_1,\dots,w_n)\in\{0,1\}^n$, at each step we subdivide $I_{w}$ into the intervals $I_{w0}$ and $I_{w1}$ with equal mass. Moreover, by construction, we have that $|I_{w0}|,|I_{w1}|\geq \rho^\prime|I_w|$ for some absolute constant $\rho^\prime$, that can be taken as
    \begin{equation*}
        \rho^\prime=\inf_{1\leq i\leq l}\inf_{x\in I(\underline{a_i})}\inf_{y\in I(\underline{b_i})}\left(|(\psi^{\ell_{p_i}}|_{I(\underline{a_i})})'(x)|\cdot|(\psi^{Nk_{e_i}}|_{I(\underline{b_i})})^\prime(y)|\right)^{-1}.
    \end{equation*}
    Since we can apply the very same argument to any triple of consecutive gaps $W_1$, $W_2$ and $W_3$ and these gaps are approaching $\max K+\max\tilde{K}$, the conclusion of the lemma follows.
\end{proof}

\begin{corollary}\label{cor:conjugating_preserves_lemma}
    Suppose that $g_1$ and $g_2$ are both increasing $C^1$ diffeomorphism defined on open neighborhoods of $K$ and $\widetilde{K}$ respectively. If $(K,\widetilde{K})$ satisfy the hypothesis of \Cref{lem:specialized_lemma}, then $(g_1(K),g_2(\widetilde{K}))$ also satisfy the hypothesis of \Cref{lem:specialized_lemma}. In particular if both $g_1$ and $g_2$ are increasing then for all $\varepsilon>0$ the set
    \begin{equation*}
        \partial(g_1(K)+g_2(\widetilde{K}))\cap[g_1(\max K)+g_2(\max\widetilde{K})-\varepsilon,g_1(\max K)+g_2(\max\widetilde{K})],
    \end{equation*}
    has positive Hausdorff dimension, and if $g_1$ and $g_2$ are decreasing then then for all $\varepsilon>0$ the set
    \begin{equation*}
        \partial(g_1(K)+g_2(\widetilde{K}))\cap[g_1(\min K)+g_2(\min\widetilde{K}),g_1(\min K)+g_2(\min\widetilde{K})+\varepsilon],
    \end{equation*}
    has positive Hausdorff dimension. 
\end{corollary}

\begin{proof}
    By bounded distortion we can replace $U_i$, $B_{U_i}$, $\widetilde{U}$, $B_{\widetilde{U}}$ by $(\psi^{k_{e_i}}|_{I(\underline{b_i})})^{-1}(U_i)$,  $(\psi^{k_{e_i}}|_{I(\underline{b_i})})^{-1}(B_{U_i})$, $(\tilde{\psi}^{\tilde{k}}|_{I(\underline{\tilde{b}})})^{-1}(\widetilde{U})$, $(\tilde{\psi}^{\tilde{k}}|_{I(\underline{\tilde{b}})})^{-1}(B_{\widetilde{U}})$ without affecting \eqref{eq:bridges_gap_quotient_2}, so we can assume that $U_i$, $B_{U_i}$, $\widetilde{U}$, $B_{\widetilde{U}}$ have sufficiently small diameter. 

    Observe that we can refine the Markov partitions of $K$ and $\widetilde{K}$ and replace $\psi$ and $\tilde{\psi}$ with $\psi^N$ and $\tilde{\psi}^N$, respectively, with arbitrarily large $N$, and replace the iterates $\tilde{\ell},\widetilde{k}$ and $\ell_{p_i},k_i$, $i=0,\dots,l-1$ appropriately, without changing the points $\tilde{p},\tilde{q}$ and $p_i,q_{e_i}$, $i=0,\dots,l-1$, so that $(K,\widetilde{K})$ still satisfy the hypothesis of \Cref{lem:main_lemma}. As consequence of \Cref{rem:smaller_copies}, we can exhibit $g_1(K)$ and $g_2(\widetilde{K})$ as dynamically defined Cantor sets defined by the expanding maps $g_1\circ\psi\circ g_1^{-1}$ and $g_2\circ\tilde{\psi}\circ g_2^{-1}$, respectively. Since multipliers at periodic points are preserved by $C^1$ conjugation, we have that
    \begin{equation*}
        \frac{\log|(g_1\circ\psi^{k_{e_i}}\circ g_1^{-1})^\prime(g_1(q_{e_i}))|}{\log|(g_2\circ\tilde{\psi}^{\widetilde{k}}\circ g_2^{-1})^\prime(g_2(\tilde{q}))|}=\frac{\log|(\psi^{k_{e_i}})^\prime(q_{e_i})|}{\log|(\tilde{\psi}^{\widetilde{k}})^\prime(\tilde{q})|}\in\R\setminus\Q, \quad\text{for } i=0,\dots,l-1.
    \end{equation*}
    Moreover, by the continuity of $g_1^\prime$ and $g_2^\prime$, and the fact that the gaps and adjacent intervals are arbitrarily small, by the mean value theorem we will have that
    \begin{equation*}
            \frac{|g_1(B_{U_{i+1}})|}{|g_1(U_{i+1})|}\frac{|g_2(B_{\widetilde{U}})|}{|g_2(\widetilde{U})|}<1, \quad\text{for } i=0,\dots,l-1.
    \end{equation*}
    Therefore the dynamically defined Cantor sets $(g_1(K),g_2(\widetilde{K}))$ satisfy all the hypothesis of \Cref{lem:main_lemma}. Finally notice that if both $g_1$ and $g_2$ are decreasing, then we can apply the result to $-g_1(K)$ and $-g_2(\widetilde{K})$ to obtain uncountability of the boundary arbitrarily close to $\max(-g_1(K))+\max(-g_2(\widetilde{K}))=-(g_1(\min K)+g_2(\min\widetilde{K}))$. 
\end{proof}

\begin{example}
    From \Cref{example_1} we know that $\tau_p(C_3,\max C_3)<1$. It can be shown similarly that $\tau_p(C_2,\max C_2)<1$ and that the corresponding multipliers at $\max C_3$ and $\max C_2$ are multiplicatively independent (the required condition \eqref{eq:bridges_gap_quotient_2}). In particular we obtain that the boundary of $C_2+C_3$ has positive Hausdorff dimension near $\max(C_2+C_3)$.
\end{example}

\section{Uncountable boundary at $\nu_F$}\label{sec:uncountable_nu_F}

In this section we will give two proofs that the boundary near $\nu_F$ has positive Hausdorff dimension. Both proofs are based on the same ideas, one is by applying directly \Cref{lem:specialized_lemma} and the second one follows the proof of \Cref{lem:specialized_lemma} by using more explicit choices and thus giving a better estimate on the Hausdorff dimension. Through this section we will denote $\underline{a}\in\{1,2,3,4\}^\Z$ to be a bi-infinite sequence. Let us recall that Freiman's gap before Hall's ray is given by
\[c_{F} = m(\overline{121313}22344^{*}3211\overline{313121}),\]
\[\nu_{F} = m(\overline{323444}313134^{*}313121133\overline{313121}).\]

The following lemmas were proven in \cite[Section 4]{Newgaps}. 

\begin{lemma}
If $4.5278295<m(\underline{a})=\lambda_0(\underline{a})<4.5278296$, then either $m(\underline{a})\geq \mu_F>\nu_F$ or $m(\underline{a})\leq\nu_F$ and, up to transposition,  
$$\underline{a} = \dots 44313134^*313121133\dots$$
\end{lemma}

\begin{lemma}
    The set $M\cap(4.5278295, \nu_F]$ is contained in $K+\widetilde{K}$ where
    $$K=\{[4;3,1,3,1,2,1,1,3,3,3+x]: x\in K_1\},$$ 
    $$\widetilde{K}=\{[0;3,1,3,1,3,4,4,4+y]: y\in K_2\},$$ 
    $K_1=\{[0;1,3,1,2,1,\underline{\theta}]\in \widehat{K}\}$, $K_2=\{[0;3,2,3,\widetilde{\underline{\theta}}]\in \widehat{K}\}$ with 
    \begin{multline*}
        \widehat{K}=\{[0;\widehat{\underline{\theta}}]:\widehat{\underline{\theta}}\in\{1,2,3,4\}^{\mathbb{N}} \textrm{ doesn't contain the strings } 14, 24, 433, 434, \\ 131313, 2343, 223444, 123444, 231313, 331313 \textrm{ or their transposes}\}.
    \end{multline*}
\end{lemma}
Observe that $\widehat{K}$ is a Gauss-Cantor set, so it is a dynamically defined Cantor set given by iterates of the Gauss map (as in \Cref{example}) and $K_1$ and $K_2$ are conjugated to $\widehat{K}$ (up to refine the Markov partition of $\widehat{K}$, see \Cref{rem:smaller_copies}). We will show that the hypothesis of \Cref{lem:specialized_lemma} hold for $(K_1,K_2)$. Since $K=g_1(K_1)$ and $\widetilde{K}=g_2(K_2)$ for the increasing $C^1$ contractions
\begin{equation*}
    g_1(x)=[4;3,1,3,1,2,1,1,3,3,3+x] \quad\text{and}\quad g_2(y) = [0;3,1,3,1,3,4,4,4+y], 
\end{equation*}
by \Cref{cor:conjugating_preserves_lemma} we will obtain that $\partial(K+\widetilde{K})$ is uncountable arbitrarily close to $\max K+\max\widetilde{K} = \nu_F$.

Consider the following inverse branches $f_1$, $f_2$ of $K_1$ and $K_2$, respectively,
$$f_1(x)= [0;1,3,1,2,1,3+x] \quad \textrm{ and } \quad f_2(y)=[0;3,2,3,4,4,4+y]$$ 
whose fixed points are 
$$\alpha=[0;\overline{1,3,1,2,1,3}] \quad \textrm{ and } \quad \beta = [0;\overline{3,2,3,4,4,4}].$$ 

\begin{lemma} 
One has $\alpha=\max K_1$, $\beta=\max K_2$, and 
$$\frac{\log|f_1'(\alpha)|}{\log|f_2'(\beta)|}\in\mathbb{R}\setminus\mathbb{Q}.$$
\end{lemma} 

\begin{proof}
We note that $f_1$ and $f_2$ can be rewritten as 
    $$f_1(x)= \frac{15x + 56}{19x + 71}, \quad f_2(y)=\frac{127y + 538}{436y + 1847}.$$
    In particular, 
    $$f_1'(x)=\frac{1}{(19x + 71)^2}, \quad \quad f_2'(y) = \frac{1}{(436y + 1847)^2}$$ 
    and 
    $$\alpha =\frac{2\sqrt{462} - 28}{19}, \quad \quad \beta = \frac{\sqrt{243542} - 430}{218}.$$
    
   A simple computation yields
   \begin{equation*}
       f_1'(\alpha)=\frac{1}{(2\sqrt{462} + 43)^2} \quad \quad f_2'(\beta)=\frac{1}{(2\sqrt{243542} + 987)^2}.
   \end{equation*}
   Since $462=2 \cdot 3 \cdot 7 \cdot 11$ and $243542=2 \cdot 13 \cdot 17 \cdot 19 \cdot 29$, these numbers generate different quadratic extensions of $\Q$.    
\end{proof}

\begin{lemma}\label{lem:34}
Consider the intervals $B_{U_1}=[\alpha_1,\alpha]$, $U_1=(\alpha_2,\alpha_1)$, $B_{\widetilde{U}}=[\beta_1,\beta]$ and $\widetilde{U}=(\beta_2,\beta_1)$, where 
\begin{itemize}
\item $\alpha_1$ is the smallest element of $K_1$ of the form $[0;1,3,1,2,1,3,\underline{\theta}]$, 
\item $\alpha_2$ is the largest element of $K_1$ of the form $[0;1,3,1,2,1,2,\underline{\theta}]$, 
\item $\beta_1$ is the smallest element of $K_2$ of the form $[0;3,2,3,4,\widetilde{\underline{\theta}}]$,
\item $\beta_2$ is the largest element of $K_2$ of the form $[0;3,2,3,3,\widetilde{\underline{\theta}}]$.
\end{itemize}
Then, 
$$\frac{|B_{U_1}|}{|U_1|}<0.98 \quad \textrm{and} \quad \frac{|B_{\widetilde{U}}|}{|\widetilde{U}|}<0.0073.$$
Moreover we have that $G^\ell(\alpha_2)=\alpha$ for $\ell=8$ where $G$ is the Gauss map.
\end{lemma}

\begin{proof}
We have
\begin{align*}
    \alpha&=[0;\overline{1,3,1,2,1,3}], & \beta &= [0;\overline{3,2,3,4,4,4}]  \\
    \alpha_1 &= [0; 1, 3, 1, 2, 1, 3, 4, 3, 1, \overline{1, 3, 1, 3, 1, 2}], &\beta_1 &= [0; 3, 2, 3, 4, 4, 3, \overline{2, 1, 3, 1, 3, 1}] \\
    \alpha_2 &= [0; 1, 3, 1, 2, \overline{1, 2, 1, 3, 1, 3}], 
    & \beta_2 &= [0; 3, 2, 3, \overline{3, 1, 3, 1, 2, 1}] 
\end{align*}

\end{proof}

The previous lemma provides us with the ingredients necessary to apply \Cref{lem:specialized_lemma} by choosing $l = 1$, $q_{e_0} = q_{e_1} = \alpha$, $p_1 = \alpha_2$, $\tilde{q} = \beta$. As consequence, we obtain that $K+\widetilde{K}=g_1(K_1)+g_2(K_2)$ contains infinitely many gaps $(a,b)$ of the form $a=x+y$ where
\begin{multline*}
    x = [4;3,1,3,1,2,1,1,3,3,3,(1,3,1,2,1,3)^{n_1},1,3,1,2,1,2,1,3,\\
    (1,3,1,2,1,3)^{n_2},1,3,1,2,1,2,1,3,\dots]
\end{multline*}
and
\begin{equation*}
    y = [0;3,1,3,1,3,\overline{4,4,4,3,2,3}]
\end{equation*}

By definition, these gaps $(a,b)$ also do not intersect $M$ and hence also not $L$. Now we show that all these points $x+y$ belong to $L^\prime$ and $M^\prime$, from which will follow that they also belong to the boundary of both $L$ and $M$. 

Observe that
\begin{equation*}
    x+y = \lambda_0\left(\overline{323444}313134^{*}3131211333(131213)^{n_1}13121213\dotsb(131213)^{n_k}\dotsb\right).
\end{equation*}

Since we have the inequalities $\lambda_0(\dots 44^*31\dots)<4.52$, $\lambda_0(\dots3234^*443\dots)<4.5277$, $\lambda_0(\dots 4313^*12\dots)<4.51$, $\lambda_0(\dots 3313^*12\dots)<4.51$,  $\lambda_0(\dots 21313^*12\dots)<4.5278$, since $x+y>4.5278$ we have that $x+y\in M$.

Now we will show that for any choice of positive integers $(n_1,n_2,\dots)$ we also have $x+y\in L$. For this end, we will show that there is a sequence of periodic bi-infinite sequences $\overline{c^{(k)}}\in\{1,2,3,4\}^\Z$ with period $c^{(k)}$ such that their Markov (and hence Lagrange values) converge to $x+y$. Since $L$ is closed this will prove that $x+y\in L$. It suffices to consider
\begin{equation*}
    c^{(k)} = (323444)^k313134^{*}313121133(131213)^{n_1}13121213\dotsb(131213)^{n_k}2.
\end{equation*}
The above inequalities show that Markov value of $\overline{c^{(k)}}$ is attained at the marked position with $*$.

Now we will provide a more explicit proof that the boundary has positive Hausdorff dimension. This proof will allow us to give better explicit lower bounds for this dimension. Given a finite non empty word $\underline{c}=(c_1,\dots,c_n)\in\{1,2,3,4\}^n$ we write $[0;\underline{c}+x]=[0;c_1,\dots,c_n+x]$. Let us define the constants $\gamma<\gamma_{12}<\gamma_{11}<\dots<\gamma_1$ by
{
\allowdisplaybreaks[4]
\begin{align*}
    \gamma &= [0;\overline{3,1,3,1,2,1}], \\
    \gamma_1 &= [0;3,2,\overline{3,1,3,1,2,1}]=[0;\underline{p}^{(1)}+\gamma] \\
    \gamma_2 &= [0;3,1,\overline{1,3,1,3,1,2}], \\
    \gamma_3 &= [0;3,1,2,1,\overline{3,1,3,1,2,1}]=[0;\underline{p}^{(2)}+\gamma], \\
    \gamma_4 &= [0;3,1,3,4,3,1,\overline{1,3,1,3,1,2}], \\
    \gamma_5 &= [0;3,1,3,2,\overline{3,1,3,1,2,1}]=[0;\underline{p}^{(3)}+\gamma], \\
    \gamma_6 &= [0;3,1,3,1,\overline{1,3,1,3,1,2}], \\
    \gamma_7 &= [0;3,1,3,1,1,\overline{1,3,1,3,1,2}]=[0;\underline{p}^{(4)}+\gamma], \\
    \gamma_8 &= [0;3,1,3,1,2,\overline{3,1,3,1,2,1}], \\
    \gamma_9 &= [0;3,1,3,1,2,2,\overline{3,1,3,1,2,1}]=[0;\underline{p}^{(5)}+\gamma], \\
    \gamma_{10} &= [0;3,1,3,1,2,1,\overline{1,3,1,3,1,2}], \\
    \gamma_{11} &= [0;3,1,3,1,2,1,2,\overline{1,3,1,3,1,2}]=[0;\underline{p}^{(6)}+\gamma], \\
    \gamma_{12} &= [0;3,1,3,1,2,1,3,4,3,1,\overline{1,3,1,3,1,2}], 
\end{align*}
}

where the pre-periods all have even length, more explicitly $\underline{p}^{(1)}=(3,2)$, $\underline{p}^{(2)}=(3,1,2,1)$, $\underline{p}^{(3)}=(3,1,3,2)$, $\underline{p}^{(4)}=(3,1,3,1,1,1)$, $\underline{p}^{(5)}=(3,1,3,1,2,2)$ and $\underline{p}^{(6)}=(3,1,3,1,2,1,2,1)$. We extend this list to a strictly decreasing infinite sequence $\gamma<\dots<\gamma_{j+1}<\gamma_j<\dots<\gamma_1$ for $j\geq 13$ by writing 
\begin{equation*}
    \gamma_j = [0;(3,1,3,1,2,1)^q+\gamma_i], \qquad j=12q+i, \quad\text{where } i\in\{1,\dots,12\},
\end{equation*}
and we also extend cyclically $\underline{p}^{(j)}=\underline{p}^{(k)}$ where $j=6p+k$ with $k\in\{1,\dots,6\}$.

\begin{lemma}
    Let $\underline{c}=(c_1,\dots,c_{2n})$ be a non empty finite word over $\{1,2,3,4\}$ of finite even length which last digit equal to $1$ or $2$ which does not contain any of the forbidden words of $\widehat{K}$. Then
\begin{itemize}
\item $[0;\underline{c}+\gamma_1]$ is the smallest element of $\widehat{K}$ of the form $[0;\underline{c},3,2,\underline{\theta}]$, 
\item $[0;\underline{c}+\gamma_2]$ is the largest element of $\widehat{K}$ of the form $[0;\underline{c},3,1,\underline{\theta}]$, 
\item $[0;\underline{c}+\gamma_3]$ is the smallest element of $\widehat{K}$ of the form $[0;\underline{c},3,1,2,\underline{\theta}]$, 
\item $[0;\underline{c}+\gamma_4]$ is the largest element of $\widehat{K}$ of the form $[0;\underline{c},3,1,3,\underline{\theta}]$, 
\item $[0;\underline{c}+\gamma_5]$ is the smallest element of $\widehat{K}$ of the form $[0;\underline{c},3,1,3,2,\underline{\theta}]$, 
\item $[0;\underline{c}+\gamma_6]$ is the largest element of $\widehat{K}$ of the form $[0;\underline{c},3,1,3,1,\underline{\theta}]$, 
\item $[0;\underline{c}+\gamma_7]$ is the smallest element of $\widehat{K}$ of the form $[0;\underline{c},3,1,3,1,1,\underline{\theta}]$, 
\item $[0;\underline{c}+\gamma_8]$ is the largest element of $\widehat{K}$ of the form $[0;\underline{c},3,1,3,1,2,\underline{\theta}]$, 
\item $[0;\underline{c}+\gamma_9]$ is the smallest element of $\widehat{K}$ of the form $[0;\underline{c},3,1,3,1,2,2,\underline{\theta}]$, 
\item $[0;\underline{c}+\gamma_{10}]$ is the largest element of $\widehat{K}$ of the form $[0;\underline{c},3,1,3,1,2,1,\underline{\theta}]$, 
\item $[0;\underline{c}+\gamma_{11}]$ is the smallest element of $\widehat{K}$ of the form $[0;\underline{c},3,1,3,1,2,1,2,\underline{\theta}]$, 
\item $[0;\underline{c}+\gamma_{12}]$ is the largest element of $\widehat{K}$ of the form $[0;\underline{c},3,1,3,1,2,1,3,\underline{\theta}]$, 
\end{itemize}
In particular, the intervals $([0;\underline{c},\gamma_{2i}],[0;\underline{c},\gamma_{2i-1}])$ are gaps of $\widehat{K}$ for all $i\geq 1$.
\end{lemma}

Let us define the constants $\delta<\delta_2<\delta_1$ by
\begin{itemize}
\item $\delta$ is the smallest element of $\widehat{K}$ of the form $[0;4,4,4,3,2,3,\underline{\theta}]$, 
\item $\delta_1$ is the smallest element of $\widehat{K}$ of the form $[0;4, 4, 4, 3, 2, 3, 3,\underline{\theta}]$, 
\item $\delta_2$ is the largest element of $\widehat{K}$ of the form $[0;4, 4, 4, 3, 2, 3, 4,\underline{\theta}]$. 
\end{itemize}
\begin{align*}
    \delta &= [0;\overline{4,4,4,3,2,3}], \\
    \delta_1 &= [0; 4, 4, 4, 3, 2, 3, \overline{3, 1, 3, 1, 2, 1}], \\
    \delta_2 &= [0; 4, 4, 4, 3, 2, 3, 4, 4, 3, \overline{2, 1, 3, 1, 3, 1}].
\end{align*}
The interval $(\delta_2,\delta_1)$ is a gap of $\widehat{K}$. Let us denote the inverse branches
\begin{equation*}
    f_1(x)=[0;3,1,3,1,2,1+x] \quad\text{and}\quad f_2(x)=[0;4,4,4,3,2,3+x].
\end{equation*}
The constants $\gamma$ and $\delta$ are fixed points of $f_1$ and $f_2$ respectively. For convenience, let us note that the derivative of $f_1^\prime(\gamma)$ is much larger than $f_2^\prime(\delta)$
\begin{equation}
 f_1^\prime(\gamma)=\frac{1}{(2\sqrt{462} + 43)^{2}}\approx 10^{-4} \quad\text{and}\quad f_2^\prime(\delta)=\frac{1}{(2\sqrt{243542} + 987)^{2}}\approx 2\cdot10^{-7}.
\end{equation}
Let us denote the decreasing maps
\begin{equation*}
    g_1(x)=[4;3,1,3,1,2,1,1,3,3+x] \quad\text{and}\quad g_2(y) = [0;3,1,3,1,3+y].
\end{equation*}

Given $\underline{c}=(c_1,\dots,c_{2n})$ a non empty finite word over $\{1,2,3,4\}$ of finite even length which last digit equal to $1$ or $2$ which does not contain any of the forbidden words of $\widehat{K}$, we denote the following gaps of $K$ and their right adjacent intervals by 
\begin{equation*}
    U_i(\underline{c}) = \left(g_1([0;\underline{c}+\gamma_{2i-1}]),g_1([0;\underline{c}+\gamma_{2i}])\right) ~\text{and}~ D_i(\underline{c}) = \left[g_1([0;\underline{c},\gamma_{2i}]),g_1([0;\underline{c},\gamma])\right].
\end{equation*}
We also define left adjacent intervals by
\begin{equation*}
    C_i(\underline{c})=\left[g_1([0;\underline{c},\underline{p}^{(i)}+\gamma_{1}]),g_1([0;\underline{c}+\gamma_{2i-1}])\right].
\end{equation*}
The interval $C_i(\underline{c})$ is chosen in such a way so that it contains all the gaps $U_j(\underline{c},\underline{p}^{(i)})$ for all $j\geq 1$. Let us denote the following gaps and adjacent intervals of $\widetilde{K}$ by
\begin{equation*}
    \widetilde{U}_0 = (g_2(\delta_1),g_2(\delta_2)), \quad D_{\widetilde{U}_0} = [g_2(\delta_2),g_2(\delta)]
\end{equation*}
and
\begin{equation*}
    \widetilde{U}_m = g_2 \circ f_2^m(\widetilde{U}), \quad D_{\widetilde{U}_m} = g_2 \circ f_2^m(D_{\widetilde{U}}).
\end{equation*}

To estimate the bounded distortion explicitly, we use the following elementary lemma.

\begin{proposition}\label{prop:bounded_distortion}
    Let $f(x)=[0;a_1,\dots,a_k+x]$ be an inverse branch of the Gauss map. Then
    \begin{equation*}
        \left|\log\frac{|f^\prime(x)|}{|f^\prime(y)|}\right|\leq 2\left|\log\frac{x}{y}\right|.
    \end{equation*}
\end{proposition}

\begin{proof}
    Write
    \begin{equation*}
        f(x)=\frac{p_{k-1}x+p_k}{q_{k-1}x+q_k}, \quad |f^\prime(x)|=\frac{1}{(q_{k-1}x+q_k)^2}.
    \end{equation*}
    Hence
    \begin{align*}
        \left|\log\frac{|f^\prime(x)|}{|f^\prime(y)|}\right|&=2|\log(q_{k-1}y+q_k)-\log(q_{k-1}x+q_k)| \\
        &=2\left|\int_x^y\frac{q_{k-1}}{q_{k-1}t+q_k}dt\right|\leq2\left|\int_x^y\frac{1}{t}dt\right| = 2\left|\log\frac{x}{y}\right|.
    \end{align*}
\end{proof}

\begin{lemma}
    Given a gap $U_i(\underline{c})$ as above with $i\geq 1$, we have
    \begin{equation*}
        0.8<\frac{|D_i(\underline{c})|}{|U_i(\underline{c})|}<2.2.
    \end{equation*}
\end{lemma}

\begin{proof}
    If $i>6$, we can replace $\underline{c}$ with $(\underline{c},(3,1,3,1,2,1)^q)$ where $q=\lfloor i/6\rfloor$, so it suffices to consider $i=1,\dots,6$. If we let $h(x)=[4;3,1,3,1,2,1,1,3,3,c_1,\dots,c_{2n}+x]$, then by the mean value theorem we have that for certain $\gamma<y_0<\gamma_{2i}<x_0<\gamma_{2i-1}$ it holds that
    \begin{equation*}
        \frac{|D_i(\underline{c})|}{|U_i(\underline{c})|} = \frac{|h^\prime(y_0)|}{|h^\prime(x_0)|}\cdot\frac{|\gamma_{2i}-\gamma|}{|\gamma_{2i-1}-\gamma_{2i}|}.
    \end{equation*}
    In particular $|h^\prime(x_0)|<|h^\prime(y_0)|$ and for all $i=1,\dots,6$ we have by the proposition above \[|h^\prime(y_0)|/|h^\prime(x_0)|<e^{2\log(x_0/y_0)}<\exp^{2\log(\gamma_{2i-1}/\gamma)}<1.212.\]  
    A direct computation yields
    \begin{equation*}
        0.8<\frac{|\gamma_{2i}-\gamma|}{|\gamma_{2i-1}-\gamma_{2i}|}<1.81.
    \end{equation*}
\end{proof}

\begin{lemma}
    Given a gap $U_i(\underline{c})$ as above with $i\geq 1$, we have
    \begin{equation*}
        0.04156<\frac{|C_i(\underline{c})|}{|U_i(\underline{c})|}<0.0457.
    \end{equation*}
\end{lemma}

\begin{proof}
    If $i>6$, we can replace $\underline{c}$ with $(\underline{c},(3,1,3,1,2,1)^q)$ where $q=\lfloor i/6\rfloor$, so it suffices to consider $i=1,\dots,6$. If we let $h(x)=[4;3,1,3,1,2,1,1,3,3,c_1,\dots,c_{2n}+x]$, then by the mean value theorem we have that for certain $\gamma_{2i}<x_0<\gamma_{2i-1}<y_0<[0;\underline{p}^{(i)}+\gamma_{1}]$ it holds that
    \begin{equation*}
        \frac{|C_i(\underline{c})|}{|U_i(\underline{c})|} = \frac{|h^\prime(y_0)|}{|h^\prime(x_0)|}\cdot\frac{|[0;\underline{p}^{(i)}+\gamma_1]-\gamma_{2i-1}|}{|\gamma_{2i-1}-\gamma_{2i}|}.
    \end{equation*}
    In particular $|h^\prime(y_0)|<|h^\prime(x_0)|$ and for all $i=1,\dots,6$ we have by the proposition above \[|h^\prime(y_0)|/|h^\prime(x_0)|>e^{2\log(-y_0/x_0)}<\exp^{-2\log([0;\underline{p}^{(i)}+\gamma_{1}]/\gamma_{2i})}>0.997.\]  
    A direct computation yields
    \begin{equation*}
        0.04169<\frac{|\gamma_{2i}-\gamma|}{|\gamma_{2i-1}-\gamma_{2i}|}<0.0457.
    \end{equation*}
\end{proof}

\begin{lemma}
    Given a gap $U_i(\underline{c})$ as above with $i\geq 1$, we have 
    \begin{equation}\label{eq:big_gap_to_first_gap}
        0.01451<\frac{|U_1(\underline{c},\underline{p}^{(i)})|}{|U_i(\underline{c})|}<0.01624.
    \end{equation}
\end{lemma}

\begin{proof}
    If $i>6$, we can replace $\underline{c}$ with $(\underline{c},(3,1,3,1,2,1)^q)$ where $q=\lfloor i/6\rfloor$, so it suffices to consider $i=1,\dots,6$. If we let \[h(x)=[4;3,1,3,1,2,1,1,3,3,c_1,\dots,c_{2n}+x],\] 
     then by the mean value theorem we have that for certain $\gamma_{2i}<x_0<\gamma_{2i-1}<[0;\underline{p}^{(i)}+\gamma_{2}]<y_0<[0;\underline{p}^{(i)}+\gamma_{1}]$ it holds that
     \begin{equation*}
         \frac{|U_j(\underline{c},\underline{p}^{(i)})|}{|U_i(\underline{c})|} = \frac{|h^\prime(y_0)|}{|h^\prime(x_0)|} \cdot \frac{|[0;\underline{p}^{(i)}+\gamma_{2j}]-[0;\underline{p}^{(i)}+\gamma_{2j-1}]|}{|\gamma_{2i}-\gamma_{2i-1}|}
     \end{equation*}
    Clearly we have $|h^\prime(y_0)|<|h^\prime(x_0)|$ and by the proposition above \[|h^\prime(y_0)|/|h^\prime(x_0)|>e^{-2\log(y_0/x_0)}>\exp^{-2\log([0;\underline{p}^{(i)}+\gamma_{1}]/\gamma_{2i})}>0.997.\]  
    A direct computation for all $i=1,\dots,6$ yields the estimate
    \begin{equation*}
        0.01456<\frac{|[0;\underline{p}^{(i)}+\gamma_{2}]-[0;\underline{p}^{(i)}+\gamma_{1}]|}{|\gamma_{2i}-\gamma_{2i-1}|}<0.01624.
    \end{equation*}
\end{proof}

\begin{lemma}
    Given an admissible word $\underline{c}$ as above, we have for $i\geq 1$ that
    \begin{equation}\label{eq:gap_j_to_gap_j_plus_1}
        0.1796<\frac{|U_{j+1}(\underline{c})|}{|U_j(\underline{c})|}<0.4945.
    \end{equation}
\end{lemma}

\begin{proof}
    Writing $j=6q+i$ where $1\leq i\leq 5$, we let
    \[h(x)=[4;3,1,3,1,2,1,1,3,3,c_1,\dots,c_{2n},(3,1,3,1,2,1)^q+x].\] 
    By the mean value theorem we have that for certain $\gamma_{2i+2}<y_0<\gamma_{2i+1}<\gamma_{2i}<x_0<\gamma_{2i-1}$ it holds that
    \begin{equation*}
        \frac{|U_{i+1}(\underline{c})|}{|U_i(\underline{c})|}=\frac{|h^\prime(y_0)|}{|h^\prime(x_0)|}\cdot\frac{|\gamma_{2i+1}-\gamma_{2i}|}{|\gamma_{2i-1}-\gamma_{2i}|}.
    \end{equation*}
    Clearly we have $|h^\prime(x_0)|<|h^\prime(y_0)|$ and by the proposition above \[|h^\prime(y_0)|/|h^\prime(x_0)|<e^{2\log(x_0/y_0)}<e^{2\log(\gamma_{2i-1}/\gamma_{2i+2})}<1.2.\]
    A direct computation yields for all $1\leq i\leq 5$
    \begin{equation*}
        0.1796<\frac{|\gamma_{2i+1}-\gamma_{2i}|}{|\gamma_{2i-1}-\gamma_{2i}|}<0.4121.
    \end{equation*}
\end{proof}

Similarly using the proposition above
\begin{equation*}
    \frac{|\delta_2-\delta|}{|\delta_1-\delta_2|}<\frac{|D_{\widetilde{U}_m}|}{|\widetilde{U}_m|}< e^{2\log(\delta_1/\delta)}\frac{|\delta_2-\delta|}{|\delta_1-\delta_2|},
\end{equation*}
hence
\begin{equation}\label{eq:relation_between_D_tilde_and_U_tilde}
    0.00722<\frac{|D_{\widetilde{U}_m}|}{|\widetilde{U}_m|}<1.00000007\cdot0.00723 < 0.00724.
\end{equation}
Analogously for $m\geq 1$
\begin{equation*}
    \frac{|(f_2)^{m}(\delta_1)-f_2^{m}(\delta_2)|}{|\delta_1-\delta_2|}<\frac{|\widetilde{U}_{m}|}{|\widetilde{U}_{0}|} < e^{2\log(\delta_1/\delta)}\frac{|(f_2)^m(\delta_1)-(f_2)^m(\delta_2)|}{|\delta_1-\delta_2|},
\end{equation*}
hence, after setting $\rho=f_2^\prime(\delta)$, we have
\begin{equation}\label{eq:estimate_rho}
    \rho^m<\frac{|\widetilde{U}_{m}|}{|\widetilde{U}_{0}|}<1.00000007\cdot\rho^m,
\end{equation}
and
\begin{equation}\label{eq:gaps_widetilde_m_to_m_plues_1}
    2.57\cdot10^{-7}<\rho<\frac{|\widetilde{U}_{m+1}|}{|\widetilde{U}_{m}|}<1.00000007\cdot\rho<2.58\cdot10^{-7}.
\end{equation}

We initialize our procedure with $\underline{c}=(3,1,3,1,2,1)$ and $m=0$. Our first task is to choose among the sequence $U_i(3,1,3,1,2,1)$, $i\geq 1$ the smallest gap such that $|U_i(3,1,3,1,2,1)|>|\widetilde{D}_0|$. We have that $|\widetilde{D}_0|\approx 7.606\cdot10^{-15}$ and
\begin{equation*}
    |U_4(3,1,3,1,2,1)|\approx 8.956\cdot10^{-15} \quad\text{and}\quad |U_5(3,1,3,1,2,1)|\approx 1.609\cdot10^{-15}.
\end{equation*}
Hence our first choice is the gap $U_4(3,1,3,1,2,1)$. We note that all gaps of $K$ we are going to construct belong to the interval $[\min C_4(3,1,3,1,2,1), \max K]$.

Suppose that $U_i(\underline{c})$, $i\geq1$ and $m\geq 1$ are such that 
\begin{equation}\label{eq:threshold_equation}
    |\widetilde{D}_m|<|U_i(\underline{c})|<42\cdot|\widetilde{D}_m|.
\end{equation}
Then
\begin{equation*}
    \frac{|D_i(\underline{c})|}{|\widetilde{U}_m|}=\frac{|D_i(\underline{c})|}{|U_i(\underline{c})|}\cdot\frac{|U_i(\underline{c})|}{|\widetilde{D}_m|}\cdot\frac{|\widetilde{D}_m|}{|\widetilde{U}_m|}<2.2\cdot 42 \cdot0.00724 < 1.
\end{equation*}

Moreover we also have
\begin{align*}
    \frac{|C_i(\underline{c})|+|U_i(\underline{c})|+|D_i(\underline{c})|}{|\widetilde{D}_m|+|\widetilde{U}_m|} &= \frac{|U_i(\underline{c})|}{|\widetilde{D}_m|}\cdot\frac{1+\frac{|C_i(\underline{c})|+|D_i(\underline{c})|}{|U_i(\underline{c})|}}{1+\frac{|\widetilde{U}_m|}{|\widetilde{D}_m|}}\\
    &<\frac{1+0.05+2.2}{1+\frac{1}{0.00722}}<42\cdot 0.0233<1.
\end{align*}

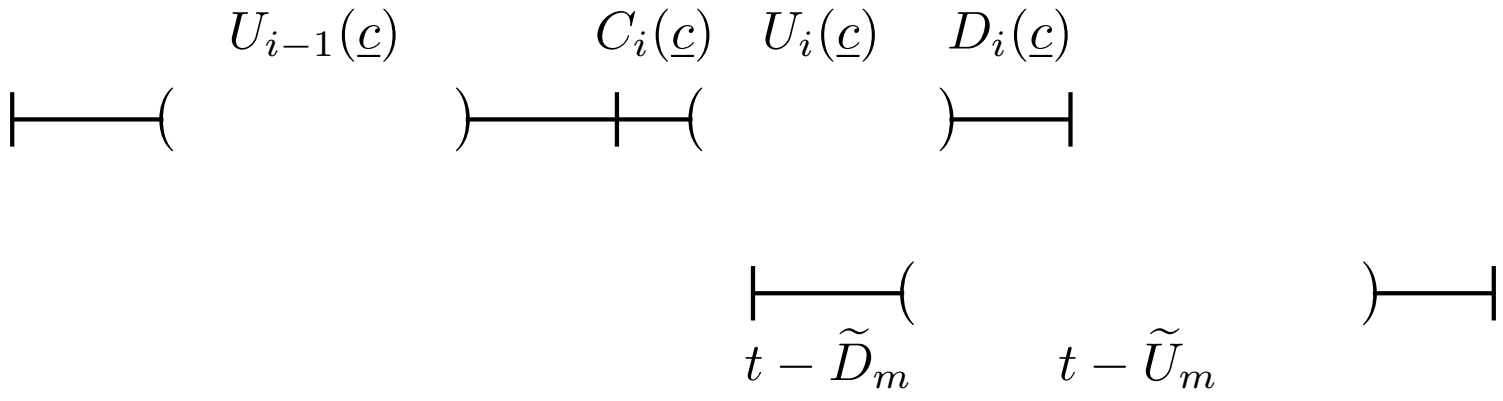
\begin{figure}
    \centering
	\begin{tikzpicture}[line width=0.8pt]

    \draw (-5,0) -- (-4,0);
    \draw (-2,0) -- (-0.5,0);
    \draw (-5,0.18) -- (-5,-0.18);
    \draw (-1,0.18) -- (-1,-0.18);
    
    \node at (-0.5,0) {\big (};
    \node at (1.2,0) {\big )};
    \node at (-4,0) {\big (};
    \node at (-2,0) {\big )};
    
    \draw (1.2,0) -- (2,0);
    \draw (2,0.18) -- (2,-0.18);
    
    \node[above] at (-3,0.25) {$U_{i-1}(\underline{c})$};
    \node[above] at (-0.75,0.25) {$C_{i}(\underline{c})$};
    \node[above] at (0.35,0.25) {$U_i(\underline{c})$};
    \node[above] at (1.6,0.25) {$D_i(\underline{c})$};
    
    \draw (-0.1,-1.15) -- (0.9,-1.15);
    \draw (-0.1,-1.33) -- (-0.1,-0.97);
    
    \node at (0.9,-1.15) {\big (};
    \node at (4,-1.15) {\big )};
    
    \draw (4,-1.15) -- (4.8,-1.15);
    \draw (4.8,-1.33) -- (4.8,-0.97);
    
    \node[below] at (0.4,-1.25) {$t-\widetilde{D}_m$};
    \node[below] at (2.45,-1.25) {$t-\widetilde{U}_m$};
    
    \end{tikzpicture}
    \caption{The exact moment $t-\widetilde{D}_m$ touches the left border of $U_i(\underline{c})$, the right border of $D_i(\underline{c})$ is far from the boundary of $t-\widetilde{U}_m$, provided that $|U_i(\underline{c})|+|D_i(\underline{c})|<|\widetilde{D}_m|+|\widetilde{U}_m|$. Moreover for that parameter $t = \min\partial U_i(\underline{c})+\max\widetilde{D}_m$ the intersection occurs at exactly one point. Moreover, the left border of $t-\widetilde{D}_m$ touches the left border of $C_i(\underline{c})$ much earlier than the right border of $t-\widetilde{U}_m$ touches the right border of $D_i(\underline{c})$ because $|C_i(\underline{c})|+|U_i(\underline{c})|+|D_i(\underline{c})|<|\widetilde{D}_m|+|\widetilde{U}_m|$.}
    \end{figure}

    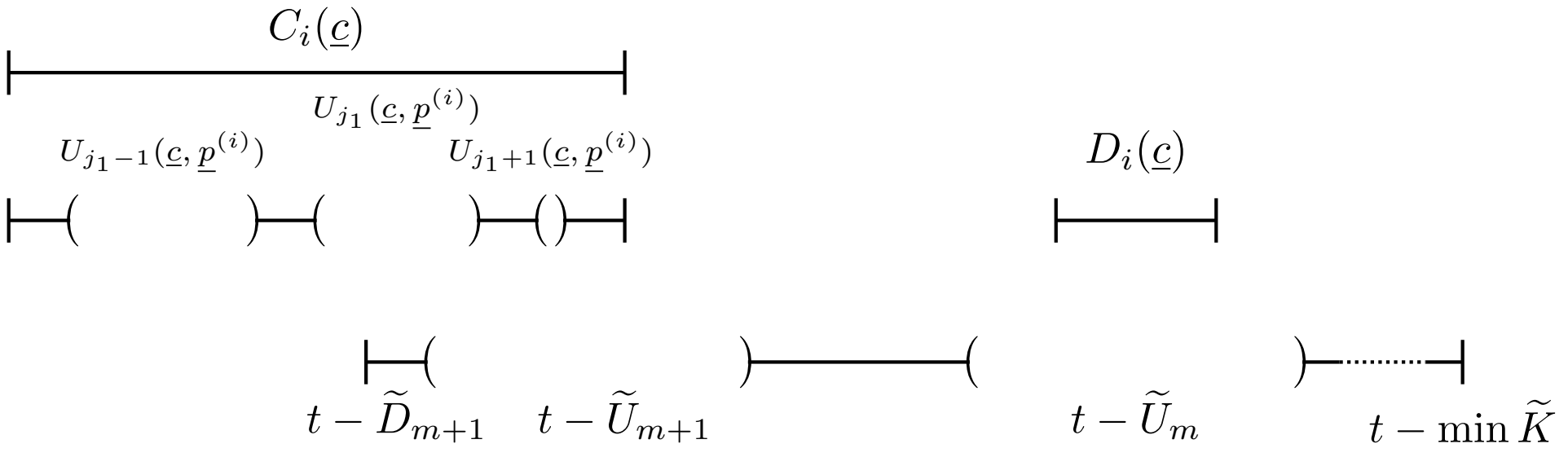
\begin{figure}
        \centering
        \begin{tikzpicture}[line width=0.8pt]

    \draw (-3,1.2) -- (2, 1.2);
    \draw (-3,0) -- (-2.5,0);
    \draw (-1,0) -- (-0.5,0);
    \draw (-3,0.18) -- (-3,-0.18);
    \draw (-3,1.38) -- (-3,1.02);
    \draw (2,1.38) -- (2,1.02);

    \draw (0.8,0) -- (1.3,0);
    \draw (1.5,0) -- (2,0);
    \draw (2,0.18) -- (2,-0.18);

    \node at (-2.5,0) {\big (};
    \node at (-1,0) {\big )};
    \node at (-0.5,0) {\big (};
    \node at (0.8,0) {\big )};
    \node at (1.3,0) {\big (};
    \node at (1.5,0) {\big )};

    \draw (5.5,0.18) -- (5.5,-0.18);
    \draw (5.5,0) -- (6.8,0);
    \draw (6.8,0.18) -- (6.8,-0.18);
    
    \node[above] at (-0.5,1.25) {$C_i(\underline{c})$};
    \node[above] at (-1.75,0.25) {\scriptsize $U_{j_1-1}(\underline{c},\underline{p}^{(i)})$};
    \node[above] at (0.15,0.6) {\scriptsize $U_{j_1}(\underline{c},\underline{p}^{(i)})$};
    \node[above] at (1.4,0.25) {\scriptsize $U_{j_1+1}(\underline{c},\underline{p}^{(i)})$};
    \node[above] at (6.15,0.25) {$D_i(\underline{c})$};

    \draw (-0.1,-1.15) -- (0.4,-1.15);
    \draw (-0.1,-1.33) -- (-0.1,-0.97);
    
    \node at (0.4,-1.15) {\big (};
    \node at (3,-1.15) {\big )};
    
    \draw (3,-1.15) -- (4.8,-1.15);
    \node at (4.8,-1.15) {\big (};

    \node at (7.5,-1.15) {\big )};
    \draw (7.5,-1.15) -- (7.8,-1.15);
    \draw[thick, densely dotted] (7.8,-1.15) -- (8.5,-1.15);
    \draw (8.5,-1.15) -- (8.8,-1.15);
    \draw (8.8,-1.33) -- (8.8,-0.97);
    
    \node[below] at (0.15,-1.25) {$t-\widetilde{D}_{m+1}$};
    \node[below] at (2,-1.25) {$t-\widetilde{U}_{m+1}$};
    \node[below] at (6.15,-1.25) {$t-\widetilde{U}_{m}$};

    \node[below,yshift=-0.5em] at (8.8,-1.15) {$t-\min\widetilde{K}$};
    
    \end{tikzpicture}
    \caption{We can push $t-\widetilde{K}$ to the left (i.e., decrease $t$), so that $t-\widetilde{D}_{m+1}$ enters into $C_i(\underline{c})$ to create new gaps. The scales of the picture do not correspond to the actual values (since $C_i(\underline{c})$ is much smaller than $D_i(\underline{c})$).}    
    \end{figure}

In particular, this creates a gap in the arithmetic sum $K+\widetilde{K}$. Indeed, $K+\widetilde{K}$ will have a gap from $\inf U_i(\underline{c})+\max\widetilde{D}_m =\inf U_i(\underline{c})+\max\widetilde{K}$ to $\min D_i(\underline{c})+\min\widetilde{D}_m$. Each time we create such gaps, we create two more gaps using the following rule. On the one hand, we will chose the smallest gap among the sequence $U_j(\underline{c})$, $j> i$ that satisfies
\begin{equation*}
    |\widetilde{D}_{m+1}|<|U_j(\underline{c})|<42\cdot|\widetilde{D}_{m+1}|,
\end{equation*}
and we will chose the smallest gap among the sequence $U_j(\underline{c},\underline{p}^{(i)})$, $j\geq 1$ that satisfies
\begin{equation*}
    |\widetilde{D}_{m+1}|<|U_j(\underline{c},\underline{p}^{(i)})|<42\cdot|\widetilde{D}_{m+1}|.
\end{equation*}

To see that such gaps exists, first notice that by \eqref{eq:gaps_widetilde_m_to_m_plues_1} we have $|\widetilde{D}_{m+1}|\approx 2.57\cdot10^{-7}\cdot|\widetilde{D}_m|$. Hence in the first case we can write
\begin{equation*}
    |U_j(\underline{c})|=\frac{|U_j(\underline{c})|}{|U_{j-1}(\underline{c})|}\cdot\dots\cdot\frac{|U_{i+1}(\underline{c})|}{|U_{i}(\underline{c})|}\cdot|U_i(\underline{c})|.
\end{equation*}
We can control each one of these products by \eqref{eq:gap_j_to_gap_j_plus_1}, so each one of these products lies between 0.17 and 0.5. In particular, since $|U_i(\underline{c})|>|\widetilde{D}_m|$, we can chose an $i+24>j>i$ such that $|U_j(\underline{c})|$ lies in the interval $(2.58\cdot10^{-7}\cdot|\widetilde{U}_m|,42\cdot2.57\cdot10^{-7}\cdot|\widetilde{U}_m|)$. Similarly, in the other situation we can write 
\begin{equation*}
    |U_j(\underline{c},\underline{p}^{(i)})|=\frac{|U_j(\underline{c},\underline{p}^{(i)})|}{|U_{j-1}(\underline{c},\underline{p}^{(i)})|}\cdot\dots\cdot\frac{|U_{2}(\underline{c},\underline{p}^{(i)})|}{|U_{1}(\underline{c},\underline{p}^{(i)})|}\cdot\frac{|U_1(\underline{c},\underline{p}^{(i)})|}{|U_i(\underline{c})|}\cdot|U_i(\underline{c})|.
\end{equation*}
Similarly, we can control the first products by \eqref{eq:gap_j_to_gap_j_plus_1} and the last product by \eqref{eq:big_gap_to_first_gap}, in order to guarantee that for some choice $j\geq 1$ the size of the gap $|U_j(\underline{c},\underline{p}^{(i)})|$ lies in  $(2.58\cdot10^{-7}\cdot|\widetilde{U}_m|,42\cdot2.57\cdot10^{-7}\cdot|\widetilde{U}_m|)$.

Finally, we will estimate explicitly the Hausdorff dimension of the previous construction using the mass distribution principle (\Cref{prop:mass_distribution}). Recall that we begun with $\underline{c}=(3,1,3,1,2,1)$ and the first gap we took was $U_{i_0}(\underline{c})$ with $i_0=4$. We distribute equal mass $1/2$ to each interval $C_{i_0}(\underline{c})$ and $D_{i_0}(\underline{c})$. Notice that each of these two interval have length at least $0.4\cdot|U_{i_0}(3,1,3,1,2,1)|$ and at most $2.2\cdot|U_{i_0}(3,1,3,1,2,1)|$, hence by \eqref{eq:threshold_equation} length at least $0.4\cdot|\widetilde{D}_0|$ and length most $93\cdot|\widetilde{D}_0|$. These two intervals are separated by a distance at least $|U_{i_0}(\underline{c})|>|\widetilde{D}_0|$. In the next step, we are going to chose two different gaps. The first gap $U_{i_1}(\underline{c})$ where $i_1>i_0$, lies inside of $D_{i_1}(\underline{c})$, so we distribute the mass of $D_{i_0}(\underline{c})$ equally into $C_{i_1}(\underline{c})$ and $D_{i_1}(\underline{c})$. On the other hand, we also chose a gap $U_{j_1}(\underline{c},\underline{p}^{(i_0)})$ where $j_1\geq 1$, which lies inside $C_{i_0}(\underline{c})$, so we distribute the mass of $C_{i_0}(\underline{c})$ equally into $C_{j_1}(\underline{c},\underline{p}^{(i_0)})$ and $D_{j_1}(\underline{c},\underline{p}^{(i_0)})$. We have now $2^2$ disjoint intervals all with sizes at least $0.4\cdot|\widetilde{D}_1|$ and at most $93\cdot|\widetilde{D}_1|$. Each pair of intervals is separated by a distance at least $|\widetilde{D}_1|$. Inductively, at stage $m$, we will have $2^{m}$ intervals with sizes at least $0.4\cdot|\widetilde{D}_m|>0.002|\widetilde{U}_0|\rho^m$ and at most $93\cdot|\widetilde{D}_m|<0.7|\widetilde{U}_0|\rho^m$, and which are separated by at least $|\widetilde{D}_m|>0.007|\widetilde{U}_0|\rho^m$, where we have used \eqref{eq:relation_between_D_tilde_and_U_tilde} and \eqref{eq:estimate_rho}. The intersection of all such sets, is a Cantor set $C$, which by the mass distribution principle has Hausdorff dimension at least $\log 2/\log(\rho^{-1})$ where $\rho=f_2^\prime(\delta)$. Since we have that $C+\max\widetilde{K}$ corresponds to all left boundaries of the gaps of $K+\widetilde{K}$ in our construction and such points belong to $\partial L$ because of the argument given after \Cref{lem:34}, we finish the proof of \Cref{thm:uncountable_nu_F}.

\section{Uncountable boundary near 4.5213} \label{sec:uncountable_near_45213}

Let us denote $\underline{b}\in\{1,2,3,4\}^\Z$ a bi-infinite sequence. Given a finite string $u=(a_i)_{i=-m}^{n}$, through this section we denote 
\begin{equation}
    \lambda^{-}_i(u):=\min \{[a_i;a_{i+1},...,a_n,\theta_1]+[a_{i-1};a_{i-2},...,a_{-m},\theta_2]: \theta_1,\theta_2\in \{1,2,3,4\}^{\mathbb{N}}\},
\end{equation}
and
\begin{equation}
    \lambda^{+}_i(u):=\max \{[a_i;a_{i+1},...,a_n,\theta_1]+[a_{i-1};a_{i-2},...,a_{-m},\theta_2];\theta_1,\theta_2\in \{1,2,3,4\}^{\mathbb{N}}\}.
\end{equation}

We will first show that $M$ has a gap from
\begin{equation*}
    \nu=\lambda_0\left(\overline{131312131311}34431213^*1312312\overline{21313121}\right)=4.5212920725004\dots
\end{equation*}
to
\begin{align*}
    \mu&=\lambda_0(\overline{131312131311}344431213^*13123123\overline{313121313111}) =4.521292947013\dots
\end{align*}

In the sequel, we will establish the following result:

\begin{proposition}\label{prop:gap_near_45213}
    One has $M\cap(\nu,\mu)=\varnothing$.
\end{proposition}

Recall that $\underline{b}\in\{1,2,3,4\}^\Z$ with $m(\underline{b})=\lambda_0(\underline{b})$ either satisfies $b_0^*=3^*$ or $b_0^*=4^*$. 

First, let us rule out the possibility that $b_0^*=4^*$.

\begin{remark} In what follows, a string leading to a Markov value $>4.5212924$ will be called forbidden and it will be marked in red. 
\end{remark}

\begin{lemma}
    $\lambda_0^{-}(\textcolor{red}{14^*})>4.75$, $\lambda_0^{-}(\textcolor{red}{24^*})>4.56$, $\lambda_0^{-}(\textcolor{red}{34^*3})>4.522$, $\lambda_0^{+}(44^*4)<4.8$.
\end{lemma}

Hence without loss of generality we can assume $b_{-2}b_{-1}b_0^*b_1\in\{444^*3,344^*3\}$. 

\begin{lemma}
    $\lambda_0^{-}(\textcolor{red}{444^*34})>\lambda_0^{-}(\textcolor{red}{444^*33})>\lambda_0^{-}(\textcolor{red}{444^*32})>4.525$, $\lambda_0^{-}(\textcolor{red}{344^*34})>\lambda_0^{-}(\textcolor{red}{344^*33})>\lambda_0^{-}(\textcolor{red}{344^*32})>4.5215$, $\lambda_0^{+}(44^*31)<4.52$. 
\end{lemma}

Therefore $b_0^*=4^*$ is not possible. Let us assume from now on that $b_0^*=3^*$. 

\begin{lemma}
    $\lambda_0^{+}(3^*4)<\lambda_0^{+}(3^*3)<\lambda_0^{+}(3^*2)<4.3$.
\end{lemma}

Hence we must have that $b_{-1}b_0^*b_1=13^*1$. 

\begin{lemma}
    $\lambda_0^{-}(\textcolor{red}{313^*13})>4.524$, $\lambda_0^{+}(213^*12)<4.48$, $\lambda_0^{+}(13^*11)<4.48$.
\end{lemma}

Hence, since 14 and 41 are forbidden, we can assume, that up to transposition, we have $b_{-2}b_{-1}b_0^*b_1b_2=213^*13$. 

\begin{lemma}
\begin{itemize}
    \item $\lambda_0^{+}(213^*134)<\lambda_0^{+}(213^*133)<\lambda_0^{+}(213^*132)<4.515$
    \item $\lambda_0^{+}(3213^*13)<4.5$, $\lambda_0^{+}(2213^*13)<4.51$.
\end{itemize}
\end{lemma}

Hence since 24 is forbidden we can assume that $b_{-3}b_{-2}b_{-1}b_0^*b_1b_2b_3=1213^*131$. 

\begin{lemma}
    $\lambda_0^{-}(\textcolor{red}{313^*13})>4.524$, $\lambda_0^{+}(1213^*1311)<4.5236$, $\lambda_0^{+}(11213^*131)<4.519$.
\end{lemma}

Since 14 and 41 are forbidden, we have that 
\begin{equation*}
    b_{-4}b_{-3}b_{-2}b_{-1}b_0^*b_1b_2b_3b_4\in\{21213^*1312,31213^*1312\}.
\end{equation*}

\begin{lemma} \hfill
\begin{itemize}
    \item $\lambda_0^{-}(\textcolor{red}{31213^*13121})>4.522$, $\lambda_0^{-}(\textcolor{red}{31213^*13122})>4.5215$, 
    \item $\lambda_0^{-}(\textcolor{red}{131213^*1312})>4.522$, $\lambda_0^{-}(\textcolor{red}{231213^*1312})>4.5214$,
    \item $\lambda_0^{-}(\textcolor{red}{331213^*13123})>4.5213$, $\lambda_0^{+}(21213^*13123)<4.52$,
    \item $\lambda_0^{+}(21213^*13122)<4.521$, $\lambda_0^{+}(221213^*1312)<4.52$, 
    \item $\lambda_0^{+}(321213^*1312)<4.519$. 
\end{itemize}
      
\end{lemma}

Since 24 and 42 are forbidden, we have that 
\begin{equation*}
    b_{-5}b_{-4}b_{-3}b_{-2}b_{-1}b_0^*b_1b_2b_3b_4b_5\in\{121213^*13121,431213^*13123\}.
\end{equation*}

\begin{lemma}
\begin{itemize}
    \item $\lambda_0^{+}(121213^*131211)<\lambda_0^{+}(121213^*131212)<4.52124$, 
    \item $\lambda_0^{+}(1121213^*13121)<\lambda_0^{+}(2121213^*13121)<4.5211$, 
    \item $\lambda_0^{+}(3121213^*131213)<4.52123$.
\end{itemize}
    
\end{lemma}

In particular we must have that 
\begin{equation*}
    b_{-5}b_{-4}b_{-3}b_{-2}b_{-1}b_0^*b_1b_2b_3b_4b_5=431213^*13123.
\end{equation*}

Since 14, 24, 343 are forbidden we must have $b_{-6}=4$. 

\begin{lemma}
\begin{itemize}
    \item $\lambda_0^{-}(\textcolor{red}{431213^*131233})>4.5213$, 
    \item $\lambda_0^{-}(\textcolor{red}{431213^*131234})>4.5213$,
    \item $\lambda_0^{-}(\textcolor{red}{4431213^*131232})>4.5213$.
\end{itemize}    
\end{lemma}

In particular we must have that 
\begin{equation*}
    b_{-6}\dots b_0^*\dots b_6=4431213^*131231.
\end{equation*}

\begin{lemma}
    $\lambda_0^{-}(\textcolor{red}{4431213^*1312311})>4.5213$, $\lambda_0^{-}(\textcolor{red}{4431213^*1312313})<4.52128$.
\end{lemma}

Since 14 and 24 are forbidden, the previous lemma implies that 
\begin{equation*}
    b_{-7}\dots b_0^*\dots b_7=\{34431213^*1312312,44431213^*1312312\}.
\end{equation*}

\begin{lemma}
\begin{itemize}
    \item $\lambda_0^{+}(44431213^*13123121)<4.52129$,
    \item $\lambda_0^{+}(44431213^*13123122)<4.5212917$.
\end{itemize}
\end{lemma}

Thus if $b_{-7}\dots b_0^*\dots b_7=44431213^*1312312$ then \[b_{-7}\dots b_0^*\dots b_7=44431213^*13123121.\] Therefore, using the forbidden words 14, 41, 42, 432, 433, 434, 31313, 1312131312, we obtain that $m(\underline{b})\geq \mu$.

\begin{lemma}
    $\lambda_0^{-}(\textcolor{red}{34431213^*13123123})>4.521293$, $\lambda_0^{+}(34431213^*13123121)<4.52129$
\end{lemma}

On the other hand, if $b_{-7}\dots b_0^*\dots b_7 = 34431213^*1312312$, the previous lemma and the fact that 432, 433, 434 are forbidden implies that \[b_{8}\dots b_0^*\dots b_8 = 134431213^*13123122.\] Therefore, using the forbidden words 14, 31313, 1312131312, 1213131213, 2213131213 we obtain that $m(\underline{b})\leq \nu$. This finishes the proof of \Cref{prop:gap_near_45213}.

A consequence of the previous analysis is the following local uniqueness results.

\begin{lemma}
If $4.5212917<m(\underline{a})=\lambda_0(\underline{a})<4.5212924$, then either $m(\underline{a})\geq \mu>\nu$ or $m(\underline{a})\leq\nu$ and, up to transposition,  
$$\underline{a} = \dots 134431213^*13123122\dots$$\end{lemma}

\begin{lemma}
    The set $M\cap(4.5212917, \nu]$ is contained in $K+\widetilde{K}$ where
    $$K=\{[3;1,3,1,2,3,1,2,2,1,3,1,3,1,2+x]: x\in K_1\},$$ 
    $$\widetilde{K}=\{[0;1,2,1,3,4,4,3,1+y]: y\in K_2\},$$ 
    and $K_1=\{[0;1,2,1,\underline{\theta}]\in \widehat{K}\}$, $K_2=\{[0;1,3,1,\widetilde{\underline{\theta}}]\in \widehat{K}\}$ with 
    \begin{multline*}
        \widehat{K}=\{[0;\widehat{\underline{\theta}}]:\widehat{\underline{\theta}}\in\{1,2,3,4\}^{\mathbb{N}} \textrm{ doesn't contain the strings } 14, 24, 432, 433, 343, 31313, \\ 1312131312, 1213131213 \textrm{ or their transposes}\}.
    \end{multline*}
\end{lemma}

Observe that $\widehat{K}$ is a Gauss-Cantor set, so it is a dynamically defined Cantor set given by iterates of the Gauss map (as in \Cref{example}) and $K_1$ and $K_2$ are conjugated to $\widehat{K}$ (up to refine the Markov partition of $\widehat{K}$, see \Cref{rem:smaller_copies}). We will show that the hypothesis of \Cref{lem:specialized_lemma} hold for $(K_1,K_2)$. Since $K=g_1(K_1)$ and $\widetilde{K}=g_2(K_2)$ for the $C^1$ contractions
\begin{equation*}
    g_1(x)=[3;1,3,1,2,3,1,2,2,1,3,1,3,1,2+x] \quad\text{and}\quad g_2(y) = [0;1,2,1,3,4,4,3,1+y], 
\end{equation*}
by \Cref{cor:conjugating_preserves_lemma} we will obtain that $\partial(K+\widetilde{K})$ is uncountable arbitrarily close to $\max K+\max\widetilde{K} = \nu$.

Consider the following inverse branches $f_1$, $f_2$ of $K_1$ and $K_2$, respectively,
\begin{equation*}
    f_1(x)=[0;1,2,1,3,1,3,1,2+x] \quad\text{and}\quad f_2(y)=[0;1,3,1,3,1,2,1,3,1,3,1,1+y].
\end{equation*}

The fixed points of $f_1$ and $f_2$ are 
$$\alpha=[0;\overline{1,2,1,3,1,3,1,2}] \quad \text{and} \quad \beta = [0;\overline{1,3,1,3,1,2,1,3,1,3,1,1}].$$ 

\begin{lemma}
One has $\alpha=\max K_1$, $\beta=\max K_2$, and 
$$\frac{\log|f_1'(\alpha)|}{\log|f_2'(\beta)|}\in\mathbb{R}\setminus\mathbb{Q}.$$
\end{lemma} 
\begin{proof}
    We note that $f_1$ and $f_2$ can be rewritten as 
    $$f_1(x)= \frac{67x + 187}{91x + 254}, \quad f_2(y)=\frac{1633y + 2925}{2064y + 3697}.$$
    In particular, 
    $$f_1'(x)=\frac{1}{(91x + 254)^2}, \quad \quad f_2'(y) = \frac{1}{(2064y + 3697)^2}$$ 
    and 
    $$\alpha =\frac{\sqrt{103037} - 187}{182}, \quad \quad \beta = \frac{\sqrt{49321} - 86}{172}.$$
    
   A simple computation yields
   \begin{equation*}
       f_1'(\alpha)=\frac{4}{(\sqrt{103037} + 321)^2} \quad \quad f_2'(\beta)=\frac{1}{(12\sqrt{49321} + 2665)^2}.
   \end{equation*}
   Since $103037=11 \cdot 17 \cdot 19 \cdot 29$ and $49321=31 \cdot 37 \cdot 43$, these numbers generate different quadratic extensions of $\Q$.
\end{proof}

\begin{lemma} 
Consider the intervals $B_{U_1}=[\alpha_1,\alpha]$, $U_1=(\alpha_2,\alpha_1)$, $B_{\widetilde{U}}=[\beta_1,\beta]$ and $\widetilde{U}=(\beta_2,\beta_1)$, where 
\begin{itemize}
\item $\alpha_1$ is the smallest element of $K_1$ of the form $[0;1,2,1,3,\underline{\theta}]$, 
\item $\alpha_2$ is the largest element of $K_1$ of the form $[0;1,2,1,2,\underline{\theta}]$, 
\item $\beta_1$ is the smallest element of $K_2$ of the form $[0;1,3,1,3,\widetilde{\underline{\theta}}]$,
\item $\beta_2$ is the largest element of $K_2$ of the form $[0;1,3,1,2,\widetilde{\underline{\theta}}]$.
\end{itemize}
Then, 
$$\frac{|B_{U_1}|}{|U_1|}<0.97 \quad \textrm{and} \quad \frac{|B_{\widetilde{U}}|}{|\widetilde{U}|}<0.89.$$

Moreover we have that $G^\ell(\alpha_2)=\alpha$ for $\ell=2$ where $G$ is the Gauss map.
\end{lemma}

\begin{proof}
A straightforward computation shows that
\begin{align*}
    \alpha&=[0;\overline{1,2,1,3,1,3,1,2}], \\
    \alpha_1 &= [0; 1, 2, 1, 3, 4, 4, 4, 3, \overline{1, 1, 3, 1, 3, 1, 2, 1, 3, 1, 3, 1}], \\
    \alpha_2 &= [0; \overline{1, 2, 1, 2, 1, 3, 1, 3}], \\
    \beta &= [0;\overline{1,3,1,3,1,2,1,3,1,3,1,1}], \\
    \beta_1 &= [0; 1, 3, 1, 3, 4, 4, 4, 3, \overline{1, 1, 3, 1, 3, 1, 2, 1, 3, 1, 3, 1}], \\
    \beta_2 &= [0; \overline{1, 3, 1, 2, 1, 3, 1, 3, 1, 1, 1, 3}].
\end{align*}
\end{proof}

The previous lemma provides us with the ingredients necessary to apply \Cref{lem:specialized_lemma} by choosing $l = 1$, $q_{e_0} = q_{e_1} = \alpha$, $p_1 = \alpha_2$, $\tilde{q} = \beta$. As consequence, we obtain that $K+\widetilde{K}=g_1(K_1)+g_2(K_2)$ contains infinitely many gaps $(a,b)$ of the form $a=x+y$ where
\begin{multline*}
    x=[3;1,3,1,2,3,1,2,2,1,3,1,3,1,2,(1,2,1,3,1,3,1,2)^{n_1}, \\
    1,2,(1,2,1,3,1,3,1,2)^{n_2},1,2,\dots],
\end{multline*}
and 
\begin{equation*}
    y = [0;1,2,1,3,4,4,3,1,\overline{1,3,1,3,1,2,1,3,1,3,1,1}].
\end{equation*}
Since we have the inequalities $\lambda_0^{+}(21213^*131212)<4.52124$, $\lambda_0^{+}(31213^*1311)<4.5212$, $\lambda_0^{+}(221313^*1212)<4.521$, $\lambda_0^{+}(2213^*13)<4.521$, $\lambda_0^{+}(313^*123)<4.52$, $\lambda_0^{+}(313^*11)<4.5$, from the fact that $x+y>4.52124$ we obtain that $x+y\in M$.

Now we will show that for any choice of positive integers $(n_1,n_2,\dots)$ we have $x+y\in L$. For this end, we will show that there is a sequence of periodic bi-infinite sequences $\overline{c^{(k)}}\in\{1,2,3,4\}^\Z$ with period $c^{(k)}$ such that their Markov (and hence Lagrange values) converge to $x+y$. Since $L$ is closed this will prove that $x+y\in L$. It suffices to consider
\begin{equation*}
    c^{(k)} = (113131213131)^{k}134431213^*13123122131312(12131312)^{n_1}12\dotsb(12131312)^{n_k}12.
\end{equation*}

By the above inequalities we obtain that the Markov value of $\overline{c^{(k)}}$ is attained at the marked position with $*$.

\section{A gap in $M$ near 4.08} \label{sec:gap_near_408}

Through this section and the next one we will work exclusively with bi-infinite sequences $\underline{b}\in\{1,2,3\}^\Z$. Given a finite string $u=(a_i)_{i=-m}^{n}$, through this section we denote 
\begin{equation}\label{eq:lambda_i_minus_123}
    \lambda^{-}_i(u):=\min \{[a_i;a_{i+1},...,a_n,\theta_1]+[a_{i-1};a_{i-2},...,a_{-m},\theta_2]: \theta_1,\theta_2\in \{1,2,3\}^{\mathbb{N}}\},
\end{equation}
and
\begin{equation}\label{eq:lambda_i_plus_123}
    \lambda^{+}_i(u):=\max \{[a_i;a_{i+1},...,a_n,\theta_1]+[a_{i-1};a_{i-2},...,a_{-m},\theta_2];\theta_1,\theta_2\in \{1,2,3\}^{\mathbb{N}}\}.
\end{equation}

Let 
$$\nu' = m\left(\overline{3133}11213^*213212\overline{1333}\right) = 4.08049296860135092\dots$$ 
and 
$$\mu' = m\left(\overline{132311}3311213^*21321\overline{3313}\right) = 4.08049996235899758\dots$$

In the sequel, we will establish the following result: 
\begin{proposition}\label{p.gap-near-408} One has $M\cap(\nu',\mu')=\emptyset$. 
\end{proposition} 

\subsection{Preliminary step: Markov values between $4.0804927$ and $4.08050$}

Recall that $\underline{b}=(a_n)_{n\in\mathbb{Z}}\in(\mathbb{N}^*)^{\mathbb{Z}}$ with $\sqrt{12}<m(\underline{b})=\lambda_0(\underline{b})<\sqrt{20}$ has the form $\underline{b}\in\{1,2,3\}^{\mathbb{Z}}$ and $\underline{b}=\dots 3^* \dots$. In particular, modulo transpositions, its possible extensions are 
$$13^*1, \quad 13^*2, \quad 13^*3, \quad 23^*2, \quad 23^*3, \quad 33^*3$$ 

\begin{remark} In what follows, a string leading to a Markov value $>4.08050$ will be called forbidden and it will be marked in red. 
\end{remark}

\begin{lemma} $\lambda_0^{-}(\textcolor{red}{13^*1})>4.11$, $\lambda_0^{+}(33^*3)<\lambda_0^{+}(23^*3)<\lambda_0^{+}(23^*2)<3.89$. 
\end{lemma}

\begin{corollary} If $3.89<m(\underline{b})=\lambda_0(\underline{b})<4.11$, then $b_{-1}b_0^*b_1 = 13^*2$ or $13^*3$. 
\end{corollary}

The left extensions of the strings above are 
\[
\begin{array}{ccc}
113^*2 & 213^*2 & 313^*2\\
113^*3 & 213^*3 & 313^*3
\end{array}
\]

\begin{lemma} $\lambda_0^{-}(\textcolor{red}{313^*2})>4.12$, $\lambda_0^{+}(113^*3)<\lambda_0^{+}(213^*3)<4.0427$. 
\end{lemma} 

\begin{corollary} If $4.0427<m(\underline{b})=\lambda_0(\underline{b})<4.11$, then $b_{-2}b_{-1}b_0^*b_1 = 113^*2, 213^*2$ or $313^*3$.
\end{corollary}

\begin{remark}\label{r.1} $131$, $3132$ and $2313$ forbidden imply that $313$ extends as $33133$. In the sequel, this is denoted $\textcolor{blue}{313\to 33133}$. 
\end{remark}  

By combining the previous corollary with this remark, we obtain: 

\begin{corollary} If $4.0427<m(\underline{b})=\lambda_0(\underline{b})<4.11$, then $b_{-2}b_{-1}b_0^*b_1 = 113^*2, 213^*2$ or $b_{-3}b_{-2}b_{-1}b_0^*b_1=3313^*3$.
\end{corollary}

The right extensions of the strings above are 
\[
\begin{array}{ccc}
113^*21 & 113^*22 & 113^*23\\
213^*21 & 213^*22 & 213^*23\\
3313^*31 & 3313^*32 & 3313^*33
\end{array}
\]

\begin{lemma} $\lambda_0^{+}(113^*21)<\lambda_0^{+}(113^*22)<4.0658$, $\lambda_0^{-}(\textcolor{red}{213^*23})>\lambda_0^{-}(\textcolor{red}{213^*22}) > 4.1$, $\lambda_0^{+}(3313^*3)<4.075$.  
\end{lemma} 

\begin{remark}\label{r.2} $21322, 21323$ forbidden imply $\textcolor{blue}{2132\to 21321}$. Moreover, $313\to 33133$ and $2132\to 21321$ imply $\textcolor{blue}{1322\to 11322}$. 
\end{remark}

\begin{corollary} If $4.075<m(\underline{b})=\lambda_0(\underline{b})<4.1$, then $b_{-2}b_{-1}b_0^*b_1b_2 = 113^*23$ or $213^*21$. 
\end{corollary} 

The left extensions of the strings above are
\[
\begin{array}{ccc}
1113^*23 & 2113^*23 & 3113^*23\\
1213^*21 & 2213^*21 & 3213^*21
\end{array}
\]

\begin{lemma} $\lambda_0^{+}(3113^*23) < \lambda_0^{+}(2113^*23) < 4.0322$.
\end{lemma} 

\begin{corollary} If $4.075<m(\underline{b})=\lambda_0(\underline{b})<4.1$, then $b_{-3}b_{-2}b_{-1}b_0^*b_1b_2 = 1113^*23$ or  
$$1213^*21, \quad 2213^*21, \quad 3213^*21$$ 
\end{corollary}

The right extensions of these strings are 
\[
\begin{array}{ccc}
1113^*231 & 1113^*232 & 1113^*233\\
1213^*211 & 1213^*212 & 1213^*213\\
2213^*211 & 2213^*212 & 2213^*213\\
3213^*211 & 3213^*212 & 3213^*213
\end{array}
\]

\begin{lemma} $\lambda_0^{+}(1113^*233)<\lambda_0^{+}(1113^*232)<4.0784$, $\lambda_0^{-}(\textcolor{red}{1213^*211})>\lambda_0^{-}(\textcolor{red}{1213^*212})>4.0844$, $\lambda_0^{-}(\textcolor{red}{2213^*211})>4.0807$, $\lambda_0^{+}(2213^*213)<4.0711$, $\lambda_0^{+}(3213^*213)<\lambda_0^{+}(3213^*212)<4.0689$. 
\end{lemma} 

\begin{remark}\label{r.3} $1213211, 1213212$ forbidden imply $\textcolor{blue}{121321\to 1213213}$. 
\end{remark}

\begin{corollary} If $4.0784<m(\underline{b})=\lambda_0(\underline{b})<4.0807$, then $b_{-3}b_{-2}b_{-1}b_0^*b_1b_2b_3=1113^*231$ or  
$$1213^*213, \quad 2213^*212, \quad 3213^*211$$
\end{corollary}

The left extensions of the strings above are 
\[
\begin{array}{ccc}
11113^*231 & 21113^*231 & 31113^*231\\
11213^*213 & 21213^*213 & 31213^*213\\
12213^*212 & 22213^*212 & 32213^*212\\
13213^*211 & 23213^*211 & 33213^*211
\end{array}
\]

\begin{lemma}
\begin{itemize}
    \item $\lambda_0^{+}(11113^*231)<\lambda_0^{+}(21113^*231)<4.0763$, 
    \item $\lambda_0^{-}(\textcolor{red}{31213^*213})>\lambda_0^{-}(\textcolor{red}{21213^*213})>4.0875$. 
    \item $\lambda_0^{+}(12213^*212)<\lambda_0^{+}(22213^*212)<4.0792$, 
\end{itemize}
\end{lemma} 

\begin{remark}\label{r.4} From Remark \ref{r.3} and $31213213, 21213213$ forbidden, we infer that $\textcolor{blue}{121321\to 11213213}$
\end{remark}

\begin{corollary} If $4.0792<m(\underline{b})=\lambda_0(\underline{b})<4.0807$, then $b_{-4}\dots b_0^*\dots b_3=31113^*231$ or 
$$11213^*213, \quad 32213^*212, \quad 13213^*211, \quad 23213^*211, \quad 33213^*211$$
\end{corollary}

Since $131, 2313, 32312$ are forbidden and $2132\to 21321$, the right extensions of the strings above are 
\[
\begin{array}{ccc}
31113^*2311 & 3111\textcolor{red}{3^*2312} & 31113^*\textcolor{red}{2313}\\
11213^*2\textcolor{red}{131} & 11213^*2132 & 11213^*2133\\
32213^*2121 & 32213^*2122 & 32213^*2123\\
13213^*2111 & 13213^*2112 & 13213^*2113\\
23213^*2111 & 23213^*2112 & 23213^*2113\\
33213^*2111 & 33213^*2112 & 33213^*2113
\end{array}
\]

Using that $131$ is forbidden, we have: 

\begin{lemma}
\begin{itemize}
    \item $m(\dots31113^*2311\dots)\leq\lambda_0^{+}(231113^*2311) < 4.0793$,
    \item $m(\dots\textcolor{red}{11213^*2133}\dots)\geq\lambda_0^{-}(\textcolor{red}{2311213^*2133})>4.0806$. 
    \item $\lambda_0^{+}(32213^*2121)<\lambda_0^{+}(32213^*2122)<4.0796$, $\lambda_0^{+}(13213^*2111)<4.0785$,
    \item $\lambda_0^{+}(23213^*2111)<4.0801$, $\lambda_0^{-}(\textcolor{red}{23213^*2112})>4.0821$,
    \item $\lambda_0^{-}(\textcolor{red}{3213^*2113})>4.0832$.
\end{itemize}
\end{lemma} 

\begin{corollary} If $4.0801<m(\underline{b})=\lambda_0(\underline{b})<4.0806$, then $b_{-4}\dots b_0^*\dots b_5=11213^*21321$ or $b_{-4}\dots b_0^*\dots b_4$ is equal to
$$32213^*2123, \quad 13213^*2112, \quad 33213^*2111$$ 
\end{corollary} 

Since $3132$ is forbidden, the left extensions of the strings above are 
\[
\begin{array}{ccc}
111213^*21321 & 211213^*21321 & 311213^*21321\\
132213^*2123 & 232213^*2123 & 332213^*2123\\
113213^*2112 & 213213^*2112 & \textcolor{red}{3132}13^*2112\\
133213^*2111 & 233213^*2111 & 333213^*2111
\end{array}
\]

\begin{lemma} $\lambda_0^{-}(\textcolor{red}{111213^*21321})>\lambda_0^{-}(\textcolor{red}{211213^*21321})>4.081$, $\lambda_0^{+}(332213^*2123) < \lambda_0^{+}(232213^*2123) < 4.0803$, $\lambda_0^{-}(\textcolor{red}{113213^*2112}) > 4.0806$. 
\end{lemma} 

\begin{corollary} If $4.0803<m(\underline{b})=\lambda_0(\underline{b})<4.0806$, then $b_{-5}\dots b_0^*\dots b_5=311213^*21321$ or $b_{-5}\dots b_0^*\dots b_4$ is equal to
$$132213^*2123, \quad 213213^*2112, \quad 133213^*2111, \quad 233213^*2111, \quad 333213^*2111$$
\end{corollary} 

Since $1322\to 11322$ (cf. Remark \ref{r.2}), the right extensions of the strings above are 
\[
\begin{array}{ccc}
311213^*213211 & 311213^*213212 & 311213^*213213\\
1132213^*21231 & 1132213^*21232 & 1132213^*21233\\
213213^*21121 & 213213^*21122 & 213213^*21123\\
133213^*21111 & 133213^*21112 & 133213^*21113\\
233213^*21111 & 233213^*21112 & 233213^*21113\\
333213^*21111 & 333213^*21112 & 333213^*21113
\end{array}
\]

Using that $131$ and $2313$ are forbidden, we have:

\begin{lemma}\hfill
\begin{itemize}
    \item $\lambda_0^{+}(311213^*213211) < 4.08041$. 
    \item \begin{align*}
        \lambda_0^{+}(1132213^*21233)&<\lambda_0^{+}(1132213^*21232)<4.0802, \\
        m(\dots 1132213^*21231\dots)&\leq\lambda_0^{+}(1132213^*212312)<4.08048,
    \end{align*}
    \item $\lambda_0^{-}(\textcolor{red}{213213^*21121})>4.0813$, $\lambda_0^{+}(133213^*21113)<\lambda_0^{+}(133213^*21112)<4.0780$, 
    \item \begin{align*}
        \lambda_0^{+}(233213^*21113)&<\lambda_0^{+}(233213^*21112)<4.078, \\
        m(\dots233213^*21111\dots)&\leq \lambda_0^{+}(233213^*2111132)<4.0804845,
    \end{align*} 
    \item \begin{align*}
        \lambda_0^{+}(333213^*21113)&<\lambda_0^{+}(333213^*21112)<4.078, \\
        m(\dots 333213^*21111\dots)&\leq \lambda_0^{+}(333213^*2111132)<4.08039.
    \end{align*}
\end{itemize}
\end{lemma} 

\begin{corollary} If $4.0804845<m(\underline{b})=\lambda_0(\underline{b})<4.0806$, then $$b_{-5}\dots b_0^*\dots b_6=311213^*213212, 311213^*213213$$ or $b_{-5}\dots b_0^*\dots b_5$ is equal to
$$213213^*21122, \quad 213213^*21123, \quad 133213^*21111$$
\end{corollary} 

Since $131$ is forbidden and $313\to 33133$, the left extensions of the strings above are 
\[
\begin{array}{ccc}
\textcolor{red}{131}1213^*213212 & 2311213^*213212 & 3311213^*213212\\
\textcolor{red}{131}1213^*213213 & 2311213^*213213 & 3311213^*213213 \\
1213213^*21122 & 2213213^*21122 & 3213213^*21122\\
1213213^*21123 & 2213213^*21123 & 3213213^*21123\\
1133213^*21111 & 2133213^*21111 & 33133213^*21111
\end{array}
\]

\begin{lemma}
\begin{itemize}
    \item $\lambda_0^{+}(2311213^*213212)<\lambda_0^{+}(2311213^*213213)<4.08042$,
    \item $\lambda_0^{-}(\textcolor{red}{3213213^*21122})>\lambda_0^{-}(\textcolor{red}{2213213^*21122})>4.0806$, 
    \item $\lambda_0^{+}(1213213^*21123)<\lambda_0^{+}(2213213^*21123)<4.08044$.
\end{itemize}
\end{lemma} 

\begin{corollary} If $4.0804845<m(\underline{b})=\lambda_0(\underline{b})<4.0806$, then $$b_{-6}\dots b_0^*\dots b_6=3311213^*213212, 3311213^*213213$$ or $b_{-6}\dots b_0^*\dots b_5$ is equal to
$$1213213^*21122, \quad 3213213^*21123, \quad 1133213^*21111,$$ 
$$2133213^*21111, \quad 33133213^*21111$$  
\end{corollary} 

Using that $121321\to 11213213$ (cf. Remark \ref{r.4}), the previous corollary can be slightly improved into 
\begin{corollary} If $4.0804845<m(\underline{b})=\lambda_0(\underline{b})<4.0806$, then $$b_{-6}\dots b_0^*\dots b_6=3311213^*213212, 3311213^*213213$$ or $b_{-7}\dots b_0^*\dots b_5$ is equal to
$$11213213^*21122, \quad 3213213^*21123, \quad 1133213^*21111,$$ 
$$2133213^*21111, \quad 33133213^*21111$$  
\end{corollary} 

Since $131$ is forbidden and $2132\to 21321$, the right extensions of the strings above are 
\[
\begin{array}{ccc}
3311213^*2132121 & 3311213^*2132122 & 3311213^*2132123\\
3311213^*2132\textcolor{red}{131} & 3311213^*21321321 & 3311213^*2132133 \\
11213213^*211221 & 11213213^*211222 & 11213213^*211223\\
3213213^*211231 & 3213213^*211232 & 3213213^*211233\\
1133213^*211111 & 1133213^*211112 & 1133213^*211113\\
2133213^*211111 & 2133213^*211112 & 2133213^*211113\\
33133213^*211111 & 33133213^*211112 & 33133213^*211113
\end{array}
\]

Using that $131$ is forbidden, we have: 
\begin{lemma}
\begin{itemize}
    \item $\lambda_0^{+}(3311213^*2132123)<\lambda_0^{+}(3311213^*2132122)<4.08048$.  
    \item $\lambda_0^{-}(\textcolor{red}{11213213^*21122})>4.08052$,
    \item $\lambda_0^{+}(3213213^*211231)<\lambda_0^{+}(3213213^*211232)< 4.08046$, 
    \item $\lambda_0^{+}(1133213^*211111)<\lambda_0^{+}(1133213^*211112)<4.08035$, 
    \item $\lambda_0^{+}(2133213^*211111)<\lambda_0^{+}(2133213^*211112)<4.08041$, 
    \item $\lambda_0^{+}(3133213^*211111)<\lambda_0^{+}(3133213^*211112)<4.08044$, 
    \item \begin{align*}
        \lambda_0^{-}(\textcolor{red}{33133213^*211113})&>\lambda_0^{-}(\textcolor{red}{2133213^*211113})>4.08058 \\
        m(\dots\textcolor{red}{1133213^*211113}\dots)&\geq \lambda_0^{-}(\textcolor{red}{231133213^*211113})>4.0805,
    \end{align*}
\end{itemize}

\end{lemma} 

\begin{corollary} If $4.0804845<m(\underline{b})=\lambda_0(\underline{b})<4.08050$, then $$b_{-6}\dots b_0^*\dots b_7=3311213^*2132121, 3311213^*2132133$$ or $b_{-6}\dots b_0^*\dots b_6=3213213^*211233$ or $b_{-6}\dots b_0^*\dots b_8=3311213^*21321321$.
\end{corollary}

The left extensions of the strings above are 
\[
\begin{array}{ccc}
13311213^*2132121 & 23311213^*2132121 & 33311213^*2132121\\
13311213^*2132133 & 23311213^*2132133 & 33311213^*2132133\\
13213213^*211233 & 23213213^*211233 & 33213213^*211233\\
13311213^*21321321 & 23311213^*21321321 & 33311213^*21321321
\end{array}
\]

\begin{lemma}
\begin{itemize}
    \item $\lambda_0^{+}(33213213^*211233)<\lambda_0^{+}(23213213^*211233)<4.08048$, 
    \item $\lambda_0^{-}(\textcolor{red}{13311213^*21321321})>4.08050$, 
    \item $\lambda_0^{+}(33311213^*21321321)<\lambda_0^{+}(23311213^*21321321)<4.0804924$, 
    \item $\lambda_0^{+}(33311213^*2132121)<\lambda_0^{+}(23311213^*2132121)<4.08046$, 
    \item $\lambda_0^{+}(33311213^*2132133)<\lambda_0^{+}(23311213^*2132133)<4.0804884$.  
\end{itemize}
\end{lemma} 

\begin{corollary} If $4.0804924<m(\underline{b})=\lambda_0(\underline{b})<4.08050$, then $$b_{-7}\dots b_0^*\dots b_7=13311213^*2132121, 13311213^*2132133$$ or $b_{-7}\dots b_0^*\dots b_6=13213213^*211233$.
\end{corollary} 

The right extensions of the strings above are 
\[
\begin{array}{ccc}
13311213^*21321211 & 13311213^*21321212 & 13311213^*21321213\\
13311213^*21321331 & 13311213^*21321332 & 13311213^*21321333\\
13213213^*2112331 & 13213213^*2112332 & 13213213^*2112333
\end{array}
\]

Using that $2313$ is forbidden and $313\to 33133$, we have: 

\begin{lemma} 
\begin{itemize}
    \item \begin{align*}
        \lambda_0^{+}(13213213^*2112333)&<\lambda_0^{+}(13213213^*2112332)<4.08049 \\
        m(\dots13213213^*2112331\dots)&\leq \lambda_0^{+}(213213213^*211233133)<4.0804927,
    \end{align*}
    \item \begin{align*}
        \lambda_0^{+}(13311213^*21321211)&<4.08049 \\
        m(\dots13311213^*21321212\dots)&\leq \lambda_0^{+}(3313311213^*21321212)<4.08049093,
    \end{align*}
    \item $\lambda_0^{-}(\textcolor{red}{13311213^*21321333})>\lambda_0^{-}(\textcolor{red}{13311213^*21321332})>4.08050$.
\end{itemize}
\end{lemma}

\begin{corollary} If $4.0804927<m(\underline{b})=\lambda_0(\underline{b})<4.08050$, then $b_{-7}\dots b_0^*\dots b_8=13311213^*21321213$ or $13311213^*21321331$. 
\end{corollary}

Since $313\to 33133$, the left extensions of the strings above are 
\[
\begin{array}{ccc}
113311213^*21321213 & 213311213^*21321213 & 3313311213^*21321213\\
113311213^*21321331 & 213311213^*21321331 & 3313311213^*21321331
\end{array}
\]

\begin{lemma} $\lambda_0^{+}(113311213^*21321213)<\lambda_0^{+}(213311213^*21321213)<4.0804914$, $\lambda_0^{-}(\textcolor{red}{3313311213^*21321331})>\lambda_0^{-}(\textcolor{red}{213311213^*21321331})>4.08051$
\end{lemma} 

\begin{corollary} If $4.0804927<m(\underline{b})=\lambda_0(\underline{b})<4.08050$, then $$b_{-9}\dots b_0^*\dots b_8=3313311213^*21321213$$ or $b_{-8}\dots b_0^*\dots b_8=113311213^*21321331$. 
\end{corollary}

Since $131$ is forbidden, $2132\to 21321$, $121321\to 11213213$ (so that \\ $\textcolor{blue}{12132\to 11213213}$) and $313\to 33133$, the right extensions of the strings above are 
\[
\begin{array}{ccc}
3313311213^*213212\textcolor{red}{131} & 3313311213^*2132\textcolor{red}{12132} & 3313311213^*213212133\\
113311213^*213213311 & 113311213^*213213312 & 113311213^*2132133133
\end{array}
\]

Using that $131$ is forbidden, we have 

\begin{lemma}
\begin{align*}
    m(\dots\textcolor{red}{113311213^*213213311}\dots)&>m(\dots\textcolor{red}{113311213^*213213312}\dots) \\
    &\geq \lambda_0^{-}(\textcolor{red}{23113311213^*213213312})> 4.0805.
\end{align*}
\end{lemma} 

\begin{corollary}\label{cor:gap_408} If $4.0804927<m(\underline{b})=\lambda_0(\underline{b})<4.08050$, then $$b_{-9}\dots b_0^*\dots b_9=3313311213^*213212133$$ or $b_{-8}\dots b_0^*\dots b_{10}=113311213^*2132133133$. 
\end{corollary}

\subsection{Proof of Proposition \ref{p.gap-near-408}} 

The desired proposition is an immediate consequence of the following next two lemmas. 

\begin{lemma}\label{l.alpha} If $b_{-9}\dots b_0^*\dots b_9=3313311213^*213212133$ has Markov value $m(\underline{b})<4.08050$, then $m(\underline{b})\leq \nu'$. 
\end{lemma}

\begin{lemma}\label{l.beta} If $b_{-8}\dots b_0^*\dots b_{10}=113311213^*2132133133$, then $m(\underline{b})\geq\mu'$. 
\end{lemma}

The remainder of this section discusses the proofs of these two statements. 

\subsubsection{Proof of Lemma \ref{l.alpha}} Let $b_{-9}\dots b_0^*\dots b_9=3313311213^*213212133$ with Markov value $m(\underline{b})<4.08050$. Given the parities of the extremal symbols in this string, the largest values of $m(\underline{b})$ are obtained by extending it as 
$$b_{-12}\dots b_0^*\dots b_{12}=\textcolor{green}{313}3313311213^*213212133\textcolor{green}{313}$$ 
Since $313\to 33133$, $\underline{b}$ extends as 
$$b_{-13}\dots b_0^*\dots b_{13}=\textcolor{green}{3313}3313311213^*213212133\textcolor{green}{3133}$$ 
At this point, it is clear that we can repeat the argument (with parities and $313\to 33133$) to conclude that the largest extension of $\underline{b}$ is 
$$\nu'=\overline{3313}3313311213^*213212133\overline{3133}$$

\subsubsection{Proof of Lemma \ref{l.beta}} Let $b_{-8}\dots b_0^*\dots b_{10}=113311213^*2132133133$. Using the same argument from the previous subsection to the right of this string, we see that the smallest Markov value of $m(\underline{b})$ is obtained by extending it as 
$$b_{-8}\dots b_0^*\dots=113311213^*2132133133\overline{3133}$$ 
Next, given the parity of the left extreme of this string and the fact that $131$ is forbidden, we have that the smallest values of $m(\underline{b})$ are obtained by extending it as 
$$b_{-10}\dots b_0^*\dots=23113311213^*2132133133\overline{3133}$$ 
Using again the parity argument, the smallest values of $m(\underline{b})$ correspond to 
$$b_{-12}\dots b_0^*\dots=\textcolor{green}{13}23113311213^*2132133133\overline{3133}$$ 
Now, we recall that $2313$ and $21323$ are forbidden to get that the unique allowed left extension of this string is 
$$b_{-13}\dots b_0^*\dots=\textcolor{green}{113}23113311213^*2132133133\overline{3133}$$ 
By applying twice the parity argument, we get that the smallest Markov values of $m(\underline{b})$ are reached by 
$$b_{-15}\dots b_0^*\dots=\textcolor{green}{31113}23113311213^*2132133133\overline{3133}$$ 
Since $131$ is forbidden, the parity argument yields that the smallest Markov values of $m(\underline{b})$ correspond to 
$$b_{-16}\dots b_0^*\dots=\textcolor{green}{231113}23113311213^*2132133133\overline{3133}$$ 
At this point, we can repeat the argument above to conclude that the smallest extension of $\underline{b}$ is 
$$\mu' = \overline{231113}23113311213^*2132133133\overline{3133}$$

\subsection{Uncountable boundary}

The following lemma is a consequence of \Cref{cor:gap_408}, \Cref{l.alpha} and \Cref{l.beta}.

\begin{lemma}
If $4.0804927<m(\underline{a})=\lambda_0(\underline{a})<4.0805$, then up to transposition, either 
$$\underline{a} = \dots 3313311213^*213212133\dots$$
and $m(\underline{a})\leq \nu'$ or 
$$\underline{a} = \dots 113311213^*2132133133\dots$$ and $m(\underline{a})\geq\mu'$.
\end{lemma}

Moreover, the proof of \Cref{l.alpha} implies the following lemma.

\begin{lemma}
    The set $M\cap[\mu,4.0805)$ is contained in $K+\widetilde{K}$ where
$$K=\{[3;2,1,3,2,1+x]: x\in K_1\},$$ 
$$\widetilde{K}=\{[0;1,2,1,1,3,3,1+y]: y\in K_2\},$$ 
and $K_1=\{[0;3,3,1,3,\underline{\theta}]\in \widehat{K}\}$, $K_2=\{[0;1, 3, 2, 3, 1, 1,\widetilde{\underline{\theta}}]\in\widehat{K}\}$ with 
\begin{multline*}
    \widehat{K}=\{[0;\widehat{\underline{\theta}}]:\widehat{\underline{\theta}}\in\{1,2,3\}^{\mathbb{N}} \textrm{ doesn't contain the strings } 131, 3132, 32312, 312132,  \\ 1112132, \textrm{ or their transposes}\}.
\end{multline*}
\end{lemma}

Observe that $\widehat{K}$ is a Gauss-Cantor set, so it is a dynamically defined Cantor set given by iterates of the Gauss map (as in \Cref{example}) and $K_1$ and $K_2$ are conjugated to $\widehat{K}$ (up to refine the Markov partition of $\widehat{K}$, see \Cref{rem:smaller_copies}). We will show that the hypothesis of \Cref{lem:specialized_lemma} hold for $(K_1,K_2)$. Since $K=g_1(K_1)$ and $\widetilde{K}=g_2(K_2)$ for the decreasing $C^1$ contractions
\begin{equation*}
    g_1(x)=[3;2,1,3,2,1+x] \quad\text{and}\quad g_2(y) = [0;1,2,1,1,3,3,1+y], 
\end{equation*}
by \Cref{cor:conjugating_preserves_lemma} we will obtain that $\partial(K+\widetilde{K})$ is uncountable arbitrarily close to $\min K+\min\widetilde{K} = \mu'$.

Consider the following inverse branches $f_1$, $f_2$ of $K_1$ and $K_2$, respectively,
$$f_1(x)=[0;3,3,1,3+x] \quad \text{ and } \quad f_2(y)=[0;1,3,2,3,1,1+y]$$ 
whose fixed points are 
$$\alpha=[0;\overline{3,3,1,3}] \quad \text{ and } \quad \beta = [0;\overline{1,3,2,3,1,1}].$$

\begin{lemma}
One has $\alpha=\max K_1$, $\beta=\max K_2$, and 
$$\frac{\log|f_1'(\alpha)|}{\log|f_2'(\beta)|}\in\mathbb{R}\setminus\mathbb{Q}.$$
\end{lemma} 
\begin{proof}
    We note that $f_1$ and $f_2$ can be rewritten as 
    $$f_1(x)= \frac{4x + 15}{13x + 49}, \quad f_2(y)=\frac{31y + 55}{40y + 71}.$$
    In particular, 
    $$f_1'(x)=\frac{1}{(13x + 49)^2}, \quad \quad f_2'(y) = \frac{1}{(40y + 71)^2}$$ 
    and 
    $$\alpha =\frac{\sqrt{2805}-45}{26}, \quad \quad \beta = \frac{\sqrt{26}-2}{4}.$$
   A simple computation yields
   \begin{equation*}
       f_1'(\alpha)=\frac{4}{(\sqrt{2805} + 53)^2} \quad \quad f_2'(\beta)=\frac{1}{(10\sqrt{26} + 51)^2}.
   \end{equation*}
   Since $2805=3 \cdot 5 \cdot 11 \cdot 17$ and $26=2\cdot 13$, these numbers generate different quadratic extensions of $\Q$.
\end{proof}

\begin{lemma} 
Consider the intervals $B_{U_1}=[\alpha_1,\alpha]$, $U_1=(\alpha_2,\alpha_1)$, $B_{\widetilde{U}}=[\beta_1,\beta]$ and $\widetilde{U}=(\beta_2,\beta_1)$, where 
\begin{itemize}
\item $\alpha_1$ is the smallest element of $K_1$ of the form $[0;3,3,1,3,\underline{\theta}]$, 
\item $\alpha_2$ is the largest element of $K_1$ of the form $[0;3,3,1,2,\underline{\theta}]$, 
\item $\beta_1$ is the smallest element of $K_2$ of the form $[0;1,3,2,3,1,1,1,3,\widetilde{\underline{\theta}}]$,
\item $\beta_2$ is the largest element of $K_2$ of the form $[0;1,3,2,3,1,1,1,2,\widetilde{\underline{\theta}}]$.
\end{itemize}
Then, 
$$\frac{|B_{U_1}|}{|U_1|}<0.072 \quad \textrm{and} \quad \frac{|B_{\widetilde{U}}|}{|\widetilde{U}|}<0.293.$$

Moreover we have that $G^\ell(\alpha_2)=\alpha$ for $\ell=6$ where $G$ is the Gauss map.
\end{lemma}

\begin{proof}
A straightforward computation shows that
\begin{align*}
    \alpha&=[0;\overline{3,3,1,3}], &\beta &= [0;\overline{1,3,2,3,1,1}], \\
    \alpha_1 &= [0; 3, 3, 1, \overline{3, 3, 1, 3}], &\beta_1 &= [0; 1, 3, 2, 3, 1, 1, 1, \overline{3, 3, 1, 3}] \\
    \alpha_2 &= [0; 3, 3, 1, 2, \overline{1, 3, 3, 3}], &\beta_2 &= [0; 1, 3, 2, 3, 1, 1, 1, 2, \overline{1, 3, 3, 3}].
\end{align*}
\end{proof}

The previous lemma provides us with the ingredients necessary to apply \Cref{lem:specialized_lemma} by choosing $l = 1$, $q_{e_0} = q_{e_1} = \alpha$, $p_1 = \alpha_2$, $\tilde{q} = \beta$. As consequence, we obtain that $K+\widetilde{K}=g_1(K_1)+g_2(K_2)$ contains infinitely many gaps $(a,b)$ of the form $b=x+y$ where
\begin{equation*}
    x=[3;2,1,3,2,1,(3,3,1,3)^{n_1},3,3,1,2,1,3,(3,3,1,3)^{n_2},\dots]
\end{equation*}
and
\begin{equation*}
    y=[0;1,2,1,1,3,3,1,\overline{1,3,2,3,1,1}]
\end{equation*}
The inequalities $\lambda_0^{+}(33^*3)<3.7$, $\lambda_0^{+}(33^*133)<4.075$, $\lambda_0^{+}(33^*12)<4.05$, $\lambda_0^{+}(3213^*213)<4.06$, $\lambda_0^{+}(33^*11)<4$, $\lambda_0^{+}(3113^*2)<4.01$ and $\lambda_0^{+}(231113^*2311)<4.08$ show that $x+y\in M$. Now we will show that for any choice of positive integers $(n_1,n_2,\dots)$ we have $x+y\in L$. For this end, we will show that there is a sequence of periodic bi-infinite sequences $\overline{c^{(k)}}\in\{1,2,3\}^\Z$ with period $c^{(k)}$ such that their Markov (and hence Lagrange values) converge to $x+y$. Since $L$ is closed this will prove that $x+y\in L$. It suffices to consider
\begin{equation*}
    c^{(k)} = (113231)^k13311213^*21321(3313)^{n_1}331213\dotsb(3313)^{n_k}311.
\end{equation*}
By the above inequalities and $\lambda_0^{+}(11113^*2)<4.07$, we obtain that the Markov value of $\overline{c^{(k)}}$ is attained at the marked position with $*$.

\section{A point of $M_3\setminus L_3$ near 4.524} \label{sec:M_3_minus_L_3}

Several papers have explicitly constructed elements of the difference set $M\setminus L$ by establishing the local uniqueness and self-replication properties of non-semisymmetric orbits; see for instance \cite{Fractalgeometryofcomplement}. In view of these precedents, we omit the lengthy and routine computations involved in such constructions, as the main point for our purposes is the existence of these examples.

Recall that $\underline{b}\in\{1,2,3\}^\Z$ and that $\lambda_0^{-}$ and $\lambda_i^{+}$ where defined on \eqref{eq:lambda_i_minus_123} and \eqref{eq:lambda_i_plus_123}.

Lets call $w:=3131213131213$ and $w^*:=3131213^*131213$. Let us denote 
\begin{equation*}
    j_0:=m\left(\overline{3131213^*131213}\right)=4.5249406323331145159085648\dots
\end{equation*}

\begin{lemma}
    Suppose that $\underline{b}\in\{1,2,3\}^\Z$ contains any of the following words
    \begin{gather*}
        31313, w2, 1w, w1, 2w, w32, 33w, w33, 23w, w311, 113w, 3213w, w3123, 2213w, \\ w3122, 11213w, 21213w, 213w312, 2313w312, 213w3132, 3313w312, 213w3133, \\
        21313w312, 213w31311, 11313w312, 213w313123, 213w313122, \\
        231213w31312, 331213w31312.
    \end{gather*}
    Then $m(\underline{b})>4.524940646>j_0+10^{-8}$.
\end{lemma}

\begin{lemma}
    Suppose that $\underline{b}\in\{1,2,3\}^\Z$ contains any of the following words
    \begin{gather*}
        1131213w31312, 2131213w31312, 213w3131211, 213w3131212, 213w31312132, \\
        213w31312133, 23131213w31312, 33131213w31312, 113131213w31312,  \\
        213w313121311, 213131213w313121312, 213131213w3131213132,  \\
        213131213w3131213133, 2213131213w31312, 3213131213w31312, \\
        11213131213w31312, 21213131213w31312, 213131213w31312131311, \\
        213131213w313121313122, 213131213w313121313123, 231213131213w31312131312, \\
        331213131213w31312131312, 1131213131213w31312131312, \\
        2131213131213w31312131312, 213131213w3131213131211.
    \end{gather*}
    Then $m(\underline{b})>4.52494063233361>j_0+5\cdot10^{-13}$. 
\end{lemma}

\begin{lemma}
    The Markov value
    \begin{align*}
        j_1 &:= m\left(\overline{21213131213131}111312131213131213ww^*\overline{31312131312121}\right) \\
        &=4.5249406323332657955554350\dots,
    \end{align*}    
    is the minimum Markov value of a bi-infinite sequence $\underline{b}\in\{1,2,3\}^\Z$ containing $213131213w3131213131212$. Moreover $j_1\in L^\prime$ and $j_1-j_0\approx1.51\cdot 10^{-13}$.
\end{lemma}

\begin{lemma}[Local uniqueness]
Suppose $\underline{b}\in\{1,2,3\}^\Z$ satisfies $j_0-10^{-9}<\lambda_0(\underline{b})=m(\underline{b})<j_0+10^{-9}$. Then $b_{-9}\dots b_0^*b_1\dots b_6=213w^*$.
\end{lemma}

\begin{lemma}[Self-replication]\label{lem:self-replication}
Suppose $\underline{b}\in\{1,2,3\}^\Z$ satisfies $b_{-9}\dots b_0^*b_1\dots b_6=213w^*$, $m(\underline{b})<j_1$. Then $b_{-12}\dots b_0^*b_1\dots b_{25}=131213w^*w313121$.
\end{lemma}

\begin{corollary}\label{cor:first_characterization}
Suppose $\underline{b}\in\{1,2,3\}^\Z$ satisfies $j_0-10^{-9}<\lambda_0(\underline{b})=m(\underline{b})<j_1$. Then $b_{-12}\dots b_0^*b_1\dots=131213w^*\overline{w}$.
\end{corollary}

\begin{corollary}
Suppose $\underline{b}\in\{1,2,3\}^\Z$ is periodic and satisfies $j_0-10^{-9}<m(\underline{b})<j_1$. Then $\underline{b}=\overline{w^*}$ and $m(\underline{b})=j_0$.
\end{corollary}

Define $F$ to be the following set of words, including its transposes:
\begin{gather*}
    31313, w2, w1, 1w, 2w, 33w, w32, w33, 23w, w312, w311, 113w,  w_r.
\end{gather*}
where $w_r:=213w$ is the self-replicating word.

\begin{theorem}\label{thm:characterization}
We have that $L\cap(j_0,j_1)=\varnothing$ and \[(M\backslash L)\cap(j_0,j_1)=C_1\cup C_2\cup C_2\cup C_4\cup D,\] where
\begin{multline*}
C_1=\{\lambda_0\left(\gamma313ww^*\overline{w}\right):\gamma\in\{1,2,3\}^{\N}~\text{and} \\ \text{$\gamma313ww$ does not contain any word from $F$}\},
\end{multline*}
\begin{multline*}
C_2=\{\lambda_0\left(\gamma21312131213131213w^*\overline{w}\right):\gamma\in\{1,2,3\}^{\N}~\text{and} \\ \text{$\gamma21312131213131213w$ does not contain any word from $F$}\},
\end{multline*}
\begin{multline*}
C_3=\{\lambda_0\left(\gamma31312131213131213w^*\overline{w}\right):\gamma\in\{1,2,3\}^{\N}~\text{and} \\ \text{$\gamma31312131213131213w$ does not contain any word from $F$}\},
\end{multline*}
\begin{multline*}
C_4=\{\lambda_0\left(\gamma111312131213131213w^*\overline{w}\right):\gamma\in\{1,2,3\}^{\N}~\text{and} \\ \text{$\gamma111312131213131213w$ does not contain any word from $F$}\},
\end{multline*}
are Cantor sets and the discrete points
\begin{multline*}
    D=\{\lambda_0\left(\gamma131213131213ww^*\overline{w}\right):\gamma\in\{1,2,3\}^{\N}, \\
    \text{$\gamma131213131213$ contains $w_r^T=w^T312$}, \\ \text{$\gamma131213131213ww$ does not contain any word from $F\cup\{w_r,w_r^T\}$} \\ 
    \text{and $\gamma131213131213ww^*\overline{w}$ attains its Markov value at $w^*$}\}
\end{multline*}
\end{theorem}

\begin{remark}
    It is possible to give a more detailed and explicit description of $D$ but we refrained from doing so to avoid excessive computations.
\end{remark}

\section{An algorithm for approximating the Markov spectrum} \label{sec:algorithm}

In this section, we describe an algorithm for constructing coverings of the Markov spectrum with arbitrarily small maximum interval length, and provide its pseudocode. Our method is an enhancement of the algorithm developed in \cite{Newgaps}. Its main advantage is that it is \emph{hereditary}: once a covering has been computed at a given precision, it can be used as the starting point for constructing a finer covering, without having to repeat the previous computations. This contrasts with the algorithm in \cite{VincentAlgorithm}, which is not hereditary: constructing the weighted graph $G_{K,Q}$ for different values of $K$ requires starting the computation from scratch. We emphasize, however, that our algorithm produces a covering of the Markov spectrum, but does not guarantee that every interval in the covering intersects the spectrum.

Moreover, since the algorithm subdivides each interval independently of the others, the computation is readily parallelizable. The same feature also allows us to approximate $M\cap[s_1,s_2]$, for any $s_1<s_2$, to arbitrary precision without first computing a covering of the entire set $M\cap(-\infty,s_2]$.

Given $A\in\N^*$, define the Gauss-Cantor set
\begin{equation*}
    C(A):=\{[0;a_1,a_2,\dots]:1\leq a_i\leq A,\forall i\in\N^*\}.
\end{equation*}

Given $A\in\N^*$ and a finite word $w=a_{-p}\dots a_q\in\{1,\dots,A\}^{[-p,q]\cap\Z}$, define the function
\begin{align*}
    \lambda_0^{\max}(w,A)&=\max\{\lambda_0(\underline{b}):\underline{b}\in\{1,\dots,A\}^\Z:b_{-p}\dots b_q = w\}, \\
    \lambda_0^{\min}(w,A)&=\min\{\lambda_0(\underline{b}):\underline{b}\in\{1,\dots,A\}^\Z:b_{-p}\dots b_q = w\}
\end{align*}

Note that these functions coincide with the following expressions: 
\begin{align*}
    \lambda_0^{\max}(w,A)&=a_0+\max(I(a_1\dots a_q)\cap C(A))+\max(I(a_{-1}\dots a_{-p})\cap C(A)) \\
    \lambda_0^{\min}(w,A)&=a_0+\min(I(a_1\dots a_q)\cap C(A))+\min(I(a_{-1}\dots a_{-p})\cap C(A)).
\end{align*}

Given $-p\leq i\leq q$, we define 
\begin{align*}
    \lambda_i^{\max}(w,A)&=\lambda_0^{\max}(a_{-p}\dots a_i^*\dots a_q,A) \\
    \lambda_i^{\min}(w,A)&=\lambda_0^{\min}(a_{-p}\dots a_i^*\dots a_q,A).
\end{align*}

We define the closed interval $J(w,A)=J(a_{-p}\dots a_p,A)$ as the convex hull of the set of values of $\lambda_0(\underline{b})$ for all bi-infinite sequences $\underline{b}\in\{1,\dots,A\}^\Z$ such that $b_{-p}\dots b_q=a_{-p}\dots a_q$. In other terms, the interval $J(w,A)$ coincides with the following interval
\begin{equation*}
    J(w,A):=[\lambda_0^{\min}(w,A),\lambda_0^{\min}(w,A)].
\end{equation*}

The following technical Lemma allows one to compute the
interval $J(a_{-p}\dots a_q)$ explicitly. 
\begin{lemma}
    \label{lem:intlim}
    Let $w=a_{-p}\dots a_q \in \{1,\dots,A\}^{[-p,q]\cap\Z}$ be a finite word. For any bi-infinite sequence $\underline{b}\in\{1,\dots,A\}^{\Z}$ such that $b_{-p}\dots b_q=a_{-p}\dots a_q$, we have the bounds
    \begin{equation*}
    	\lambda_0^{\min}(w,A)\leq\lambda_0(\underline{b})\leq\lambda_0^{\max}(w,A).
    \end{equation*}
    
    Moreover
    $$
    \lambda_0^{\max}(a_{-p}\dots a_q,A) =
    \begin{cases}
    [a_0; a_1, \dots, a_q, \overline{1,A}] + 
    [0;a_{-1},, \dots, a_{-p}, \overline{1,A}], &\mbox{ if } p \mbox{ and }
    q \mbox{ are even,} \\
    [a_0; a_1, \dots, a_q, \overline{A,1}] + 
    [0;a_{-1}, \dots, a_{-p}, \overline{1,A}], &\mbox{ if } p \mbox{ is
    even and } q \mbox{ is odd,} \\ 
    [a_0; a_1, \dots, a_q, \overline{1,A}] + 
    [0;a_{-1}, \dots, a_{-p}, \overline{A,1}], &\mbox{ if } p \mbox{ is
    odd and } q \mbox{ is even,} \\ 
    [a_0; a_1, \dots, a_q, \overline{A,1}] + 
    [0;a_{-1}, \dots, a_{-p}, \overline{A,1}], &\mbox{ if } p \mbox{ and }
    q \mbox{ are odd;} 
    \end{cases}
    $$
    and 
    $$
    \lambda_0^{\min}(a_{-p}\dots a_q,A) =  
    \begin{cases}
    [a_0; a_1, \dots, a_q, \overline{A,1}] + 
    [0;a_{-1}, \dots, a_{-p}, \overline{A,1}], &\mbox{ if } p \mbox{ and }
    q \mbox{ are even,} \\
    [a_0; a_1, \dots, a_q, \overline{1,A}] + 
    [0;a_{-1}, \dots, a_{-p}, \overline{A,1}], &\mbox{ if } p \mbox{ is
    even and } q \mbox{ is odd,} \\ 
    [a_0; a_1, \dots, a_q, \overline{A,1}] + 
    [0;a_{-1}, \dots, a_{-p}, \overline{1,A}], &\mbox{ if } p \mbox{ is
    odd and } q \mbox{ is even,} \\ 
    [a_0; a_1, \dots, a_q, \overline{1,A}] + 
    [0;a_{-1}, \dots, a_{-p}, \overline{1,A}], &\mbox{ if } p \mbox{ and }
    q \mbox{ are odd.} 
    \end{cases}
    $$
\end{lemma}

Let us recall \cite[Lemma 6, Chapter 1]{Cusick-Flahive}.

\begin{lemma}\label{lem:canonical_sequence}
Given $m\in M$, there is a bi-infinite sequence $\underline{a}\in\{1,\dots,\lfloor m\rfloor\}^\Z$ such that $\lambda_0(\underline{a})=m(\underline{a})=m$.
\end{lemma}

The previous lemma is our starting point for approximate the Markov spectra. The basic idea is that we can obtain approximations of $M$ up to a certain value of $t$, just by finding all ``admissible" intervals $J$. More precisely, given $t\in(0,\infty)$ define
\begin{multline}\label{eq:B_c}
	B_c(t,A):=\{a_{-N}\dots a_N\in\{1,\dots,A\}^{[-N,N]\cap\Z}:m(\underline{b})=\lambda_0(\underline{b})\leq t, \\ \text{for some } \underline{b}\in\{1,\dots,A\}^\Z,b_{-N}\dots b_N=a_{-N}\dots a_N \}.
\end{multline}

Observe that if $m(\underline{b})\leq t$ then any symbol larger than $\lfloor t\rfloor$ can not appear in $\underline{b}$. In particular, we can cover the Markov spectrum with the intervals
\begin{equation*}
	M\cap(-\infty,t]\subset\bigcup_{a_{-N}\dots a_N\in B_c(t,\lfloor t\rfloor)} J(a_{-N}\dots a_N).
\end{equation*}

However, the previous covering could be using much more symbols than needed. The optimal number of symbols needed, can be were found in \cite{DynamicalCharacterization}, Given any $A\in\N^*$, $A\neq 3$, we have that $T(A)=\sqrt{A^2+4A}>0$ are the heights at which the equality 
\begin{equation*}
    M\cap(-\infty,\sqrt{A^2+4A}]=\{m(\underline{a}):\underline{a}\in\{1,\dots,A\}^\Z\}
\end{equation*}
holds. The authors also proved that the equality does not hold for $A=3$. For our purpose, we do not need such a precise information. Indeed, the only information we need is the following:
\begin{enumerate}
\item For values of $m\in M$ with $m\leq\sqrt{12}$, we only need to consider bi-infinite sequences in $\{1,2\}^\Z$. In other words, one has the equality 
\begin{equation}\label{eq:alphabet_12}
    M\cap(-\infty,\sqrt{12}]=\{m(\underline{a}):\underline{a}\in\{1,2\}^\Z\}.
\end{equation}
To prove this, we only need to consider that if a bi-infinite sequence $\underline{a}\in\{1,2,3\}^\Z$ contains $3$, the one has the inequality $m(\underline{a})\leq 3+2\cdot[0;\overline{3,1}]=\sqrt{7/3} + 2>\sqrt{12}$. 
\item For values in $m\in M$ with $\sqrt{13}\leq m\leq 5$, we have that
\begin{equation}\label{eq:alphabet_1234}
    M\cap(-\infty,5]=\{m(\underline{a}):\underline{a}\in\{1,2,3,4\}^\Z\}.
\end{equation}
This follows trivially from the fact that if a bi-infinite sequence $\underline{a}\in(\N^*)^\Z$ contains $5$ then $m(\underline{a})>5$. 
\end{enumerate}

Using the main result of \cite{DynamicalCharacterization} and from the proof of \Cref{lem:canonical_sequence}, we have

\begin{lemma}\label{lem:naive_covering}
Given any positive integer $N$ and $A\in\N^*$, $A\neq 3$ we have for any $t\leq\sqrt{A^2+4A}$,
\begin{equation*}
    M\cap(-\infty,t]\subset\bigcup_{a_{-N}\dots a_N\in B_c(t,A)} J(a_{-N}\dots a_N).
\end{equation*}
\end{lemma}

\begin{remark}
As consequence of \eqref{eq:alphabet_12} and \eqref{eq:alphabet_1234}, the case $A=2$ and $A=4$ with $t\leq 5$ do not need \cite{DynamicalCharacterization}.
\end{remark}

\begin{definition}[Allowed word]
Given a word $a_{-p}\dots a_0\dots a_q\in\{1,\dots,A\}^{[-p,q]\cap\Z}$, we will save the word is \emph{allowed} in the alphabet $\{1,\dots,A\}$ if for all $i\in[-p,q]\cap\Z$ one has
\begin{equation*}
	\lambda_0^{\min}(a_{-p}\dots a_i^*\dots a_q,A)\leq \lambda_0^{\max}(a_{-p}\dots a_0^*\dots a_q,A).
\end{equation*}
\end{definition}

In other terms, the word $a_{-p}\dots a_q$ is not allowable if there does not exists a bi-infinite sequence $\underline{b}\in(\N^*)^\Z$ with $b_{-p}\dots b_q=a_{-p}\dots a_q$ that attains its Markov value at the position of $a_0$, that is, such that $m(\underline{b})=\lambda_0(\underline{b})$. 

\begin{definition}[Critical positions]\label{def:critical_positions}
Given an allowed word $a_{-p}\dots a_0\dots a_q\in\{1,\dots,A\}^{[-p,q]\cap\Z}$ in the alphabet $\{1,\dots,A\}$, we define its \emph{critical positions} as the indices $i\in[-p,q]\cap\Z$ such that
\begin{equation*}
	\lambda_0^{\max}(a_{-p}\dots a_i^*\dots a_q,A)\geq\lambda_0^{\min}(a_{-p}\dots a_0^*\dots a_q).
\end{equation*}
We use $\Cp(a_{-p}\dots a_0\dots a_q,A)$ to denote the set of critical positions or just \\ $\Cp(a_{-p}\dots a_0\dots a_q)$ when the alphabet $A$ is clear from the context.
\end{definition}

\subsection{Description of the algorithm to approximate the Markov spectra}

Given a finite word $w$, we will denote $\pal(w)$ a boolean value that indicates whether the word $w$ is a palindrome or not, equivalently $\pal(w)=1$ iff $w$ is palindrome.

\begin{enumerate}
\item The initial inputs of the algorithm are a finite alphabet $\cA=\{1,\dots,A\}$, an approximation error $\rho>0$ and an interval $[l,r]\subset(0,\infty)$. 
\item For our initial setup, we initialize $\mathcal{F}\gets\varnothing$ and set \[\cW_0=\{(1^{*},\{0\},1),(2^{*},\{0\},1),\dots,(A^{*},\{0\},1)\}.\]
At any step of the algorithm, the lists $\cW_n$ and $\mathcal{F}$ will consist of tuples $(w^*,\Cp(w,A),\pal(w))$ where $w=a_{-m}\dots a_0^*\dots a_m$, $\Cp(w,A)$ corresponds to the critical positions of $w$ (see \Cref{def:critical_positions}) and $\pal(w)$ is a boolean value that indicates whether $w$ is palindrome or not. For $\cW_n$ all such words have the same lateral length $m=n$. We set $\Delta=0$ which will be the current approximation error and we set $n\gets 1$ to be the length of the lateral subwords we will iterate in the next step, namely the length of $a_1\dots a_{n}$.
\item We initialize $\cW_n\gets\varnothing$ and we will iterate over the elements of $\cW_{n-1}$. 
\item Given an element $(w,\Cp(w,A),\pal(w))\in \cW_{n-1}$, we will iterate over the extensions
\begin{equation*}
    \mathcal{E}(w) = \begin{cases}
        \cA\times\cA, & \pal(w) = 0, \\
        \{(x,y)\in \cA\times\cA:x\leq y\}, & \pal(w) = 1.
    \end{cases}
\end{equation*}
\item We iterate over $(x,y)\in\mathcal{E}(w)$. We set the trial word $w^\prime = xw^*y$ and we set $D(w^\prime) = \{-n\}\cup \Cp(w,A)\cup \{n\}$. Observe that $\Cp(w^\prime,A)\subset D(w^\prime)$. 
\item We test if the trial word actually could lie on the range $[l,r]$. For this we check if $\lambda_0^{\max}(w^\prime,A)<l$ or $\lambda_0^{\min}(w^\prime,A)>r$ hold, in which case we go back to Step 5.
\item We set $U(w^\prime)=\lambda_0^{\max}(w^\prime,A)$ and $L(w^\prime)=\lambda_0^{\min}(w^\prime,A)$. We will use the first quantity to determine if the $xw^*y$ is allowable and use the second to determine the new dangerous positions. 
\item We initialize the new critical positions with $\{0\}$. We iterate over $i\in D(w^\prime)$, $i\neq 0$ to compute $\Cp(w^\prime,A)$. First, we test whether $\lambda_i^{\min}(w^\prime,A)>U(w^\prime)$. If this condition holds then the word is not allowable, so it is discarded and we go back to Step 5. If it does not hold, then we test if $\lambda_i^{\max}(w^\prime,A)\geq L(w^\prime)$. If this last condition holds, then we keep $i$ in the new list of dangerous positions and then we continue with the next $i$. In other words, if we never hit the break condition $\lambda_i^{\min}(w^\prime,A)>U(w^\prime)$, then we have that
\begin{equation*}
    \Cp(w^\prime,A) = \{i \in D(w^\prime): \lambda_i^{\max}(w^\prime,A)\geq L(w^\prime)\}.
\end{equation*}
\item The previous test ensures that $w^\prime=xw^*y$ is allowable and also computes $\Cp(w^\prime,A)$. If both $\pal(w)=1$ and $x=y$ we have $\pal(xwy)=1$, otherwise $\pal(xwy)=0$. 
\item We set $\delta(w^\prime)=U(w^\prime)-L(w^\prime)=\lambda_0^{\max}(w^\prime,A)-\lambda_0^{\min}(w^\prime,A)$.
\item If $\delta(w^\prime)<\rho$, then we append $(xw^*y,\Cp(xw^*y),\pal(xwy))$ to the $\mathcal{F}$ list, otherwise we append it to $\cW_{n}$.
\item If $\delta(w^\prime)<\rho$, then we update $\Delta \gets \max\{\Delta,\delta(w^\prime)\}$. 
\item After completing the extensions of all words in $\cW_{n-1}$, if $\cW_{n} \neq \{\}$, we set $n \gets  n + 1$ and we go back to Step 3. Otherwise the algorithm has finished and we return $\mathcal{F}$ as our approximation data. Since the intervals decrease by $\delta(w^\prime)\leq(3-\sqrt{5})\cdot\delta(w)$ at each step, the algorithm eventually terminates.
\end{enumerate}

The resulting upper approximation of the Markov spectrum is obtained through
\begin{equation*}
    \bigcup_{(w,\Cp(w,A),\pal(w))\in\mathcal{F}} J(w,A)\cap[l,r].
\end{equation*}
The resulting covering contains $M\cap[l,r]$, and the maximum length of the intervals in this covering is bounded by $\Delta$. We can not, in general, guarantee that its Hausdorff distance from $M\cap[l,r]$ is at most \(\Delta\), since some of the intervals $J(w,A)$ may not intersect the Markov spectrum.  

\begin{remark}
    Since it is possible to prove that  a bi-infinite sequence $\underline{b}\in(\N^*)^\Z$ satisfies $m(\underline{b})=\lambda_0(\underline{b})$ with $b_0=1$ if and only if $\underline{b}=\overline{1}$, in the above algorithm is not necessary to test positions with digits equal to 1.
\end{remark}

\subsection{Example of the algorithm}

\subsubsection{Initial setup}

We want to find a covering of $M\cap[0,\sqrt{12}]$. Consider the alphabet $\cA=\{1,2\}$, so that $A=2$. The initial collection of words is
\begin{equation*}
    \cW_0=\{(1^*,\{0\},1),(2^*,\{0\},1)\}.
\end{equation*}
We also initialize $\mathcal F=\varnothing$ and $\Delta=0$.

\subsubsection{First iteration}

Since the initial words are palindromes, the set of extensions is
\begin{equation*}
    \mathcal E(1^*)=\mathcal E(2^*)
    =\{(1,1),(1,2),(2,2)\}.
\end{equation*}

We first consider the extensions of $(1^*,\{0\},1)$.

\begin{itemize}
\item Consider the trial word $w'=11^*1$. We have
\begin{align*}
    U(11^*1)&=\lambda_0^{\max}(11^*1,2)=[1;1,\overline{2,1}]+[0;1,\overline{2,1}]=2.4641\ldots,\\
    L(11^*1)&=\lambda_0^{\min}(11^*1,2)=[1;1,\overline{1,2}]+[0;1,\overline{1,2}]=2.1547\ldots.
\end{align*}

We now test the positions $-1$ and $1$. For $i=-1$, we obtain
\begin{align*}
    \lambda_{-1}^{\max}(11^*1,2)
    &=[1;\overline{1,2}]+[0;1,1,\overline{1,2}]
      =2.36\ldots,\\
    \lambda_{-1}^{\min}(11^*1,2)
    &=[1;\overline{2,1}]+[0;1,1,\overline{2,1}]
      =1.94\ldots.
\end{align*}
Thus,
\[
    \lambda_{-1}^{\min}(11^*1,2)<U(11^*1)
    \quad\text{and}\quad
    \lambda_{-1}^{\max}(11^*1,2)>L(11^*1).
\]
Consequently, $-1$ is a critical position. Since $11^*1$ is a palindrome, the same computation applies to the position $1$. Hence $\Cp(11^*1,2)=\{-1,0,1\}$. The word is allowable and, since it is a palindrome, we retain the tuple
\[
    (11^*1,\{-1,0,1\},1).
\]

\item Next, consider $w'=11^*2$. We have
\begin{equation*}
    U(11^*2)
    =[1;2,\overline{2,1}]+[0;1,\overline{2,1}]
    =2.098\ldots.
\end{equation*}
However,
\begin{equation*}
    \lambda_1^{\min}(11^*2,2)
    =[2;\overline{2,1}]+[0;1,1,\overline{2,1}]
    =2.943\ldots>U(11^*2).
\end{equation*}
Therefore, the word is not allowable and we discard this extension.

\item The extension $(2,2)$ gives the trial word $21^*2$. An analogous computation shows that it is not allowable either, so we discard it.
\end{itemize}

We now consider the extensions of $(2^*,\{0\},1)$.

\begin{itemize}
\item Consider the trial word $w'=12^*1$. We have
\begin{align*}
    U(12^*1)&=\lambda_0^{\max}(12^*1,2)=[2;1,\overline{2,1}]+[0;1,\overline{2,1}] =3.4641\dots,\\
    L(12^*1)&=\lambda_0^{\min}(12^*1,2)=[2;1,\overline{1,2}]+[0;1,\overline{1,2}]=3.1547\ldots.
\end{align*}

We test the positions $-1$ and $1$. Since $12^*1$ is a palindrome, the corresponding computations give the same bounds. In particular,
\begin{equation*}
    \lambda_{-1}^{\max}(12^*1,2)
    =\lambda_1^{\max}(12^*1,2)
    =[1;\overline{1,2}]+[0;2,1,\overline{1,2}]
    =2.120\ldots<L(12^*1).
\end{equation*}
Thus neither position is critical, and we obtain $\Cp(12^*1,2)=\{0\}$. The word is allowable, and we retain the tuple
\[
    (12^*1,\{0\},1).
\]

\item The extensions $(1,2)$ and $(2,2)$ give the trial words $12^*2$ and $22^*2$, respectively. Analogous computations show that both words are allowable and yield
\[
    \Cp(12^*2,2)=\{0,1\},
    \qquad
    \Cp(22^*2,2)=\{-1,0,1\}.
\]
Moreover, $12^*2$ is not a palindrome, whereas $22^*2$ is a palindrome. We therefore retain the tuples
\[
    (12^*2,\{0,1\},0)
    \quad\text{and}\quad
    (22^*2,\{-1,0,1\},1).
\]
\end{itemize}

Thus, the four allowable extensions $\cW_1$ obtained after the first iteration are
\begin{gather*}
    (11^*1,\{-1,0,1\},1), \quad (12^*1,\{0\},1), \quad (12^*2,\{0,1\},0), \quad (22^*2,\{-1,0,1\},1).
\end{gather*}
Their associated interval lengths are
\begin{align*}
    \delta(11^*1)&=\delta(12^*1)=0.3094\ldots, &\delta(12^*2)&=0.2113\ldots, & \delta(22^*2)&=0.1132\ldots.
\end{align*}
In particular, if $\rho>0.3094\ldots$, all four tuples are added to $\mathcal F$, and the first iteration terminates with $\Delta=0.3094\ldots$.

The resulting covering is
\begin{equation*}
    J(11^*1,2)\cup J(12^*1,2)
    \cup J(12^*2,2)\cup J(22^*2,2).
\end{equation*}

The maximum interval length is attained by the two words $11^*1$ and $12^*1$. 

This example illustrates a useful feature of the algorithm: intervals associated with words containing more occurrences of the digit $2$ can have substantially smaller lengths. Consequently, such intervals may reach the prescribed precision after fewer iterations and need not be subdivided as many times.

\subsubsection{Second iteration}

After another iteration, the surviving candidate list is
\begin{multline*}
\cW_2 =
\bigl\{
(111^*11,\{-2,-1,0,1,2\},1),
(112^*11,\{0\},1),
(112^*12,\{0,2\},0),\\
(212^*12,\{-2,0,2\},1),
(112^*21,\{0,1\},0),
(112^*22,\{0,2\},0),\\
(212^*21,\{-2,0,1\},0),
(212^*22,\{-2,0,2\},0),
(222^*22,\{-2,-1,0,1,2\},1)
\bigr\}.
\end{multline*}
Here, as before, each tuple records the word, its set of critical positions, and its palindrome indicator. The interval sizes are
\begin{align*}
    \delta(111^*11) &= 0.1132\dots & \delta(212^*12) &= 0.0582\dots &  \delta(212^*21) &= 0.0511\dots  \\
    \delta(112^*11) &= 0.1132\dots & \delta(112^*21) &= 0.0785\dots & \delta(212^*22) &= 0.0390\dots \\
    \delta(112^*12) &= 0.0857\dots & \delta(112^*22) &= 0.0665\dots &  \delta(222^*22) &= 0.0197\dots
\end{align*}

The largest interval in the resulting covering has diameter
\begin{equation*}
\Delta=\delta(112^*11)=\lambda_0^{\max}(112^*11,2)-\lambda_0^{\min}(112^*11,2)=0.1132\dots.
\end{equation*}

The same diameter is attained by $\delta(111^*11,2)$. Thus, after this iteration, the maximum interval length has decreased from $0.3094\ldots$ to $0.1132\ldots$.

At the other end of the range of interval sizes, the interval associated with $222^*22$ is the smallest among those in the new covering. Its diameter is, by definition,
\begin{equation*}
\delta(222^*22)
=\lambda_0^{\max}(222^*22,2)
-\lambda_0^{\min}(222^*22,2)=0.0197\dots.
\end{equation*}
The remaining interval diameters lie between these two extremes. The rest of the intervals sizes appear to distribute in a ``normal'' way from the largest interval size $0.1132\dots$ up to $0.0197\dots$. 

\subsubsection{Third iteration}

After another iteration, the surviving candidate list is

\begin{multline*}
\cW_3=
\bigl\{
(1111^*111,\{-3,-2,0,2,3\},1),(1112^*111,\{0\},1),(1112^*112,\{0,3\},0), \\
(2112^*112,\{-3,0,3\},0),(1112^*121,\{0,2\},0),(1112^*122,\{0\},0), \\
(2112^*121,\{-3,0,2\},0),(2112^*122,\{-3,0\},0),(1212^*121,\{-2,0,2\},1), \\
(1212^*122,\{-2,0\},0),(2212^*122,\{0\},1),(1112^*211,\{0\},0), \\
(2112^*211,\{-3,0\},0),(1112^*221,\{0,2\},0),(1112^*222,\{0,3\},0), \\
(2112^*221,\{-3,0,2\},0),(2112^*222,\{-3,0,3\},0),(2212^*211,\{-3,-2,0\},0), \\
(2212^*212,\{-3,-2,0,3\},0),(2212^*221,\{-2,0,2\},0),(2212^*222,\{-2,0,3\},0), \\
(2222^*222,\{-3,-2,-1,0,1,2,3\},1)
\bigr\}.
\end{multline*}

The largest interval (which is also attained by $\delta(1111^*111,2)$) in the resulting covering has diameter
\begin{equation*}
\Delta=\delta(1112^*111)=\lambda_0^{\max}(1112^*111,2)-\lambda_0^{\min}(1112^*111,2)=0.04394\dots.
\end{equation*}

\newpage

\begin{figure}[H]
\includegraphics[angle=90,width=0.25\textheight,keepaspectratio]{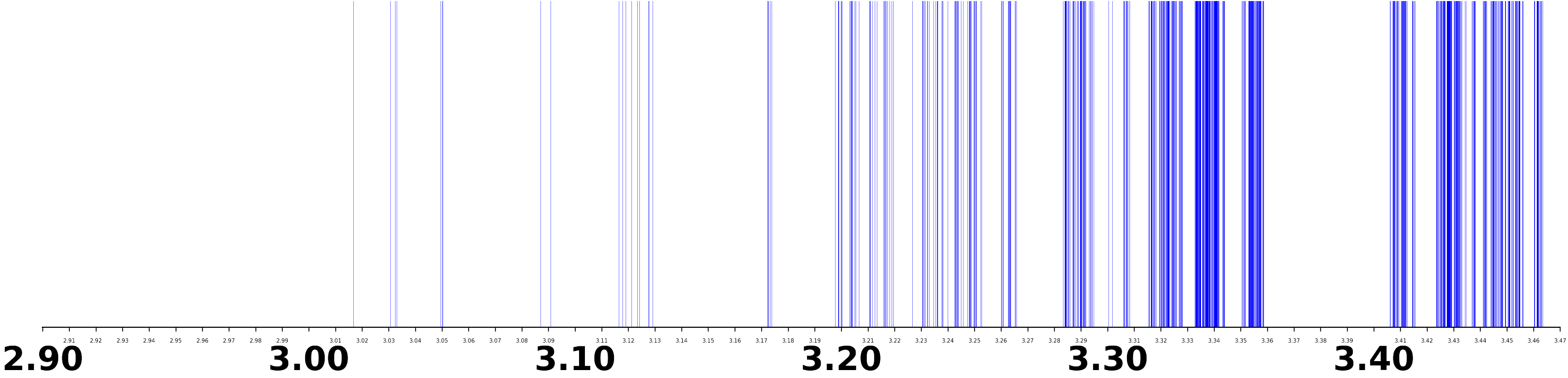}
\hspace*{2cm}
\includegraphics[angle=90,width=0.25\textheight,keepaspectratio]{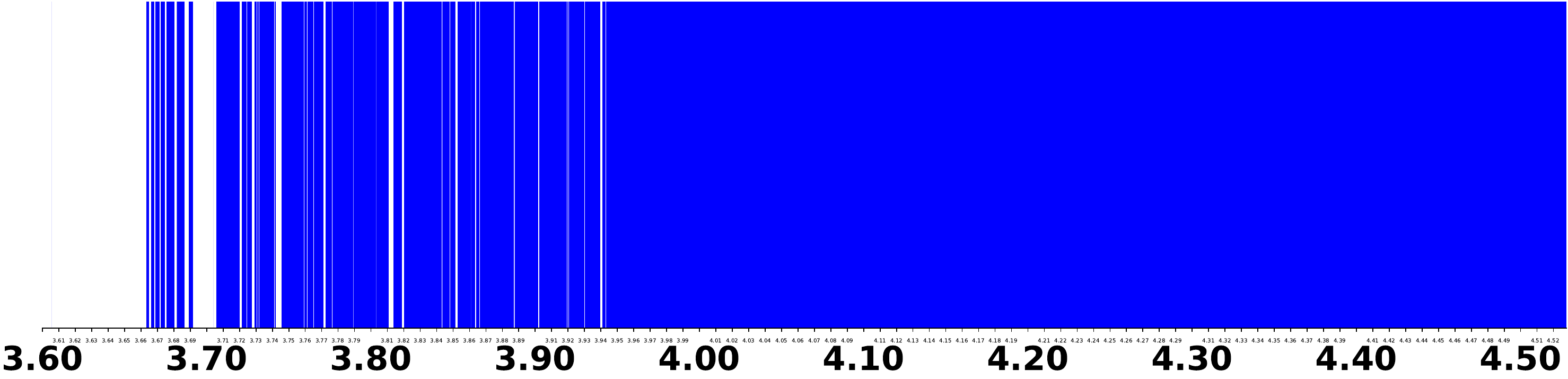}
\caption{The plot on the left is a covering of the Markov spectra between $[0,\sqrt{12}]$ with error at most $9.109424371\cdot10^{-8}$. The plot on the right is a covering of the Markov spectra between $[\sqrt{13},\nu_F]$ with error at most  $9.2\cdot 10^{-5}$.}
\end{figure}

\begin{figure}[H]
\includegraphics[angle=90,width=0.25\textheight,keepaspectratio]{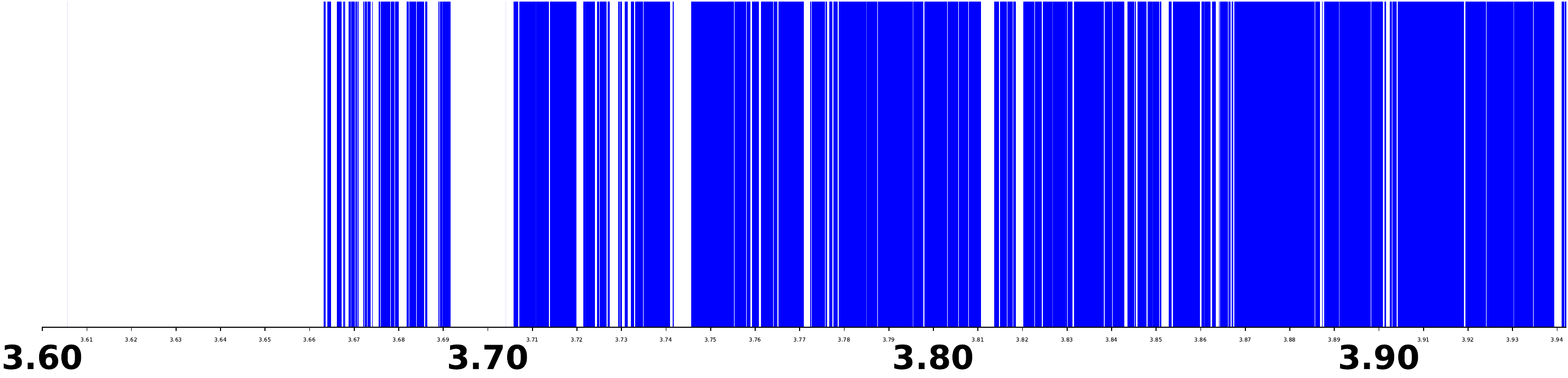}
\hspace*{2cm}
\includegraphics[angle=90,width=0.25\textheight,keepaspectratio]{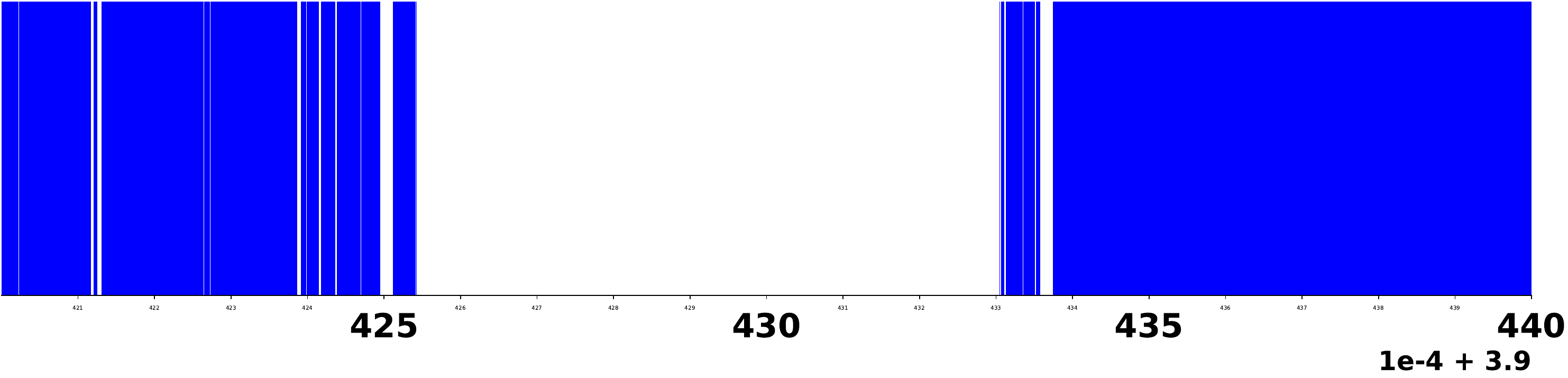}
\caption{Plots of coverings of the Markov spectra between $[3.6,3.942]$ and $[3.942,3.944]$ with error at most  $9.2\cdot 10^{-5}$.}
\end{figure}

\begin{figure}[H]
\includegraphics[angle=90,width=0.25\textheight,keepaspectratio]{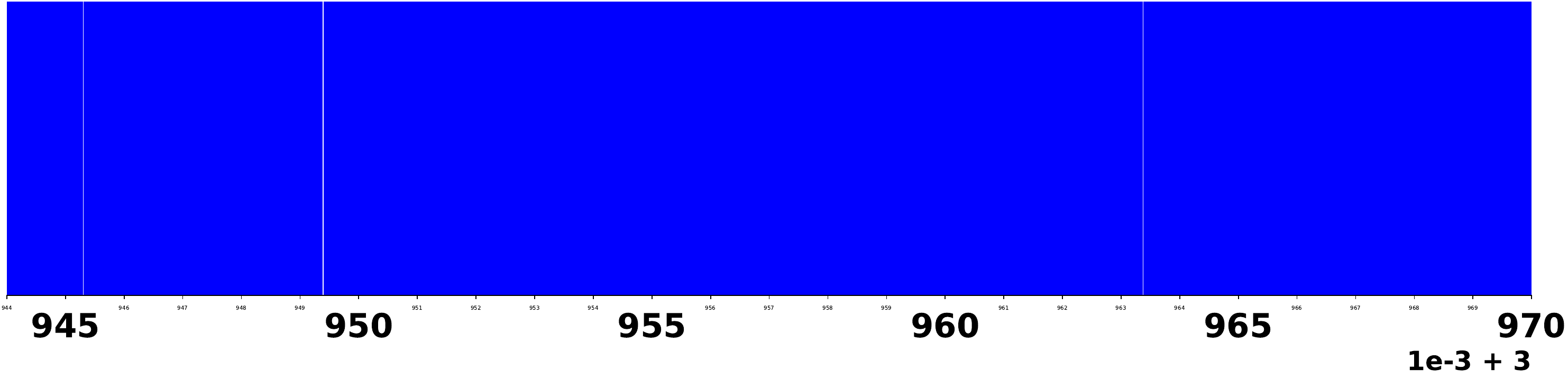}
\hspace*{2cm}
\includegraphics[angle=90,width=0.25\textheight,keepaspectratio]{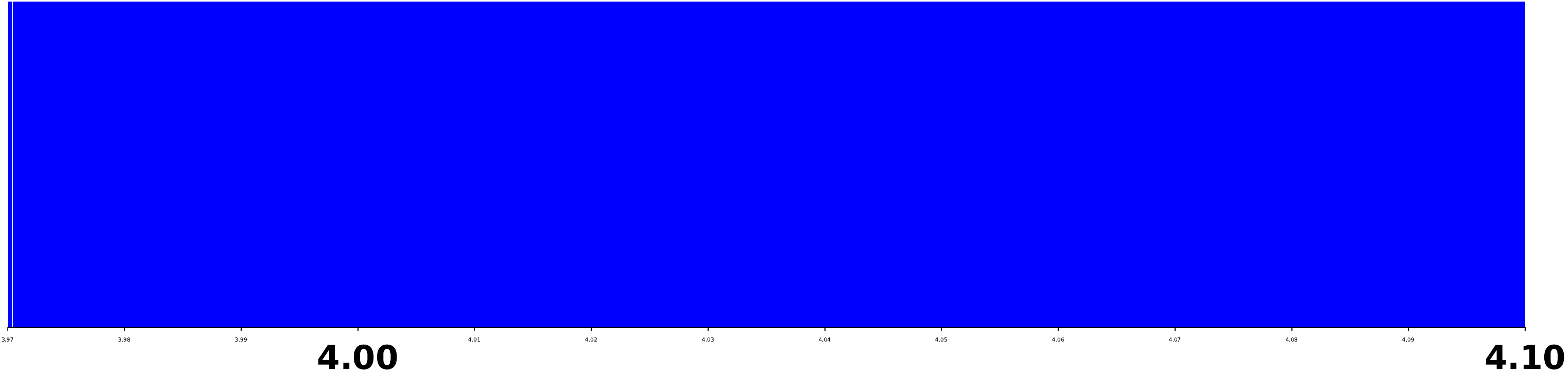}
\caption{Plots of coverings of the Markov spectra between $[3.944,3.97]$ and $[3.97,4.1]$ with error at most  $9.2\cdot 10^{-5}$.}
\end{figure}

\begin{figure}[H]
\includegraphics[angle=90,width=0.25\textheight,keepaspectratio]{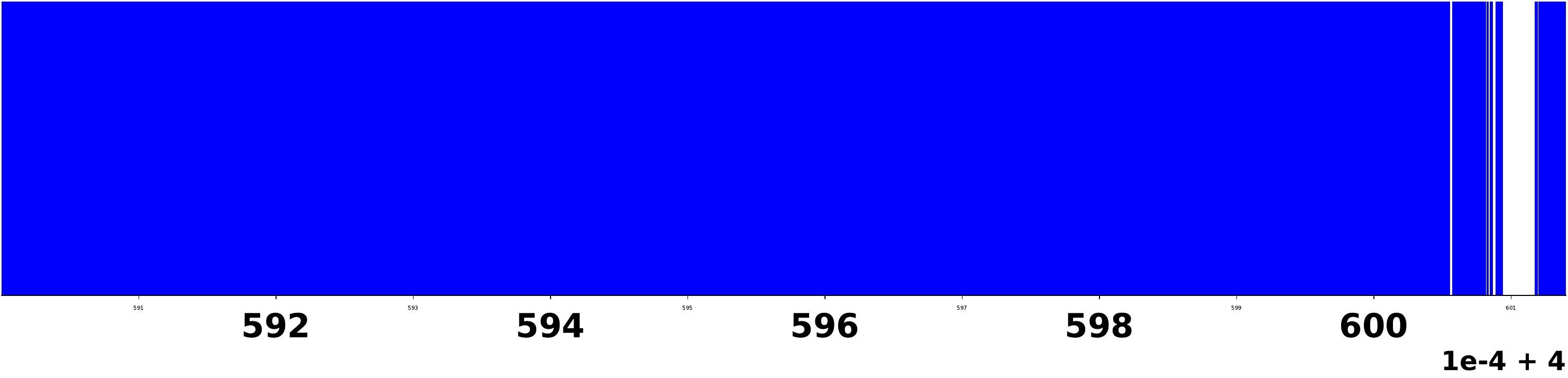}
\hspace*{2cm}
\includegraphics[angle=90,width=0.25\textheight,keepaspectratio]{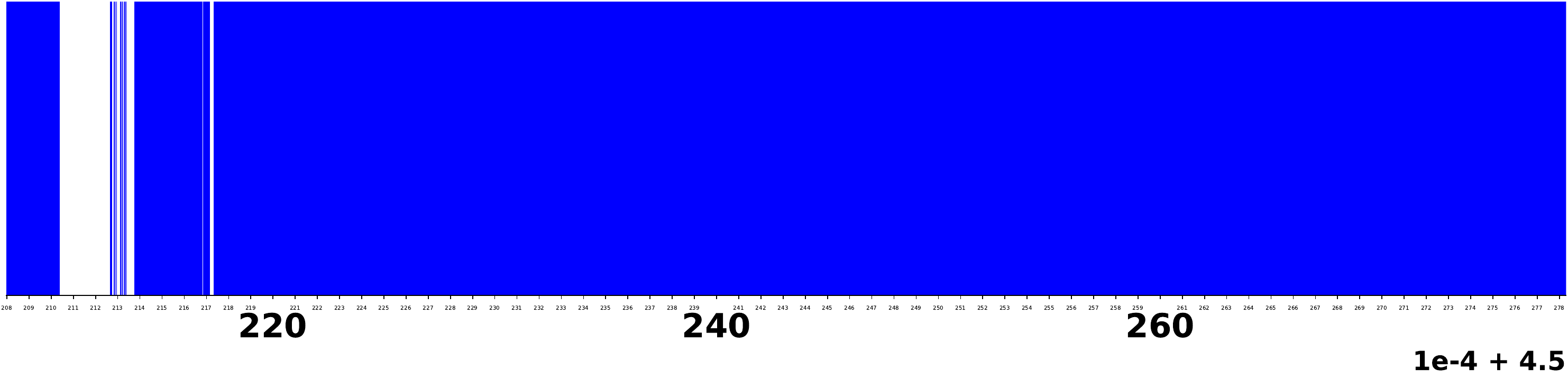}
\caption{Plots of coverings of the Markov spectra between $[4.509,4.601014]$ and $[4.52079673,4.527829566]$ with error at most  $9.2\cdot 10^{-5}$.}
\end{figure}

\begin{figure}[H]
\includegraphics[angle=90,width=0.25\textheight,keepaspectratio]{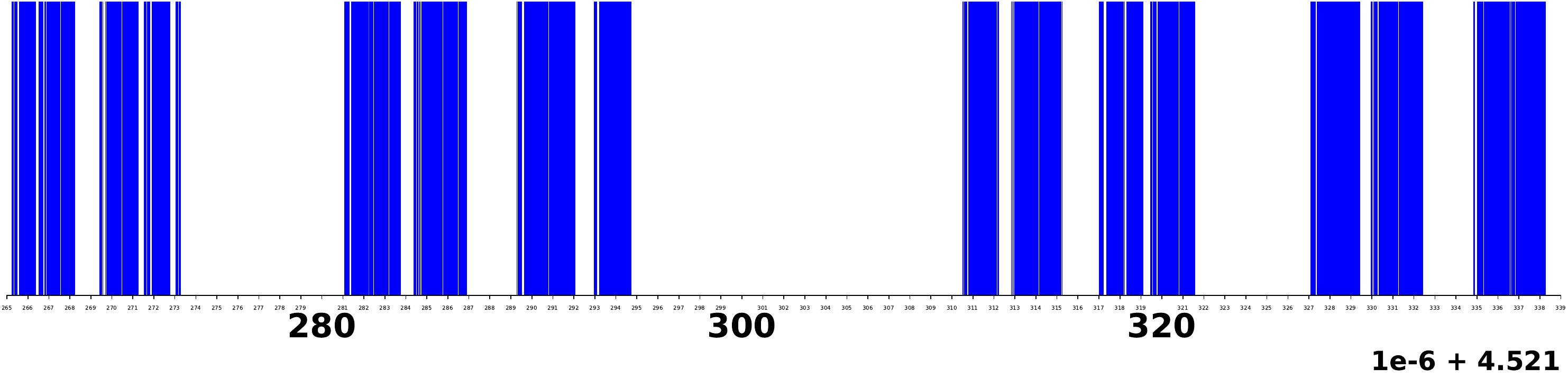}
\hspace*{2cm}
\includegraphics[angle=90,width=0.25\textheight,keepaspectratio]{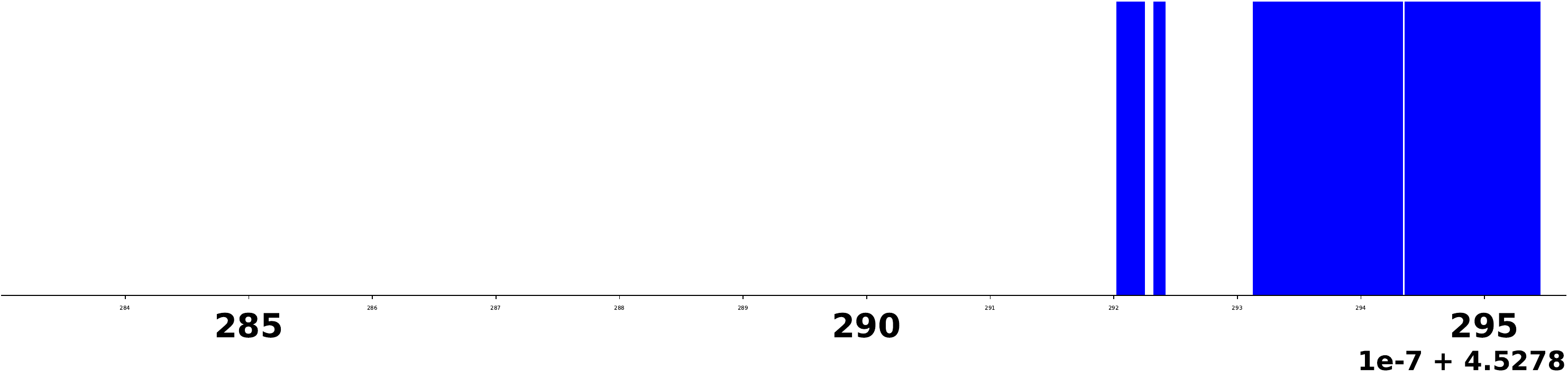}
\caption{The plot on the left is a covering of the Markov spectra between $[4.521265,4.521339]$ with error at most $3.96\cdot 10^{-8}$. The plot on the right is a covering of the Markov spectra between $[4.5278283,4.527829566]$ with error at most  $9.2\cdot 10^{-5}$.}
\end{figure}

\appendix

\section{Pseudo-code of the algorithm to cover the Markov spectrum}

Below is the pseudocode for the algorithm described in \Cref{sec:algorithm}, which determines the central portion $w$ of bi-infinite sequence $\underline{a}\in\{1,\dots,A\}^\Z$ for which there might exists $m(\underline{a})=\lambda_0(\underline{a})$ for some interval $[l,r]$ and for which their associated interval $J(w,A)$ has length at most some input error $\rho>0$. The code and the data produced with it that was used to produce the plots of the Markov spectrum is publicly available at \url{https://github.com/harolderazo/Lagrange_and_Markov_spectra}

\begin{algorithm}[!ht]\label{al:search}
	\caption{Find central blocks whose Markov values could lie in the range $[l,r]$.}
	\begin{algorithmic}
    \setlength{\baselineskip}{11pt}
		\Require{Alphabet $\{1,\dots,A\}$, Approximation error $\rho>0$, Interval $[l,r]\subset(0,\infty)$}
		\Let{$candidates$}{$[(1^{*},[0],True),(2^{*},[0],True),\dots,(A^{*},[0],True)]$}
		\Let{$all\_extensions$}{$(alphabet\times alphabet)$}
		\Let{$nonredundant\_extensions$}{$\{(x,y)\in(alphabet\times alphabet):x\leq y\}$}
		\Let{$word\_side\_len$}{1} \Comment{Length of each side of the central word}
		\Let{$\Delta$}{$0$}	\Comment{The current approximation error}
		\Let{$final\_blocks$}{$[]$}
		\While{$candidates \neq []$} 
		\Let{$new\_candidates$}{$[]$}
		\For{$(word,dan\_pos,is\_palindrome)\in candidates$}
		\Let{$extended\_dan\_pos$}{$[-word\_side\_len]+dan\_pos+[word\_side\_len]$}
		\If{$is\_palindrome$}
		\Let{$extensions$}{$nonredundant\_extensions$}
		\Else
		\Let{$extensions$}{$all\_extensions$}
		\EndIf
		\For{$(x,y)\in extensions$}
		\Let{$allowable$}{True}
		\Let{$trial\_word$}{$x+word+y$}
		\If{$\lambda_0^{\max}(trial\_word,A)<l$}
		\State \textbf{continue} \Comment{The trial word is below the range}
		\EndIf
		\If{$\lambda_0^{\min}(trial\_word,A)>r$}
		\State \textbf{continue} \Comment{The trial word is above the range}
		\EndIf
		\Let{$upper\_bound\_dan\_pos$}{$\lambda_0^{\max}(trial\_word,A)$}
		\Let{$lower\_bound\_dan\_pos$}{$\lambda_0^{\min}(trial\_word,A)$}
		\Let{$new\_dan\_pos$}{$[]$}
		\For{$i \in extended\_dan\_pos$}
        \If{$i \neq 0$ \textbf{and} $\lambda_i^{\min}(trial\_word,A)>upper\_bound\_dan\_pos$}
		\Let{$allowable$}{False}
		\State \textbf{break}
		\EndIf
		\If{$\lambda_{i}^{\max}(trial\_word,A)\geq lower\_bound\_dan\_pos$}
		\State append $i$ to $new\_dan\_pos$
		\EndIf
		\EndFor
		
		\algstore{bkbreak}
		
	\end{algorithmic}
\end{algorithm}

\setcounter{algorithm}{0}

\begin{algorithm}[!ht]
	\caption{Continued}
	\begin{algorithmic}
		\algrestore{bkbreak}

		\If{$allowable$}
		\If{$is\_palindrome$}
		\If{$x=y$}
		\Let{$is\_trial\_palindrome$}{True}
		\Else
		\Let{$is\_trial\_palindrome$}{False}
		\EndIf
        \Else
        \Let{$is\_trial\_palindrome$}{False}
		\EndIf

		\Let{$interval\_size$}{$\lambda_0^{\max}(trial\_word,A)-\lambda_0^{\min}(trial\_word,A)$}
		\If{$interval\_size<\rho$}
		\State append $(trial\_word,new\_dan\_pos,is\_trial\_palindrome)$ to $final\_blocks$\
        \Let{$\Delta$}{$\max\{\Delta,interval\_size\}$}
		\Else
		\State append $(trial\_word,new\_dan\_pos,is\_trial\_palindrome)$ to $new\_candidates$
		\EndIf
		\EndIf
		\EndFor
		\EndFor
		\Let{$candidates$}{$new\_candidates$}
		\If{$candidates \neq []$}
		\Let{$word\_side\_len$}{$word\_side\_len+1$}
		\EndIf
		\EndWhile
		\State \Return $final\_blocks$
\end{algorithmic}
\end{algorithm}

\printbibliography

@book{Fr75,
 Author = {Freiman, G. A.},
 Title = {Diophantine approximation and the geometry of numbers (The Markoff spectrum)},
 Year = {1975},
 Language = {(in Russian)},
 Publisher = {{Kalininski{\u{\i}}} {Gos}. {Univ}., Moscow},
 zbMATH = {3540889},
 Zbl = {0347.10025}
}

@article{Ha47,
 ISSN = {0003486X},
 URL = {http://www.jstor.org/stable/1969389},
 author = {Hall Jr., M.},
 journal = {Annals of Mathematics},
 number = {4},
 pages = {966--993},
 publisher = {Annals of Mathematics},
 title = {On the Sum and Product of Continued Fractions},
 volume = {48},
 year = {1947}
}

@article{BerConj,
title = {On the structure of the Markov spectrum},
journal = {Chapter III in Teorija cisel (Number Theory), Kalininskii Gosudarstvennyi
Universitet, Moscow},
language = {(in Russian)},
year = {1973},
pages = {50-78},
author = {A. A. Berstein}
}

@article{Fr68,
	author = {Freiman, G.  A. },
	da = {1968/02/01},
	doi = {10.1007/BF01094333},
	id = {Freiman1968},
	isbn = {1573-8876},
	journal = {Mathematical notes of the Academy of Sciences of the USSR},
	number = {2},
	pages = {125--128},
	title = {Noncoincidence of the Markov and Lagrange spectra},
	ty = {JOUR},
	url = {https://doi.org/10.1007/BF01094333},
	volume = {3},
	year = {1968}
}

@article {Newgaps,
    AUTHOR = {Jeffreys, Luke and Matheus, Carlos and Moreira, Carlos
              Gustavo},
     TITLE = {New gaps on the {L}agrange and {M}arkov spectra},
   JOURNAL = {J. Th\'eor. Nombres Bordeaux},
  FJOURNAL = {Journal de Th\'eorie des Nombres de Bordeaux},
    VOLUME = {36},
      YEAR = {2024},
    NUMBER = {1},
     PAGES = {311--338},
      ISSN = {1246-7405,2118-8572},
   MRCLASS = {11J06 (11A55)},
  MRNUMBER = {4788375},
}

@article {DynamicalCharacterization,
    AUTHOR = {Lima, Davi and Moreira, Carlos Gustavo},
     TITLE = {Dynamical characterization of initial segments of the {M}arkov
              and {L}agrange spectra},
   JOURNAL = {Monatsh. Math.},
  FJOURNAL = {Monatshefte f\"ur Mathematik},
    VOLUME = {199},
      YEAR = {2022},
    NUMBER = {4},
     PAGES = {817--852},
      ISSN = {0026-9255,1436-5081},
   MRCLASS = {11J06 (11A55 37B10)},
  MRNUMBER = {4497171},
MRREVIEWER = {Agamemnon\ Zafeiropoulos},
       DOI = {10.1007/s00605-022-01765-3},
       URL = {https://doi.org/10.1007/s00605-022-01765-3},
}

@article {Fractalgeometryofcomplement,
    AUTHOR = {Matheus, C. and Moreira, C. G.},
     TITLE = {Fractal geometry of the complement of {L}agrange spectrum in
              {M}arkov spectrum},
   JOURNAL = {Comment. Math. Helv.},
  FJOURNAL = {Commentarii Mathematici Helvetici. A Journal of the Swiss
              Mathematical Society},
    VOLUME = {95},
      YEAR = {2020},
    NUMBER = {3},
     PAGES = {593--633},
      ISSN = {0010-2571,1420-8946},
   MRCLASS = {11J06 (37B10 37C45 37C50)},
  MRNUMBER = {4152626},
MRREVIEWER = {Thomas\ Garrity},
       DOI = {10.4171/CMH/498},
       URL = {https://doi.org/10.4171/CMH/498}
}

@article {LGH,
    AUTHOR = {Erazo, Harold and Jeffreys, Luke and Moreira, Carlos Gustavo},
     TITLE = {On the classical {L}agrange and {M}arkov spectra: new results
              on the local dimension and the geometry of the difference set},
   JOURNAL = {Acta Arith.},
  FJOURNAL = {Acta Arithmetica},
    VOLUME = {220},
      YEAR = {2025},
    NUMBER = {1},
     PAGES = {29--100},
      ISSN = {0065-1036,1730-6264},
   MRCLASS = {11J06 (11A55 28A78 37B10)},
  MRNUMBER = {4940533},
       DOI = {10.4064/aa240821-26-5},
       URL = {https://doi.org/10.4064/aa240821-26-5},
}

@article {GuguThesis,
    AUTHOR = {de A. Moreira, Carlos Gustavo T.},
     TITLE = {Stable intersections of {C}antor sets and homoclinic
              bifurcations},
   JOURNAL = {Ann. Inst. H. Poincar\'e{} C Anal. Non Lin\'eaire},
  FJOURNAL = {Annales de l'Institut Henri Poincar\'e{} C. Analyse Non
              Lin\'eaire},
    VOLUME = {13},
      YEAR = {1996},
    NUMBER = {6},
     PAGES = {741--781},
      ISSN = {0294-1449,1873-1430},
   MRCLASS = {58F14 (58F15)},
  MRNUMBER = {1420497},
MRREVIEWER = {P.\ Mendes},
       DOI = {10.1016/S0294-1449(16)30122-6},
       URL = {https://doi.org/10.1016/S0294-1449(16)30122-6},
}

@article {M:formes_quadratiques1,
    AUTHOR = {Markoff, A.},
     TITLE = {Sur les formes quadratiques binaires ind\'{e}finies},
   JOURNAL = {Math. Ann.},
  FJOURNAL = {Mathematische Annalen},
    VOLUME = {15},
      YEAR = {1879},
    NUMBER = {3-4},
     PAGES = {381--406},
      ISSN = {0025-5831},
   MRCLASS = {11E16 (11A55)},
  MRNUMBER = {4788527},
	DOI = {10.1007/BF02086269},
	URL = {https://doi.org/10.1007/BF02086269},
}

@article {M:formes_quadratiques2,
	AUTHOR = {Markoff, A.},
	TITLE = {Sur les formes quadratiques binaires ind\'{e}finies II},
	JOURNAL = {Math. Ann.},
	FJOURNAL = {Mathematische Annalen},
	VOLUME = {17},
	YEAR = {1880},
	NUMBER = {3},
	PAGES = {379--399},
	ISSN = {0025-5831},
	MRCLASS = {DML},
	MRNUMBER = {1510073},
	DOI = {10.1007/BF01446234},
	URL = {https://doi.org/10.1007/BF01446234},
}

@article {geometricproperties,
    AUTHOR = {Moreira, C. G.},
     TITLE = {Geometric properties of the {M}arkov and {L}agrange spectra},
   JOURNAL = {Ann. of Math. (2)},
  FJOURNAL = {Annals of Mathematics. Second Series},
    VOLUME = {188},
      YEAR = {2018},
    NUMBER = {1},
     PAGES = {145--170},
      ISSN = {0003-486X,1939-8980},
   MRCLASS = {11J06 (11J70 28A78 37D20)},
  MRNUMBER = {3815461},
MRREVIEWER = {Thomas\ W.\ Cusick},
       DOI = {10.4007/annals.2018.188.1.3},
       URL = {https://doi.org/10.4007/annals.2018.188.1.3}
}

@article {MendesOliveira,
    AUTHOR = {Mendes, Pedro and Oliveira, Fernando},
     TITLE = {On the topological structure of the arithmetic sum of two
              {C}antor sets},
   JOURNAL = {Nonlinearity},
  FJOURNAL = {Nonlinearity},
    VOLUME = {7},
      YEAR = {1994},
    NUMBER = {2},
     PAGES = {329--343},
      ISSN = {0951-7715,1361-6544},
   MRCLASS = {58F14 (28A80 34C37 58F15)},
  MRNUMBER = {1267692},
MRREVIEWER = {Feliks\ Przytycki},
       DOI = {10.1088/0951-7715/7/2/002},
       URL = {https://doi.org/10.1088/0951-7715/7/2/002},
}

@article {MMPV,
    AUTHOR = {Matheus, C. and Moreira, C. G. and Pollicott,
              M. and Vytnova, P.},
     TITLE = {Hausdorff dimension of {G}auss-{C}antor sets and two
              applications to classical {L}agrange and {M}arkov spectra},
   JOURNAL = {Adv. Math.},
  FJOURNAL = {Advances in Mathematics},
    VOLUME = {409},
      YEAR = {2022},
    NUMBER = {part B},
     PAGES = {Paper No. 108693, 70},
      ISSN = {0001-8708},
   MRCLASS = {11K50 (11A55 11J06 37D35 37M25)},
  MRNUMBER = {4483249},
MRREVIEWER = {Thomas Ward},
       DOI = {10.1016/j.aim.2022.108693},
       URL = {https://doi.org/10.1016/j.aim.2022.108693}
}

@book {Cusick-Flahive,
    AUTHOR = {Cusick, T. W. and Flahive, M. E.},
     TITLE = {The {M}arkoff and {L}agrange spectra},
    SERIES = {Mathematical Surveys and Monographs},
    VOLUME = {30},
 PUBLISHER = {American Mathematical Society, Providence, RI},
      YEAR = {1989},
     PAGES = {x+97},
      ISBN = {0-8218-1531-8},
   MRCLASS = {11J06},
  MRNUMBER = {1010419},
MRREVIEWER = {Richard T. Bumby},
       DOI = {10.1090/surv/030},
       URL = {https://doi.org/10.1090/surv/030}
}

@article {VincentAlgorithm,
    AUTHOR = {Delecroix, Vincent and Matheus, Carlos and Moreira, Carlos
              Gustavo},
     TITLE = {Approximations of the {L}agrange and {M}arkov spectra},
   JOURNAL = {Math. Comp.},
  FJOURNAL = {Mathematics of Computation},
    VOLUME = {89},
      YEAR = {2020},
    NUMBER = {325},
     PAGES = {2521--2536},
      ISSN = {0025-5718,1088-6842},
   MRCLASS = {11J06 (05A15 37-04 37B10)},
  MRNUMBER = {4109576},
MRREVIEWER = {Hemangi\ Madhusudan\ Shah},
       DOI = {10.1090/mcom/3513},
       URL = {https://doi.org/10.1090/mcom/3513},
}

@incollection {Sullivan,
    AUTHOR = {Sullivan, Dennis},
     TITLE = {Differentiable structures on fractal like sets, determined by
              intrinsic scaling functions on dual {C}antor sets},
 BOOKTITLE = {Nonlinear evolution and chaotic phenomena ({N}oto, 1987)},
    SERIES = {NATO Adv. Sci. Inst. Ser. B: Phys.},
    VOLUME = {176},
     PAGES = {101--110},
 PUBLISHER = {Plenum, New York},
      YEAR = {1988},
      ISBN = {0-306-42909-8},
   MRCLASS = {58F12 (28A80)},
  MRNUMBER = {1109972},
       DOI = {10.1007/978-1-4613-1017-4\_7},
       URL = {https://doi.org/10.1007/978-1-4613-1017-4_7},
}

@article {Gazette,
    AUTHOR = {Matheus, C.},
     TITLE = {Renormalisation: une puissante technique dans la th\'eorie des
              syst\`emes dynamiques},
   JOURNAL = {Gaz. Math.},
  FJOURNAL = {Gazette des Math\'ematiciens},
    NUMBER = {180},
      YEAR = {2024},
     PAGES = {22--35},
      ISSN = {0224-8999,2275-0622},
   MRCLASS = {37E20},
  MRNUMBER = {4745400},
MRREVIEWER = {Dieter\ H.\ Mayer},
}

@article {PalisTakensAnnals1987,
    AUTHOR = {Palis, J. and Takens, F.},
     TITLE = {Hyperbolicity and the creation of homoclinic orbits},
   JOURNAL = {Ann. of Math. (2)},
  FJOURNAL = {Annals of Mathematics. Second Series},
    VOLUME = {125},
      YEAR = {1987},
    NUMBER = {2},
     PAGES = {337--374},
      ISSN = {0003-486X,1939-8980},
   MRCLASS = {58F10 (58F15)},
  MRNUMBER = {881272},
MRREVIEWER = {Genesio\ dos Reis},
       DOI = {10.2307/1971313},
       URL = {https://doi.org/10.2307/1971313},
}

@article {NewhouseIHES1979,
    AUTHOR = {Newhouse, Sheldon E.},
     TITLE = {The abundance of wild hyperbolic sets and nonsmooth stable
              sets for diffeomorphisms},
   JOURNAL = {Inst. Hautes \'Etudes Sci. Publ. Math.},
  FJOURNAL = {Institut des Hautes \'Etudes Scientifiques. Publications
              Math\'ematiques},
    NUMBER = {50},
      YEAR = {1979},
     PAGES = {101--151},
      ISSN = {0073-8301,1618-1913},
   MRCLASS = {58F13},
  MRNUMBER = {556584},
MRREVIEWER = {Zbigniew\ Nitecki},
       URL = {http://www.numdam.org/item?id=PMIHES_1979__50__101_0},
}

@article {MoreiraYoccoz2010,
    AUTHOR = {Moreira, Carlos Gustavo and Yoccoz, Jean-Christophe},
     TITLE = {Tangences homoclines stables pour des ensembles hyperboliques
              de grande dimension fractale},
   JOURNAL = {Ann. Sci. \'Ec. Norm. Sup\'er. (4)},
  FJOURNAL = {Annales Scientifiques de l'\'Ecole Normale Sup\'erieure.
              Quatri\`eme S\'erie},
    VOLUME = {43},
      YEAR = {2010},
    NUMBER = {1},
     PAGES = {1--68},
      ISSN = {0012-9593,1873-2151},
   MRCLASS = {37D05 (37C45 37E30 37G25)},
  MRNUMBER = {2583264},
MRREVIEWER = {Isabel\ Lug\~ao\ Rios},
       DOI = {10.24033/asens.2115},
       URL = {https://doi.org/10.24033/asens.2115},
}

\end{document}